\documentclass [11pt,a4paper]{article}
\usepackage[ansinew]{inputenc}
\usepackage{amsthm,amsmath,amssymb}
\usepackage{graphicx}
\usepackage{cite}
\usepackage{mathrsfs}
\usepackage{bm}

\usepackage{color}

\newtheorem{theorem}{\textbf{Theorem}}[section]
\newtheorem{lemma}{\textbf{Lemma}}[section]
\newtheorem{proposition}{\textbf{Proposition}}[section]
\newtheorem{corollary}{\textbf{Corollary}}[section]
\newtheorem{remark}{\textbf{Remark}}[section]
\newtheorem{definition}{\textbf{Definition}}[section]

\allowdisplaybreaks[4]

\def\be{\begin{equation}}
\def\ee{\end{equation}}
\def\bea{\begin{eqnarray}}
\def\eea{\end{eqnarray}}
\def\bt{\begin{theorem}}
\def\et{\end{theorem}}
\def\bl{\begin{lemma}}
\def\el{\end{lemma}}
\def\br{\begin{remark}}
\def\er{\end{remark}}
\def\bp{\begin{proposition}}
\def\ep{\end{proposition}}
\def\bc{\begin{corollary}}
\def\ec{\end{corollary}}
\def\bd{\begin{definition}}
\def\ed{\end{definition}}

\def\p{\partial}

\def\non{\nonumber}
\def\wha{\widehat}
\def\wt{\widetilde}

\def\bfv{\mathbf{v}}
\def\bfd{\mathbf{d}}
\def\bphi{\bm{\phi}}
\def\bpsi{\bm{\psi}}
\def\bom{\bm{\omega}}
\def\bxi{\bm{\xi}}
\def\bfe{\mathbf{e}}
\def\bfw{\mathbf{w}}

\definecolor{coloras}{rgb}{0.,0.67,0}
\begin{document}

\title{Optimal boundary control of a simplified Ericksen--Leslie
system for nematic liquid crystal flows in $2D$}

\author{
  Cecilia Cavaterra
  \footnote{ Dipartimento di Matematica,
Universit\`a degli Studi di Milano, Via Saldini 50, 20133 Milano,
Italy. \texttt{cecilia.cavaterra@unimi.it}
 }
 \and
  Elisabetta Rocca
  \footnote{Dipartimento di Matematica,
  Universit\`a degli Studi di Pavia, Via Ferrata 5, 27100, Pavia, Italy.
  \texttt{elisabetta.rocca@unipv.it}
  }
  \and
  Hao Wu
  \footnote{School of Mathematical Sciences and Shanghai
Key Laboratory for Contemporary Applied Mathematics, Fudan
University, 200433 Shanghai, China.
    \texttt{haowufd@yahoo.com}}
}

\date{\today}

\maketitle


\begin{abstract}

In this paper, we investigate an optimal boundary control problem for a two dimensional simplified Ericksen--Leslie
system modelling the incompressible nematic liquid crystal flows.
The hydrodynamic system consists of the Navier--Stokes equations for the fluid velocity coupled with a convective Ginzburg--Landau
type equation for the averaged molecular orientation.
The fluid velocity is assumed to satisfy a no-slip boundary condition,
while the molecular orientation is subject to a time-dependent Dirichlet boundary condition that corresponds to the strong
anchoring condition for liquid crystals. We first establish the existence of optimal boundary controls.
Then we show that the control-to-state operator is Fr\'echet differentiable between appropriate Banach spaces
and derive first-order necessary optimality conditions in terms of a variational
inequality involving the adjoint state variables.
 \medskip

\noindent \textbf{Keywords}: Optimal boundary control, first-order necessary optimality conditions, nematic liquid crystal flow, Navier--Stokes equations.\\
\textbf{AMS Subject Classification}: 49J20, 35Q35, 76A15, 76D05.
\end{abstract}

\section{Introduction}
\label{Intro}
We consider the following hydrodynamic system for
incompressible liquid crystal flows of nematic type:
 \bea
 \p_t\mathbf{v}+\mathbf{v}\cdot\nabla \mathbf{v}-\nu \Delta \mathbf{v}
 +\nabla P &=&-\lambda
 \nabla\cdot(\nabla \mathbf{d}\odot\nabla \mathbf{d}),\label{1}\\
 \nabla \cdot \mathbf{v} &=& 0,\label{2}\\
 \p_t\mathbf{d} +\mathbf{v}\cdot\nabla \mathbf{d}&=&\eta(\Delta \mathbf{d}-\mathbf{f}(\mathbf{d})),\label{3}
 \eea
in $\Omega \times\mathbb{R}^+$, where $\Omega \subset \mathbb{R}^n$ $(n=2,3)$
is a bounded domain with smooth boundary $\Gamma$. In the system \eqref{1}--\eqref{3},
$\mathbf{v}$ stands for the velocity field of the fluid, $\mathbf{d}$ represents the averaged
macroscopic/continuum molecular orientation and $P$ is a scalar function representing the pressure
(including both the hydrostatic and the induced elastic part from
the orientation field). The positive constants $\nu, \lambda$ and $\eta$ stand
for the fluid viscosity, the competition between kinetic energy and elastic potential
energy, and the elastic relaxation time (Deborah number) for
the molecular orientation field. The $n\times n$ matrix $\nabla \mathbf{d}\odot \nabla
\mathbf{d}$ denotes the Ericksen stress tensor whose $(i,j)$-th entry is
given by $\nabla_i \mathbf{d}\cdot \nabla_j \mathbf{d}$, for $1\leq
i,j\leq n$. The vector valued nonlinear function
$\mathbf{f}(\mathbf{d})$ is the gradient of certain smooth scalar potential function $F(\mathbf{d}):\mathbb{R}^n\rightarrow \mathbb{R}$ such that
$\mathbf{f}(\mathbf{d})=\nabla_{\mathbf{d}} F(\mathbf{d})$. In this paper, we take $F$ to be the Ginzburg--Landau approximation, i.e.,
$$F(\bfd)=\frac{1}{4\epsilon^2}(|\bfd|^2-1)^2,\quad \mathbf{f}(\mathbf{d})=\frac{1}{\epsilon^2}(|\mathbf{d}|^2-1)\mathbf{d},\quad \epsilon>0,$$
which has been used to relax the nonlinear constraint $|\mathbf{d}|=1$ on the molecule length in the literature (cf.
\cite{LL95,LLW}). Besides, we assume that the system
\eqref{1}--\eqref{3} is subject to the following boundary and initial conditions
 \begin{align}
 & \mathbf{v}(x,t)=\mathbf{0},\quad \mathbf{d}(x,t)=\mathbf{h}(x,t),\qquad (x, t)\in \Gamma\times
 \mathbb{R}^+,
 \label{4} \\
 &\mathbf{v}|_{t=0}=\mathbf{v}_0(x) \ \ \text{with}\ \nabla\cdot \mathbf{v}_0=0,\quad
 \mathbf{d}|_{t=0}=\mathbf{d}_0(x),\qquad x\in \Omega.\label{5}
 \end{align}

System \eqref{1}--\eqref{3} was firstly proposed in \cite{lin1} as a
simplified approximate system of the original Ericksen--Leslie model
for nematic liquid crystal flows (cf. \cite{E1,Le}).
Well-posedness of the autonomous initial boundary value problem of system
\eqref{1}--\eqref{3} (namely, with the no-slip
boundary condition for $\mathbf{v}$ and a time-independent Dirichlet
boundary condition for $\mathbf{d}$) was first analyzed in
\cite{LL95} (see also \cite{LL96} for the result on partial regularity).
Concerning the long-time behavior of global solutions to the
autonomous system, a natural question
on the uniqueness of asymptotic limit as $t\to +\infty$ was raised in \cite{LL95}.
This question was answered later in \cite{W10},
where it was proven that each trajectory converges
to a single steady state by using the \L ojasiewicz--Simon approach
(see \cite{PRS} for some generalizations).
We also refer to \cite{DMM1} for the asymptotic behavior of system
\eqref{1}--\eqref{3} in the whole space $\mathbb{R}^3$ and to \cite{DMM2,HW13} for the inhomogeneous
case with non-constant density. Some generalizations of system \eqref{1}--\eqref{3} have been considered in \cite{CR,CRW,sunliu,WXL13},
where the stretching effects are taken into account.

The technically more challenging case of time-dependent
Dirichlet boundary conditions for $\mathbf{d}$ has been recently
analyzed in \cite{B,C06,C09,G09, GW13}. Under proper regularity
assumptions of the time-dependent boundary datum $\mathbf{h}(x,t)$, existence of global weak solutions,
existence of global regular solutions under large viscosity, and weak/strong uniqueness were obtained in
\cite{C09}. Regularity criteria for local strong solutions in the three dimensional case can be
found in \cite{G09}. Concerning the long-time behavior, existence of global and
exponential attractors was proven in \cite{B} when the spatial dimension is two,
allowing the presence of a time-dependent external force. Besides,
in \cite{GW13}, stability of local energy minimizers and convergence to
a single equilibrium for any bounded trajectory were established.

With the well-posedness results of \cite{B,GW13} at hand, the road is paved
for studying optimal control problems associated with the system \eqref{1}--\eqref{5}
at least when the spatial dimension is two. This is the goal of this paper.
We note that in this case, the velocity field
$\mathbf{v}=(v_1,v_2)^{\mathrm{tr}}: \Omega\to \mathbb{R}^2$ is reduced to a two dimensional vector, while
the molecular director $\bfd=(d_1,..., d_n)^{\mathrm{tr}}: \Omega\to \mathbb{R}^n$ ($n=2,3$)
is allowed to be either two or three dimensional.

In particular, we are interested in the optimal boundary control problem for system \eqref{1}--\eqref{5}.
From the physical point of view, nonhomogeneous Dirichlet boundary condition
for the director $\mathbf{d}$ corresponds to the strong anchoring condition for liquid crystals.
The so-called anchoring refers to the description of how the molecular director is aligned on the boundary surfaces,
which is an important issue both in theoretic studies and applications of liquid crystals \cite{Ste}.
One typical type of those anchoring conditions is the strong anchoring that
occurs when the surface energy is sufficiently large. In this case, the
molecular orientation on the boundary can be simply fixed in
a preferred orientation determined by proper alignment techniques (see \cite{ZCW}).
Strong anchoring conditions are widely used due to their simplicity.
On the other hand, the well-posedness results in \cite{B,C09,GW13} imply that the dynamics of liquid crystal flow in the whole domain $\Omega$
can be determined by its boundary conditions.
This motives us to study the optimal boundary control problem for system \eqref{1}--\eqref{5}.

Throughout this paper, we assume that $T\in (0,+\infty)$ is a given finite final time and we set for
convenience
 \be
 Q:=\Omega \times (0,T),\quad \Sigma:= \Gamma \times (0,T).\nonumber
 \ee
Moreover, we make the following basic assumptions:
\begin{itemize}
\item[(A1)] $\beta_i\geq 0$ ($i=1,2,3,4$) and $\gamma\geq 0$ are given constants that do not vanish simultaneously.
\item[(A2)] The vector-valued functions
$$  \mathbf{v}_Q\in L^2(0,T; \mathbf{H}),\quad \mathbf{d}_Q \in \mathbf{L}^2(Q),\quad \mathbf{v}_\Omega\in \mathbf{H},
\quad \mathbf{d}_\Omega \in \mathbf{L}^2(\Omega)$$
are given target functions where the space $\mathbf{H}$ is given as in \eqref{spaces} below.
\end{itemize}
Then the optimal boundary control problem under investigation reads as follows:
\begin{itemize}
\item[(\textbf{CP})] \emph{Minimize the tracking type cost functional:}
\bea
\mathcal{J}((\mathbf{v}, \mathbf{d}), \mathbf{h}) &:=&
\frac{\beta_1}{2}\|\mathbf{v}-\mathbf{v}_Q\|^2_{\mathbf{L}^2(Q)}+\frac{\beta_2}{2}\|\mathbf{d}-\mathbf{d}_Q\|^2_{\mathbf{L}^2(Q)}
+\frac{\beta_3}{2}\|\mathbf{v}(T)-\mathbf{v}_\Omega\|^2_{\mathbf{L}^2(\Omega)}\nonumber\\
&& + \frac{\beta_4}{2}\|\mathbf{d}(T)-\mathbf{d}_\Omega\|^2_{\mathbf{L}^2(\Omega)}+\frac{\gamma}{2}\|\mathbf{h}\|^2_{\mathbf{L}^2(\Sigma)},\label{costJ}
\eea
\emph{subject to the boundary control constraint} $\mathbf{h}$ \emph{as well as the state constraint due to the initial boundary value problem} \eqref{1}--\eqref{5}.
\end{itemize}
Here, the vector $\mathbf{h}$ plays the role of a boundary control,
which is postulated to belong to a suitable closed, bounded and convex set $\widetilde{\mathcal{U}}^M_{\mathrm{ad}}$ in the space of controls $\widetilde{\mathcal{U}}$
(which will be specified later, see \eqref{U}, \eqref{Uad}). Besides, in the cost functional \eqref{costJ}, the pair $(\mathbf{v}, \mathbf{d})$ is the
unique global strong solution to the state problem \eqref{1}--\eqref{5} subject to the time-dependent Dirichlet boundary condition
$\mathbf{d}|_{\Gamma}=\mathbf{h}(x,t)$.

The main results of this paper are summarized as follows:
\begin{itemize}
\item[(1)] we establish the existence of optimal boundary controls for  problem (\textbf{CP}) (see Theorem \ref{existence});
\item[(2)] we show that the control-to-state operator $\mathcal{S}$ defined by problem \eqref{1}--\eqref{5}
is Fr\'echet differentiable between appropriate Banach spaces (see Theorem \ref{SD});
\item[(3)] we derive the first-order necessary optimality condition (see Theorem \ref{nece1}) and in particular, in terms of a variational
inequality involving the adjoint states (see Theorem \ref{nece2}).
\end{itemize}

Before ending the introduction, we would like to make some comments on the results of this paper and related literature.

Despite its physical and mathematical interests, to the best of our knowledge,
the optimal control problem (\textbf{CP}) for nematic liquid crystal flows has never been tackled in the literature.
We refer to \cite{FS92,FS94,Fu82,FGH98,FGH05} and the references therein for extensive studies on various optimal control
problems of the time-dependent Navier--Stokes equations for single viscous Newtonian fluids.
On the other hand, optimal control problems related to complex fluids have been studied in the literature (for instance, two-phase flows),
but never in the framework of liquid crystal flows.
For example, concerning the coupled Navier--Stokes--Cahn--Hilliard (or Allen--Cahn) system for two-phase fluids,
there exist recent contributions on optimal control problems for the time-discretized local version of the system (see \cite{HW2,HW3})
and on numerical aspects of the control problem (see \cite{HK}).
It seems that a rigorous analysis for the problem without time discretization has never been performed before until recently
 a diffuse interface model for incompressible isothermal mixtures of two immiscible fluids coupling the Navier--Stokes
system with a convective nonlocal Cahn--Hilliard equation in two spatial dimensions has been analyzed in \cite{FRS} from the optimal
control point of view. There the control was distributed (i.e., in terms of a body force in the bulk) and was not located on the boundary.
Even for the much simpler case of the convective Cahn--Hilliard equation, where
the velocity is prescribed so that the Navier--Stokes equation is not present, only very few contributions exist that deal with optimal control problems.
In this direction, we refer to \cite{ZL1,ZL2} for local models in one and two spatial
dimensions and to the recent paper \cite{RS}, in which the first-order necessary optimality conditions
were derived for the nonlocal convective Cahn--Hilliard system in three dimensions, in the case of degenerate mobilities and singular potentials.
At last, regarding the Allen--Cahn type equations (i.e., the scalar version of the director equation \eqref{3} with zero velocity), distributed
and boundary optimal control problems with various types of dynamic boundary conditions have been studied in a number of recent papers, in particular
for the case of double obstacle potentials (see \cite{Hassan1,Hassan2,HKR,CS}).

Our control problem (\textbf{CP}) deals with boundary controls of Dirichlet type,
which are important in many practical applications such as the active boundary control of single Newtonian fluids
(see, e.g., \cite{FGH98,FGH05,HK04}).
For instance, in \cite{FGH98} the authors studied optimal boundary control problems for the two-dimensional time-dependent
Navier--Stokes equations in an unbounded domain, where the control is effected through the Dirichlet boundary condition for
the fluid velocity.
They established the existence of optimal boundary controls in a suitable subset of
the trace space for velocity fields with almost minimal regularity. Moreover, they derived the
optimality system from which the optimal states and boundary controls can be determined.
We note that, in the practice, boundary controls with low regularity are excepted to be admissible
since one may be interested in blowing and suction as controls on part of the boundary, which
could possibly allow jumps and satisfy point-wise bounds (see \cite{GHZ}).
Thus, one usually works with very weak solutions of the evolutionary system, see e.g.,
\cite{AR02, KV07} for examples of Dirichlet boundary control problems for parabolic equations.
However, the situation is different when the Navier--Stokes equations are concerned.
It was pointed out in \cite{FGH98} that an important feature of the Dirichlet boundary control problem is that one
can derive an optimality system only in spaces of sufficiently smooth functions for
which the nonlinear terms of the Navier--Stokes system are subordinate to the linear
terms. For this reason, in \cite{FGH98}, the authors worked with sufficiently regular spaces of
boundary data that allow to obtain finite energy weak solutions for the Navier--Stokes
equations. Moreover, in contrast to the classical parabolic boundary control problems \cite{AR02, KV07}, it is necessary to
fulfill certain compatibility conditions for the boundary and initial values in this case (see \cite{FGH98} for further details).

For our current problem (\textbf{CP}) the situation is even more complicated.
The state system \eqref{1}--\eqref{5} consists of the Navier--Stokes equations coupled with a convective Ginzburg--Landau type equation,
which involves highly nonlinear multi-scale interactions between the macroscopic fluid velocity $\bfv$ and the microscopic molecular director $\bfd$.
Different from the simple fluid case as in \cite{FGH98, FGH05}, the boundary control $\mathbf{h}$ is now imposed on the director $\bfd$,
which influences the dynamics of the fluid through the higher-order nonlinear Ericksen stress tensor (i.e., in terms of a nonlinear bulk force).
We observe that the existence of an optimal boundary control to problem (\textbf{CP}) could be proven in less regular spaces
(see e.g., \eqref{hyp2}--\eqref{hyp4} that yield finite energy weak solutions
of system \eqref{1}--\eqref{5}). Nevertheless, in order to show the Fr\'echet differentiability of the control-to-state operator $\mathcal{S}$, as well as the
first-order necessary optimality conditions, we have to work with more regular trace spaces with a compatibility condition between the boundary value of the initial datum $\mathbf{d}_0$ and the initial value of the boundary control $\mathbf{h}$, which ensure the existence of a unique global strong solution to the state problem
\eqref{1}--\eqref{5} (see assumptions \eqref{hyp2s}--\eqref{hyp3s} with \eqref{hyp4}, see also \eqref{U}).

On the other hand, the liquid crystal system \eqref{1}--\eqref{5} satisfies a weak maximum principle on the molecular length $|\bfd|_{\mathbb{R}^n}$
(see Corollary \ref{oexe}) that plays an essential role in the mathematical analysis on its well-posedness and long-time behavior (see
\cite{B,C09,GW13} and also \cite{LL95} for the autonomous case). However, the validity of the weak maximum principle requires some additional constraints
on the $L^\infty$-norms of the boundary and initial data $\mathbf{h}$, $\bfd_0$ (see \eqref{hyp5} and \eqref{hyp5bis}), which will bring extra difficulties to the
study of the control problem (\textbf{CP}). For instance, the differentiability of the control-to-state mapping $\mathcal{S}$ with respect to the boundary
control $\mathbf{h}$ should be considered in a certain bounded convex set.
Besides, in addition to the control constraints, a state constraint is also to be respected. In order to avoid these difficulties, the existing well-posedness and regularity results in
\cite{B,GW13} must be refined. Here, we shall show that the constraints \eqref{hyp5} and \eqref{hyp5bis} are indeed not necessary for the well-posedness of
problem \eqref{1}--\eqref{5}, with a cheap price to be paid on the slightly higher integrability in time for the time derivative of $\mathbf{h}$
(compare \eqref{hyp3} and \eqref{hyp2w}, see also Theorem \ref{exe2d}).
The well-posedness result obtained in Theorem \ref{exe2d} implies that our assumptions on the regularity of trace spaces for the boundary control $\mathbf{h}$
seem to be almost minimal. It will be an interesting problem to reduce the regularity requirements on $\mathbf{h}$ and thus include less regular boundary
controls to problem (\textbf{CP}).

Our contribution can be viewed as a first step towards the study on optimal control problems related to liquid crystal flows.
We note that the PDE system \eqref{1}--\eqref{3} is highly simplified and it only keeps some essential properties of the original Ericksen--Leslie system \cite{LL95,E1,Le}.
It will be interesting to investigate optimal boundary (and also distributed) control problems for more comprehensive liquid crystal systems with important
physical considerations as, for instance, the rotation/streching effects and the nonlinear constraint $|\bfd|_{\mathbb{R}^n}=1$.
These will be the subjects of our future study.

The plan of the paper is as follows: in Section 2, we present some preliminary results concerning
the well-posedness of problem \eqref{1}--\eqref{5} under suitable regularity and compatibility assumptions on the Dirichlet boundary data $\mathbf{h}$.
Then we prove some stability estimates for global strong solutions in higher-order Sobolev norms
when the spatial dimension is two, which are crucial for the analysis of the optimal control problem (\textbf{CP}).
In Section 3, we show the existence of an optimal boundary control over a suitable admissible set.
Section 4 is devoted to the Fr\'echet differentiability of the control-to-state operator $\mathcal{S}$.
In Section 5 and Section 6, some first-order necessary optimality conditions for the problem (\textbf{CP}) are derived.

\section{Preliminaries}
\label{Sec2} \setcounter{equation}{0}

\subsection{Functional settings}
Without loss of generality, throughout the paper, we simply set
 $$\nu=\lambda=\eta=\epsilon=1,$$
because the values of those coefficients do not play a role in the subsequent analysis.

Let us introduce the function spaces we shall work with.  As
usual, $L^p(\Omega)$ and $W^{k,p}(\Omega)$ stand for the Lebesgue
and the Sobolev spaces of real valued functions, with the convention
that $H^k(\Omega)= W^{k,2}(\Omega)$. The spaces of vector-valued
functions are denoted by bold letters, correspondingly.
We set
\be
\mathbf{H}=\overline{\mathcal{V}}^{\mathbf{L}^2(\Omega)},\quad
\mathbf{V}=\overline {\mathcal{V}}^{\mathbf{H}_0^1(\Omega)},\quad
\text{where}\ \mathcal{V}=\left\{\mathbf{v}\in
C_0^\infty(\Omega,\mathbb{R}^n)
 :\, \nabla\cdot \mathbf{v}=0\right\}
 \label{spaces}
 \ee
that are the classical Hilbert spaces for the incompressible Navier--Stokes equations subject to no-slip boundary
conditions (see \cite[Chapter 1]{Te1}). Their norms are given by $\|\cdot\|_{\mathbf{L}^2(\Omega)}$ and $\|\cdot\|_{\mathbf{H}^1(\Omega)}$, respectively.
For any Banach space $B$, we denote its dual space by $B'$.
The notations $\langle\cdot, \cdot\rangle_{B', B}$ and $\|\cdot\|_{B}$ will stand for the duality pairing between $B$ and its dual $B'$,
and for the norm of $B$, respectively.
In particular, we denote the dual space of $\mathbf{H}_0^1(\Omega)$ by  $\mathbf{H}^{-1}(\Omega)$.
We recall that the operator  $-\Delta$ with homogeneous Dirichlet boundary condition is strictly positive and self-adjoint
in $\mathbf{L}^2(\Omega)$. So that the spectral theorem allows us to define the
powers $(-\Delta)^s$ and the associated spaces $D((-\Delta)^s)$, for $s \in\mathbb{R}$.
Denote $\mathbf{H}^s(\Omega)=D((-\Delta)^\frac{s}{2})$. We know that
$$D((-\Delta)^1)={\mathbf H}^2(\Omega)\cap {\mathbf H}^1_0(\Omega),\quad D((-\Delta)^\frac{1}{2}={\mathbf H}^1_0(\Omega),
\quad D((-\Delta)^0)={\mathbf L}^2(\Omega).$$
If it is not misleading, we will use the shorthand notations $\|\cdot\|_{{\mathbf L}^2}, \,  \|\cdot\|_{{\mathbf H}^1}, \dots$ to indicate the norms
defined in the domain $\Omega$, that is $\|\cdot\|_{{\mathbf L}^2(\Omega)}, \,  \|\cdot\|_{{\mathbf H}^1(\Omega)}, \dots$

We recall here the regularity result for the Stokes problem (see, e.g., \cite[Chapter 1, Proposition 2.2]{Te1}):
  \bl \label{S}
  Consider  the Stokes operator $S: D(S)= \mathbf{ V}\cap\mathbf{H}^2(\Omega) \to
  \mathbf{ H} $ defined by
  $$S\mathbf{u}=-\Delta \mathbf{u} +\nabla \pi \in \ \mathbf{H}, \quad
  \forall\, \mathbf{u} \in D(S),$$
where $\pi \in H^1(\Omega)$.
Then it holds
  $$\|\mathbf{u}\|_{\mathbf{H}^2}+\|\pi\|_{H^1 / \mathbb{R} }
  \leq C\|S \mathbf{u}\|_{{\mathbf L}^2}, \quad \forall\, \mathbf{u} \in D(S),$$
for some positive constant $C$ only depending on the domain $\Omega$ and the spatial dimension $n$.
 \el


\subsection{Well-posedness of the state problem \eqref{1}--\eqref{5}}
We start with an existence result for global weak solutions to problem \eqref{1}--\eqref{5} in both two and three dimensions.
 \bt [Global weak solutions] \label{we} Let $n=2,3$.  For any given $T>0$, assume that
 \begin{align}
  \label{hyp2}
 &\mathbf{h}\in L^2(0, T; \mathbf{H}^{\frac32}(\Gamma)),\\
  \label{hyp3}
 &\partial_t\mathbf{h}\in L^4 (0,T; \mathbf{H}^{-\frac12}(\Gamma)).
\end{align}
Then for any initial datum $(\mathbf{v}_0, \mathbf{d}_0)\in \mathbf{H}\times
 \mathbf{H}^1(\Omega)$
satisfying the compatibility condition
 \begin{align}
 \label{hyp4}
&\mathbf{d}_0|_\Gamma=\mathbf{h}|_{t=0},
\end{align}
 problem \eqref{1}--\eqref{5} admits a global
 weak solution $(\mathbf{v}, \mathbf{d})$  such that
\begin{align}
&\mathbf{v}\in L^\infty(0, T; \mathbf{H})\cap L^2(0, T; \mathbf{V}),\non\\
&\mathbf{d}\in L^\infty(0, T; \mathbf{H}^1(\Omega)) \cap L^2(0, T; \mathbf{H}^2(\Omega)).\non
\end{align}
\et
If one considers global weak solutions with $L^\infty$-constraint on the length of molecular director $\bfd$, it holds
\bc [Global weak solutions with constraint $|\bfd|_{\mathbb{R}^n}\leq 1$]\label{oexe}
Let $n=2,3$ and $(\mathbf{v}_0, \mathbf{d}_0)\in \mathbf{H}\times \mathbf{H}^1(\Omega)$ such that \eqref{hyp4} is satisfied.
For any given $T>0$, assume \eqref{hyp2} and
\begin{align}
& \partial_t\mathbf{h}\in L^2 (0,T; \mathbf{H}^{-\frac12}(\Gamma)),         \label{hyp2w}\\
&|\mathbf{h}(x,t)|_{\mathbb{R}^n}\leq 1 \quad\ \text{a.e. on } \Sigma,      \label{hyp5}\\
&|\mathbf{d}_0(x)|_{\mathbb{R}^n}\leq 1\quad \ \  \text{a.e. in } \Omega.   \label{hyp5bis}
\end{align}
Then problem \eqref{1}--\eqref{5} admits a global
 weak solution $(\mathbf{v}, \mathbf{d})$ with the same regularity as in Theorem \ref{we}. Moreover, $\bfd$ satisfies the weak maximum principle
\begin{align}
 &|\mathbf{d}(x,t)|_{\mathbb{R}^n}\leq 1, \quad \textrm{ a.e. in } Q. \label{max}
\end{align}
 \ec

\begin{remark}
We note that the statement of Theorem \ref{we} is different from the existing results obtained in \cite[Corollary 1]{B} or
\cite[Theorem 7]{C09} that are similar to Corollary \ref{oexe}.
In particular, in Theorem \ref{we}, we avoid using the assumptions \eqref{hyp5}--\eqref{hyp5bis} by taking a slightly stronger assumption
on $\partial_t\mathbf{h}$ (compare \eqref{hyp3} and \eqref{hyp2w}). Theorem \ref{we} can be proved by using a semi-Galerkin approximation scheme similar to
\cite{B,C09}. However, the argument in \cite{B,C09} cannot be applied directly, because it essentially relies on the assumptions \eqref{hyp5}--\eqref{hyp5bis}
that lead to the weak maximum principle for $\bfd$ and the estimate \eqref{max}  plays a crucial role in controlling
the nonlinear term $\mathbf{f}(\bfd)$. Without these two assumptions, necessary modifications should be made
in order to overcome the difficulty due to the lack of control on $\|\bfd\|_{L^\infty(0,T;\mathbf{L}^\infty(\Omega))}$.
The proof of Theorem \ref{we} will be sketched in the Appendix.
\end{remark}

Due to the time-dependent boundary condition \eqref{4}, the system \eqref{1}--\eqref{5} no longer satisfies the dissipative energy law like the autonomous case in \cite{LL95}.
However, with the help of a suitable lifting function $\bfd_E$ (see \eqref{LE}), one can still derive a specific energy inequality (see Lemma \ref{BEL}).
Combining it with Lemmas \ref{Ap1}, \ref{dpe}, we can obtain some
uniform estimates for global weak solutions to problem \eqref{1}--\eqref{5} on arbitrary time interval $[0,T]$:
 \bp  \label{lowe} Let the assumptions of Theorem \ref{we} hold. Then, any global weak solution $(\mathbf{v}, \mathbf{d})$
 to problem \eqref{1}--\eqref{5} fulfills the following estimate
 \bea
 &&
 \|\mathbf{v}(t)\|_{\mathbf{L}^2}^2+ \|\mathbf{d}(t)\|_{\mathbf{H}^1}^2+\int_0^t(\|\mathbf{v}(\tau)\|_{\mathbf{H}^1}^2
 + \|\mathbf{d}(\tau)\|_{\mathbf{H}^2}^2) d\tau\leq C_T, \quad \forall\, t\in[0, T],
 \eea
where the positive constant $C_T$ depends on
$\|\mathbf{v}_0\|_{\mathbf{L}^2}$, $\|\mathbf{d}_0\|_{\mathbf{H}^1}$,  $\Vert \mathbf{h}\Vert_{L^2(0,T;\mathbf{H}^{\frac32}(\Gamma))}$,
  $\|\partial_t\mathbf{h}\|_{ L^4(0,T; \mathbf{H}^{\frac12}(\Gamma))}$, $\Omega$ and $T$.
 \ep

When the spatial dimension is two, further conclusions can be obtained.
First, we can prove the following continuous dependence result on initial and boundary data in the lower-order energy space
$\mathbf{L}^2(\Omega)\times \mathbf{H}^1(\Omega)$, which  easily yields the uniqueness of global weak solutions when $n=2$.

\bp[Continuous dependence in $\mathbf{H}\times \mathbf{H}^1(\Omega)$]  \label{uniw}
Let the assumptions of Theorem \ref{we} hold. If $n=2$, then problem \eqref{1}--\eqref{5} admits a unique  weak solution $(\mathbf{v}, \mathbf{d})$.
Moreover, let $(\mathbf{v}^{(i)}, \mathbf{d}^{(i)})$ $(i=1,2)$ be two weak solutions corresponding to the initial data
$(\mathbf{v}_0^{(i)}, \mathbf{d}_0^{(i)})$ as well as the boundary data $\mathbf{h}^{(i)}$.
Denoting the differences $\bar\bfv=\bfv^{(1)}-\bfv^{(2)}$, $\bar\bfd=\bfd^{(1)}-\bfd^{(2)}$, $\bar\bfv_0=\bfv^{(1)}_0-\bfv^{(2)}_0$,
$\bar\bfd_0=\bfd^{(1)}_0-\bfd^{(2)}_0$
and $\bar{\mathbf{h}}= \mathbf{h}^{(1)}-\mathbf{h}^{(2)}$,
then for $t\in [0,T]$ the following estimate holds:
\begin{align}
&\|\bar\bfv(t)\|_{\mathbf{L}^2}^2+\|\nabla \bar\bfd(t)\|_{\mathbf{L}^2}^2
+\int_0^t(\|\nabla \bar\bfv(\tau)\|_{\mathbf{L}^2}^2+\|\Delta \bar\bfd(\tau)\|^2_{\mathbf{L}^2}) d\tau\non\\
 &\quad \leq C_T\Big[\|\bar\bfv_0\|_{\mathbf{L}^2}^2+\|\nabla \bar\bfd_0\|_{\mathbf{L}^2}^2
 + \int_0^t \Big(\|\bar{\mathbf{h}}(\tau)\|_{\mathbf{H}^\frac32(\Gamma)}^2
 +\|\partial_t\bar{\mathbf{h}}(\tau)\|_{\mathbf{H}^{-\frac12}(\Gamma)}^2\Big) d\tau\Big],\label{Psi}
\end{align}
where the constant $C_T>0$ depends on
$\|\mathbf{v}_0\|^{(i)}_{\mathbf{L}^2}$, $\|\mathbf{d}_0\|^{(i)}_{\mathbf{H}^1}$,  $\Vert \mathbf{h}^{(i)}\Vert_{L^2(0,T;\mathbf{H}^{\frac32}(\Gamma))}$,
  $\|\partial_t\mathbf{h}^{(i)}\|_{ L^4(0,T; \mathbf{H}^{\frac12}(\Gamma))}$, $\Omega$ and $T$.
 \ep
 \begin{proof}
 We recall that a similar result was proven in \cite[Theorem 2.4]{B} under the additional assumptions \eqref{hyp5}--\eqref{hyp5bis} that are not valid in our case.
 Nevertheless, the only difference in the proof is related to the treatment of the nonlinear term $\mathbf{f}(\bfd)$, namely,
 $\int_\Omega (\mathbf{f}(\bfd^{(1)})-\mathbf{f}(\bfd^{(2)}))\cdot \Delta \bar{\bfd} dx$.

 We take $\bfd_E^{(i)}$ $(i=1,2)$ as the elliptic lifting functions (see Appendix) given by
\be
 \begin{cases}
 -\Delta \mathbf{d}^{(i)}_E=\mathbf{0},\qquad \text{ in } \Omega \times (0,T),\\
 \mathbf{d}^{(i)}_E=\mathbf{h}^{(i)},\qquad\  \ \text{ on } \Gamma \times (0,T).
 \end{cases}
 \label{LEi}
 \ee
 Then we set $\wha{\bfd}^{(i)}=\bfd^{(i)}-\bfd_E^{(i)}$ ($i=1,2$), which satisfy the homogeneous Dirichlet boundary condition $\wha{\bfd}^{(i)}|_\Gamma=\mathbf{0}$.
Instead of using the weak maximum principle for $\bfd$ (cf. \eqref{max}), we deduce from Proposition \ref{lowe}, the Sobolev embedding theorem ($n=2$),
 the H\"older, Poincar\'e and Young's inequalities that
 \begin{align}
 &\left |\int_\Omega (\mathbf{f}(\bfd^{(1)})-\mathbf{f}(\bfd^{(2)}))\cdot \Delta \bar{\bfd} dx\right |\non\\
 &\quad \leq 2(\|\bfd^{(1)}\|_{\mathbf{L}^6}^2+\|\bfd^{(2)}\|_{\mathbf{L}^6}^2)\|\bar{\bfd}\|_{\mathbf{L}^6}\|\Delta \bar{\bfd}\|_{\mathbf{L}^2}
 +\|\bar{\bfd}\|_{\mathbf{L}^2}\|\Delta \bar{\bfd}\|_{\mathbf{L}^2}\non\\
 &\quad \leq C_T\|\bar{\bfd}\|_{\mathbf{H}^1}\|\Delta \bar{\bfd}\|_{\mathbf{L}^2}\non\\
 &\quad \leq C_T(\|\wha{\bfd}^{(1)}-\wha{\bfd}^{(2)}\|_{\mathbf{H}^1}+\|\bfd_E^{(1)}-\bfd_E^{(2)}\|_{\mathbf{H}^1})\|\Delta \bar{\bfd}\|_{\mathbf{L}^2}\non\\
 &\quad \leq C_T(\|\nabla(\wha{\bfd}^{(1)}-\wha{\bfd}^{(2)})\|_{\mathbf{L}^2}
 +\|\bfd_E^{(1)}-\bfd_E^{(2)}\|_{\mathbf{H}^1})\|\Delta \bar{\bfd}\|_{\mathbf{L}^2}\non\\
 &\quad \leq C_T(\|\nabla \bar{\bfd}\|_{\mathbf{L}^2}+\|\nabla (\bfd_E^{(1)}-\bfd_E^{(2)})\|_{\mathbf{L}^2}
 +\|\bfd_E^{(1)}-\bfd_E^{(2)}\|_{\mathbf{H}^1})\|\Delta \bar{\bfd}\|_{\mathbf{L}^2}\non\\
 &\quad \leq \varepsilon \|\Delta \bar{\bfd}\|_{\mathbf{L}^2}^2+C_T(\|\nabla \bar{\bfd}\|_{\mathbf{L}^2}^2
 +\|\bar{\mathbf{h}}\|^2_{\mathbf{H}^\frac32(\Gamma)}),\non
 \end{align}
 for certain sufficiently small constant $\varepsilon>0$.
 Keeping the above estimate in mind, we can follow the  argument  as \cite[Theorem 2.4]{B} to prove our conclusion.
 \end{proof}
 Next, if the initial data is more regular, namely, $(\mathbf{v}_0, \mathbf{d}_0)\in \mathbf{V}\times
 \mathbf{H}^2(\Omega)$, we can further prove the existence of a unique global strong solution to problem
\eqref{1}--\eqref{5} in two spatial dimensions:
\bt[Global strong solutions in $2D$] \label{exe2d}
 Let $n=2$. For any $T >0$, assume
  \begin{align}
  \label{hyp2s}
 &\mathbf{h}\in L^2(0, T; \mathbf{H}^{\frac52}(\Gamma)),\\
 \label{hyp3s}
 &\partial_t\mathbf{h}\in L^4 (0,T; \mathbf{H}^{\frac12}(\Gamma)).
\end{align}
 Let  $(\mathbf{v}_0, \mathbf{d}_0)\in \mathbf{V}\times
 \mathbf{H}^2(\Omega)$ such that the compatibility condition \eqref{hyp4} holds.
  Then problem \eqref{1}--\eqref{5} admits a unique global strong solution $(\mathbf{v}, \mathbf{d})$ such that
 \begin{align}
&\mathbf{v}\in L^\infty(0, T; \mathbf{V})\cap L^2(0, T; \mathbf{H}^2(\Omega)),\non\\
&\mathbf{d}\in L^\infty(0, T; \mathbf{H}^2(\Omega)) \cap L^2(0, T; \mathbf{H}^3(\Omega)).\non
\end{align}
Moreover, the following estimates hold
 \bea
 &&\|\mathbf{v}(t)\|_{\mathbf{H}^1}
   + \|\mathbf{d}(t)\|_{\mathbf{H}^2}\leq C_T, \quad \forall\, t\in[0,T], \label{globstrong1}
 \\
 &&\int_0^t\Big(\|\mathbf{v}(\tau)\|_{\mathbf{H}^2}^2+\|\mathbf{d}(\tau)\|_{\mathbf{H}^3}^2\Big)d\tau
 \leq C_T,\quad \forall\, t \in [0,T], \label{globstrong2}
 \eea
 where $C_T>0$ depends on $\|\mathbf{v}_0\|_{\mathbf{H}^1}$,
 $\|\mathbf{d}_0\|_{\mathbf{H}^2}$,
    $\|\mathbf{h}\|_{L^2(0,T;\mathbf{H}^\frac52(\Gamma))}$, $\|\partial_t\mathbf{h}\|_{L^4(0,T; \mathbf{H}^\frac12(\Gamma))}$, $\Omega$, and $T$.
 \et
 \begin{remark} \label{rmexe2d}
 The proof of Theorem \ref{exe2d} is similar to \cite[Theorem 2.6]{B} and \cite[Theorem 2.7 (ii)]{GW13} with some necessary refined estimates
 without using the weak maximum principle \eqref{max} (see Proposition \ref{esAP}). A sketch of the proof will be given in the Appendix.

 Besides, by exploiting the equations \eqref{1} and \eqref{3}, one can also verify that in the two dimensional case the global strong solution
 $(\mathbf{v}, \mathbf{d})$ satisfies the following regularity in time:
 \begin{align}
&\partial_t\mathbf{v}\in L^2(0, T; \mathbf{H}),\quad \partial_t \mathbf{d}\in L^\infty(0, T; \mathbf{L}^2(\Omega)) \cap L^2(0, T; \mathbf{H}^1(\Omega)),\non
\end{align}
which further implies the continuity property $\mathbf{v}\in C([0, T]; \mathbf{V})$, $\mathbf{d}\in C([0, T]; \mathbf{H}^2(\Omega))$, by means of the
interpolation (see e.g., \cite{SI}). Moreover, it holds
\be
 \|\partial_t\bfd(t)\|_{\mathbf{L}^2}^2
 +\int_0^t\Big(\|\partial_t \mathbf{v}(\tau)\|_{\mathbf{L}^2}^2+\|\partial_t\mathbf{d}(\tau)\|_{\mathbf{H}^1}^2\Big)d\tau
 \leq C_T,\quad \forall\, t \in [0,T], \label{globstrong3}
 \ee
 where $C_T>0$ depends on $\|\mathbf{v}_0\|_{\mathbf{H}^1}$, $\|\mathbf{d}_0\|_{\mathbf{H}^2}$,
    $\|\mathbf{h}\|_{L^2(0,T;\mathbf{H}^\frac52(\Gamma))}$, $\|\partial_t\mathbf{h}\|_{L^4(0,T; \mathbf{H}^\frac12(\Gamma))}$, $\Omega$ and $T$.
 \end{remark}

Finally, in the two dimensional case,  we can deduce the following continuous dependence result on initial and boundary data
 for global strong solutions to problem \eqref{1}--\eqref{5} in the higher-order space $\mathbf{H}^1(\Omega)\times \mathbf{H}^2(\Omega)$, which will be crucial to
 prove the differentiability of the state-to-control operator $\mathcal{S}$ (see \eqref{ctsS}).
\bp[Continuous dependence in $\mathbf{V}\times \mathbf{H}^2(\Omega)$] \label{strconti}
Suppose that the assumptions of Theorem \ref{exe2d} are satisfied. Let $(\mathbf{v}^{(i)}, \mathbf{d}^{(i)})$ $(i=1,2)$ be two strong solutions
corresponding to the initial data $(\mathbf{v}_0^{(i)}, \mathbf{d}_0^{(i)})$ and the boundary data $\mathbf{h}^{(i)}$ on $[0,T]$.
Denoting the differences $\bar\bfv=\bfv^{(1)}-\bfv^{(2)}$, $\bar\bfd=\bfd^{(1)}-\bfd^{(2)}$, $\bar\bfv_0=\bfv^{(1)}_0-\bfv^{(2)}_0$, $\bar\bfd_0=\bfd^{(1)}_0-\bfd^{(2)}_0$
and $\bar{\mathbf{h}}= \mathbf{h}^{(1)}-\mathbf{h}^{(2)}$, then the following estimate holds:
\begin{align}
&\|\bar\bfv(t)\|_{\mathbf{H}^1}^2+\|\bar\bfd(t)\|_{\mathbf{H}^2}^2
+\int_0^t(\|\bar\bfv(\tau)\|_{\mathbf{H}^2}^2+\|\bar\bfd(\tau)\|_{\mathbf{H}^3}^2) d\tau\non\\
&\quad \leq C_T\Big[ \|\bar\bfv_0\|_{\mathbf{H}^1}^2+\|\bar\bfd_0\|_{\mathbf{H}^2}^2
+ \int_0^t \Big(\|\bar{\mathbf{h}}(\tau)\|_{\mathbf{H}^\frac52(\Gamma)}^2
+\|\partial_t\bar{\mathbf{h}}(\tau)\|_{\mathbf{H}^{\frac12}(\Gamma)}^2\Big) d\tau\Big],\quad \forall\, t\in [0,T],\label{conti}
\end{align}
where the positive constant $C_T$ depends on the norms $\|\mathbf{v}^{(i)}_0\|_{\mathbf{H}^1}$,
 $\|\mathbf{d}^{(i)}_0\|_{\mathbf{H}^2}$, $\|\mathbf{h}^{(i)}\|_{L^2(0,T;\mathbf{H}^\frac52(\Gamma))}$,
 $\|\partial_t\mathbf{h}^{(i)}\|_{L^4(0,T; \mathbf{H}^\frac12(\Gamma))}$, $\Omega$ and $T$.
\ep
\begin{proof}
Let $\bfd_P^{(i)}$ $(i=1,2)$ be the parabolic lifting functions satisfying
\be
 \begin{cases}
 \p_t\mathbf{d}^{(i)}_P-\Delta \mathbf{d}^{(i)}_P=\mathbf{0},\qquad \text{ in } \Omega\times \mathbb{R}^+,\\
 \mathbf{d}^{(i)}_P=\mathbf{h}^{(i)},\qquad  \qquad \quad \ \text{ on } \Gamma\times \mathbb{R}^+,\\
 \mathbf{d}^{(i)}_P(0)= \mathbf{d}^{(i)}_{E0}, \qquad \qquad \text{ in } \Omega,
 \end{cases}\label{LPi}
 \ee
where the initial data $\mathbf{d}^{(i)}_{E0}$ $(i=1,2)$ are respectively given by
 \be
 \begin{cases}
 -\Delta \mathbf{d}^{(i)}_{E0}=\mathbf{0},\qquad \text{ in } \Omega,\\
 \mathbf{d}^{(i)}_{E0}=\mathbf{d}^{(i)}_0,\qquad \ \ \text{ on } \Gamma.
 \end{cases}\label{iE0i}
 \ee
Then we set
$$\bar{\bfd}_P=\bfd_P^{(1)}-\bfd^{(2)}_P\quad \text{and}\quad \wt\bfd^{(i)}=\bfd^{(i)}-\mathbf{d}^{(i)}_P,\quad i=1,2.$$
In particular, it holds that
$$\wt\bfd^{(i)}=\mathbf{0}\quad \text{ on }\Gamma,\quad i=1,2.$$

Using the estimates \eqref{globstrong1}, \eqref{globstrong2} and \eqref{globstrong3} on $[0,T]$ for the global strong solutions
$(\mathbf{v}^{(i)}, \mathbf{d}^{(i)})$ $(i=1,2)$, after performing similar calculations as in the proof of \cite[Theorem 4.7, pp. 427-432]{B},
we can obtain the following differential inequality
\begin{align}
& \frac{d}{dt}\mathcal{Y}(t)+\frac12 \mathcal{Z}(t)
   \leq C_T\mathcal{Y}(t) +C_T\Big(\|\bar{\bfd}\|_{\mathbf{L}^2}^2
            +\|\bar{\mathbf{h}}\|_{\mathbf{H}^\frac52(\Gamma)}^2
            +\|\partial_t\bar{\mathbf{h}}\|_{\mathbf{H}^\frac12(\Gamma)}^2\Big),\label{contiH1H2}
\end{align}
where
\begin{align}
&\mathcal{Y}(t)=\|\nabla \bar{\bfv}\|_{\mathbf{L}^2}^2
+ \|\Delta (\wt\bfd^{(1)}-\wt\bfd^{(2)})-\mathbf{f}(\bfd^{(1)})+\mathbf{f}(\bfd^{(2)})\|_{\mathbf{L}^2}^2
+ \|\Delta \bar{\mathbf{d}}_P\|_{\mathbf{L}^2}^2,\non\\
&\mathcal{Z}(t)=\|\Delta \bar{\bfv}\|_{\mathbf{L}^2}^2
+ \|\nabla [\Delta (\wt\bfd^{(1)}-\wt\bfd^{(2)})-\mathbf{f}(\bfd^{(1)})+\mathbf{f}(\bfd^{(2)})]\|_{\mathbf{L}^2}^2
+ \|\bar{\mathbf{d}}_P\|_{\mathbf{H}^3}^2,\non
\end{align}
and $C_T>0$ depends on $\|\mathbf{v}^{(i)}_0\|_{\mathbf{H}^1}$,
 $\|\mathbf{d}^{(i)}_0\|_{\mathbf{H}^2}$, $\|\mathbf{h}^{(i)}\|_{L^2(0,T;\mathbf{H}^\frac52(\Gamma))}$,
 $\|\partial_t\mathbf{h}^{(i)}\|_{L^4(0,T; \mathbf{H}^\frac12(\Gamma))}$, $\Omega$ and $T$.

Then, from \eqref{contiH1H2} and the Gronwall lemma it follows that
\begin{align}
\mathcal{Y}(t)&\leq \mathcal{Y}(0) e^{C_T t}+ C_T\int_0^t \Big(\|\bar{\bfd}(\tau)\|_{\mathbf{L}^2}^2
            +\|\bar{\mathbf{h}}(\tau)\|_{\mathbf{H}^\frac52(\Gamma)}^2
            +\|\partial_t\bar{\mathbf{h}}(\tau)\|_{\mathbf{H}^\frac12(\Gamma)}^2\Big)d\tau, \quad \forall \, t \in [0,T].\label{YY}
\end{align}
Similar to the estimate \eqref{A3h1}, we have the following bound for the difference $\bar{\mathbf{d}}_P$ of lifting functions
\begin{align}
 \|\bar{\mathbf{d}}_P(t)\|^2_{\mathbf{H}^1}
  &   \leq \|\bar{\bfd}_{0}\|_{\mathbf{H}^1}^2
       + c\int_0^t\Big(\|\bar{\mathbf{h}}(\tau)\|_{\mathbf{H}^{\frac32}(\Gamma)}^2
       +\|\partial_t\bar{\mathbf{h}}(\tau)\|_{\mathbf{H}^{-\frac12}(\Gamma)}^2\Big)d\tau,\label{A3h1d}
 \end{align}
which together with the lower-order continuous dependence result in Proposition \ref{uniw} implies, for $ t\in [0,T]$,
\begin{align}
\|\bar{\bfd}(t)\|_{\mathbf{H}^1}^2
&  \leq 2\|(\wt\bfd^{(1)}-\wt\bfd^{(2)})(t)\|_{\mathbf{H}^1}^2+2\|\bar{\mathbf{d}}_P(t)\|^2_{\mathbf{H}^1}\non\\
& \leq C\|\nabla(\wt\bfd^{(1)}-\wt\bfd^{(2)})(t)\|_{\mathbf{L}^2}^2+2\|\bar{\mathbf{d}}_P(t)\|^2_{\mathbf{H}^1} \non\\
&  \leq C\|\nabla\bar{\bfd}(t)\|_{\mathbf{L}^2}^2+C\|\bar{\mathbf{d}}_P(t)\|^2_{\mathbf{H}^1}\non\\
&  \leq C_T\Big[\|\bar\bfv_0\|_{\mathbf{L}^2}^2+\|\bar\bfd_0\|_{\mathbf{H}^1}^2
+ \int_0^t \Big(\|\bar{\mathbf{h}}(\tau)\|_{\mathbf{H}^\frac32(\Gamma)}^2
+\|\partial_t\bar{\mathbf{h}}(\tau)\|_{\mathbf{H}^{-\frac12}(\Gamma)}^2\Big) d\tau\Big],\label{contiL2a}
\end{align}
\begin{align}
\int_0^t \|\bar{\bfd}(\tau)\|_{\mathbf{H}^2}^2 d\tau&\leq  C\int_0^t \Big(\|\Delta \bar{\bfd}(\tau)\|_{\mathbf{L}^2}^2
+ \|\bar{\mathbf{h}}(\tau)\|_{\mathbf{H}^\frac32(\Gamma)}^2\Big) d\tau \non\\
&  \leq C_T\Big[\|\bar\bfv_0\|_{\mathbf{L}^2}^2+\|\nabla \bar\bfd_0\|_{\mathbf{L}^2}^2
       + \int_0^t \Big(\|\bar{\mathbf{h}}(\tau)\|_{\mathbf{H}^\frac32(\Gamma)}^2
       +\|\partial_t\bar{\mathbf{h}}(\tau)\|_{\mathbf{H}^{-\frac12}(\Gamma)}^2\Big) d\tau\Big],\label{contiL2b}
\end{align}
where the constant $C_T>0$ depends on
$\|\mathbf{v}_0\|^{(i)}_{\mathbf{L}^2}$, $\|\mathbf{d}_0\|^{(i)}_{\mathbf{H}^1}$,  $\Vert \mathbf{h}^{(i)}\Vert_{L^2(0,T;\mathbf{H}^{\frac32}(\Gamma))}$,
  $\|\partial_t\mathbf{h}^{(i)}\|_{ L^4(0,T; \mathbf{H}^{\frac12}(\Gamma))}$, $\Omega$ and $T$.

Next, using Proposition \ref{lowe},  estimate \eqref{globstrong1} and the Sobolev embedding theorem ($n=2$), we obtain that
\begin{align}
&\|\Delta (\wt\bfd^{(1)}-\wt\bfd^{(2)})-\mathbf{f}(\bfd^{(1)})+\mathbf{f}(\bfd^{(2)})\|_{\mathbf{L}^2}\non\\
&\quad \leq \|\Delta (\wt\bfd^{(1)}-\wt\bfd^{(2)})\|_{\mathbf{L}^2}+\|\mathbf{f}(\bfd^{(1)})-\mathbf{f}(\bfd^{(2)})\|_{\mathbf{L}^2}\non\\
&\quad \leq C\|\bar{\bfd}\|_{\mathbf{H}^2}+C\|\bar{\mathbf{d}}_P\|_{\mathbf{H}^2}+C_T\|\bar{\bfd}\|_{\mathbf{L}^2},\label{esy1a}
\end{align}
\begin{align}
\|\bar{\bfd}\|_{\mathbf{H}^2}
 &\leq \|\wt\bfd^{(1)}-\wt\bfd^{(2)}\|_{\mathbf{H}^2}+ \|\bar{\mathbf{d}}_P\|_{\mathbf{H}^2}\non\\
 &\leq C\|\Delta (\wt\bfd^{(1)}-\wt\bfd^{(2)})\|_{\mathbf{L}^2}+ \|\bar{\mathbf{d}}_P\|_{\mathbf{H}^2}\non\\
 &\leq C\|\Delta (\wt\bfd^{(1)}-\wt\bfd^{(2)})-\mathbf{f}(\bfd^{(1)})+\mathbf{f}(\bfd^{(2)})\|_{\mathbf{L}^2}\non\\
 &\quad       +C\|\mathbf{f}(\bfd^{(1)})-\mathbf{f}(\bfd^{(2)})\|_{\mathbf{L}^2} +\|\bar{\mathbf{d}}_P\|_{\mathbf{H}^2}\non\\
 &\leq C\|\Delta (\wt\bfd^{(1)}-\wt\bfd^{(2)})-\mathbf{f}(\bfd^{(1)})+\mathbf{f}(\bfd^{(2)})\|_{\mathbf{L}^2}
       +C_T\|\bar{\bfd}\|_{\mathbf{L}^2}+\|\bar{\mathbf{d}}_P\|_{\mathbf{H}^2},\label{esy1b}
\end{align}
and
\begin{align}
\|\bar{\bfd}\|_{\mathbf{H}^3}
 &\leq \|\wt\bfd^{(1)}-\wt\bfd^{(2)}\|_{\mathbf{H}^3}+ \|\bar{\mathbf{d}}_P\|_{\mathbf{H}^3}\non\\
 &\leq C\|\nabla\Delta (\wt\bfd^{(1)}-\wt\bfd^{(2)})\|_{\mathbf{L}^2}+C\|\Delta (\wt\bfd^{(1)}-\wt\bfd^{(2)})\|_{\mathbf{L}^2}
      + \|\bar{\mathbf{d}}_P\|_{\mathbf{H}^3}\non\\
 &\leq C\|\nabla [\Delta (\wt\bfd^{(1)}-\wt\bfd^{(2)})-\mathbf{f}(\bfd^{(1)})+\mathbf{f}(\bfd^{(2)})]\|_{\mathbf{L}^2}
       +C\|\mathbf{f}(\bfd^{(1)})-\mathbf{f}(\bfd^{(2)})\|_{\mathbf{H}^1}\non\\
 &\qquad    + \|\bar{\mathbf{d}}_P\|_{\mathbf{H}^3} \non\\
 &\leq C\|\nabla [\Delta (\wt\bfd^{(1)}-\wt\bfd^{(2)})-\mathbf{f}(\bfd^{(1)})+\mathbf{f}(\bfd^{(2)})]\|_{\mathbf{L}^2}
       +C_T\|\bar{\bfd}\|_{\mathbf{H}^2}+\|\bar{\mathbf{d}}_P\|_{\mathbf{H}^3},\label{esy1bb}
\end{align}
where $C>0$ is a constant depending on $\Omega$ and $C_T>0$ is a constant depending on
$\|\mathbf{v}_0\|_{\mathbf{H}^1}$, $\|\mathbf{d}_0\|_{\mathbf{H}^2}$,  $\Vert \mathbf{h}\Vert_{L^2(0,T;\mathbf{H}^{\frac52}(\Gamma))}$,
  $\|\partial_t\mathbf{h}\|_{ L^4(0,T; \mathbf{H}^{\frac12}(\Gamma))}$, $\Omega$ and $T$.

On the other hand, we infer from the estimates \eqref{LPi}, \eqref{iE0i} and \eqref{esy1a} that
\begin{align}
\mathcal{Y}(0)&\leq \|\nabla \bar{\bfv}_0\|_{\mathbf{L}^2}^2+3\|\bar{\bfd}_0\|_{\mathbf{H}^2}^2
+ 3\|\mathbf{d}^{(1)}_{E0}-\mathbf{d}^{(2)}_{E0}\|_{\mathbf{H}^2}^2+3C\|\bar{\bfd}_0\|_{\mathbf{L}^2}^2
+\|\Delta(\mathbf{d}^{(1)}_{E0}-\mathbf{d}^{(2)}_{E0})\|_{\mathbf{L}^2}^2\non\\
&\leq \|\nabla \bar{\bfv}_0\|_{\mathbf{L}^2}^2+ C\|\bar{\bfd}_0\|_{\mathbf{H}^2}^2,\label{esy0}
\end{align}
where $C>0$ is a constant depending on $\Omega$ and the coefficients of the system \eqref{1}--\eqref{3}.

As a consequence, from \eqref{YY} and \eqref{contiL2b} it follows that, for $t\in [0,T]$,
\begin{align}
\mathcal{Y}(t)\leq C_T\Big[\|\bar{\bfv}_0\|_{\mathbf{H}^1}^2+ \|\bar{\bfd}_0\|_{\mathbf{H}^2}^2
+\int_0^t \Big(\|\bar{\mathbf{h}}(\tau)\|_{\mathbf{H}^\frac52(\Gamma)}^2
            +\|\partial_t\bar{\mathbf{h}}(\tau)\|_{\mathbf{H}^\frac12(\Gamma)}^2\Big)d\tau\Big],\label{Ytes}
\end{align}
while integrating \eqref{contiH1H2} with respective to time, we also have, for $t\in [0,T]$,
\begin{equation}
\int_0^t \mathcal{Z}(\tau) d\tau\leq C_T\Big[\|\bar{\bfv}_0\|_{\mathbf{H}^1}^2+ \|\bar{\bfd}_0\|_{\mathbf{H}^2}^2
+\int_0^t \Big(\|\bar{\mathbf{h}}(\tau)\|_{\mathbf{H}^\frac52(\Gamma)}^2
            +\|\partial_t\bar{\mathbf{h}}(\tau)\|_{\mathbf{H}^\frac12(\Gamma)}^2\Big)d\tau\Big].\label{Ztes}
\end{equation}

Finally, applying the estimate \eqref{A4h2} to the difference $\mathbf{d}^{(1)}_P-\mathbf{d}^{(2)}_P$,
which satisfies a linear parabolic equation similar to \eqref{LP}, we get
\begin{align}
  &\|\bar{\mathbf{d}}_P(t)\|^2_{\mathbf{H}^2}
    +\int_0^t \|\bar{\mathbf{d}}_P(\tau)\|^2_{\mathbf{H}^3} d\tau  \leq \|\bar{\bfd}_0\|_{\mathbf{H}^2}^2
       + c\int_0^t\Big(\|\bar{\mathbf{h}}(\tau)\|_{\mathbf{H}^{\frac52}(\Gamma)}^2
       +\|\partial_t\bar{\mathbf{h}}(\tau)\|_{\mathbf{H}^{\frac12}(\Gamma)}^2\Big)d\tau.\label{A4h212}
 \end{align}
Then our conclusion \eqref{conti} follows from \eqref{Ytes} and \eqref{Ztes} together with the estimates
\eqref{contiL2a}, \eqref{contiL2b}, \eqref{esy1b}, \eqref{esy1bb} and \eqref{A4h212}.

The proof is complete.
\end{proof}

\section{The Optimal Control Problem: Existence}
\label{OPCO}
\setcounter{equation}{0}

In this section we investigate the optimal control problem (\textbf{CP}) with the cost functional $\mathcal{J}$ given by \eqref{costJ}
when the spatial dimension $n=2$.

\subsection{Space of control functions and admissible sets}
Let $T\in (0,+\infty)$ be an arbitrary but fixed time.
Motivated by the existence and continuous dependence results on global strong solutions to problem \eqref{1}--\eqref{5}
(see Theorem \ref{exe2d} and Proposition \ref{strconti}), the boundary control $\mathbf{h}$ should be measured in a norm
that is not weaker than those required therein.
As a result, we introduce the space
\begin{align}
\mathcal{U}& :=\big\{\mathbf{h}(x,t)\, \mid\, \mathbf{h}\in L^2(0,T;\mathbf{H}^\frac52(\Gamma)),\
\partial_t\mathbf{h}\in L^{4}(0,T; \mathbf{H}^\frac12(\Gamma))\big\},\label{U1}
\end{align}
whose norm is given by
$$
\|\mathbf{h}\|_{\mathcal{U}}:=\|\mathbf{h}\|_{L^2(0,T;\mathbf{H}^\frac52(\Gamma))}
+\|\partial_t\mathbf{h}\|_{L^{4}(0,T; \mathbf{H}^\frac12(\Gamma))}, \quad \forall\,\mathbf{h}\in \mathcal{U}.
$$
Next, assume that $(\mathbf{v}_0, \mathbf{d}_0)\in \mathbf{V}\times
 \mathbf{H}^2(\Omega)$ is an arbitrary given initial datum.
  Then we define the \emph{affine space of control functions} $\wt{\mathcal{U}}$ (associated with $(\mathbf{v}_0, \mathbf{d}_0)$) such that
\begin{align}
\wt{\mathcal{U}}& :=\big\{\mathbf{h}(x,t)\, \mid\, \mathbf{h}\in\mathcal{U}  \  \text{with}\ \mathbf{h}(x,t)|_{t=0}=\bfd_0(x)|_\Gamma\big\}.\label{U}
\end{align}
\begin{remark}
By the continuous embeddings ($n=2$)
$$L^2(0,T;\mathbf{H}^\frac52(\Gamma))\cap W^{1,4}(0,T; \mathbf{H}^\frac12(\Gamma))\hookrightarrow C([0,T]; \mathbf{H}^\frac32(\Gamma))
\hookrightarrow C(\Gamma\times [0,T]),$$
we see that any $\mathbf{h}\in \wt{\mathcal{U}}$ is a continuous vector defined on $\Gamma\times [0,T]$
and thus the compatibility condition $\mathbf{h}(x,t)|_{t=0}=\bfd_0(x)|_\Gamma$ in \eqref{U} makes sense.
\end{remark}
It is natural that the size of the boundary controls should be constrained from both the physical and mathematical point of view.
To this end, we introduce the admissible set for problem (\textbf{CP}):
\begin{definition}
Let $M\in (0,+\infty)$ be a prescribed positive constant. The set of \emph{admissible boundary control functions}
$\wt{\mathcal{U}}_{\mathrm{ad}}^M$ is defined as follows
\begin{align}
\wt{\mathcal{U}}_{\mathrm{ad}}^M:=
&\big\{\mathbf{h}\in \mathcal{U} \, \mid\, \|\mathbf{h}\|_{L^2(0,T;\mathbf{H}^\frac52(\Gamma))}\leq M,
\quad \|\partial_t\mathbf{h}\|_{L^4(0,T; \mathbf{H}^\frac12(\Gamma))}\leq M \big\}.
 \label{Uad}
\end{align}
\end{definition}
\begin{remark}\label{adrem}
(1) Due to the compatibility condition $\mathbf{h}(x,t)|_{t=0}=\bfd_0(x)|_\Gamma$, the admissible set $\wt{\mathcal{U}}_{\mathrm{ad}}^M$
is non-empty, provided that $M>0$ is sufficiently large. For instance, if $\bfd_0|_\Gamma\in \mathbf{H}^\frac52(\Gamma)$, then $\mathbf{h}=\bfd_0|_\Gamma\in \wt{\mathcal{U}}_{\mathrm{ad}}^M$ provided  that $M\geq T^\frac12 \|\bfd_0|_\Gamma\|_{\mathbf{H}^\frac52(\Gamma)}$.

(2) For any fixed admissible $M>0$, the set $\wt{\mathcal{U}}_{\mathrm{ad}}^M$ is a bounded, convex
and closed subset of the Banach space $\mathcal{U}$.
\end{remark}

\subsection{The control-to-state operator $\mathcal{S}$}
Denote by
\begin{align}
\mathcal{H}:= & \big[C([0, T]; \mathbf{V})\cap L^2(0, T; \mathbf{H}^2(\Omega))\cap H^1(0,T; \mathbf{H})\big]\non\\
&\times  \big[C([0, T]; \mathbf{H}^2(\Omega)) \cap L^2(0, T; \mathbf{H}^3(\Omega))\cap H^1(0,T; \mathbf{H}^1(\Omega))\big]
\label{Hcal}
\end{align}
the function space for global strong solutions to problem \eqref{1}--\eqref{5} (cf. Theorem \ref{exe2d} and Remark \ref{rmexe2d}). Then we introduce

\begin{definition}\label{ctsdef}
The \emph{control-to-state} mapping $\mathcal{S}$ associated with the given initial data $(\mathbf{v}_0, \mathbf{d}_0)$ is defined as follows:
\begin{align}
\mathcal{S}:\wt{\mathcal{U}}\to \mathcal{H},\quad \mathbf{h}\in \wt{\mathcal{U}}\mapsto \mathcal{S}(\mathbf{h}):=(\bfv,\bfd)\in \mathcal{H},
\label{ctsS}
\end{align}
where $(\bfv,\bfd)$ is the unique global strong solution to problem \eqref{1}--\eqref{5} on $[0,T]$ subject to the initial data $(\mathbf{v}_0, \mathbf{d}_0)$.
\end{definition}

 The following property for $\mathcal{S}$ is a direct consequence of Theorem \ref{exe2d}, Remark \ref{rmexe2d}
 and the continuous dependence result of Proposition \ref{strconti}:
\bp[Lipschitz continuity of $\mathcal{S}$] \label{STC}
Let $n=2$ and $T\in(0,+\infty)$. Assume that the hypotheses of Theorem \ref{exe2d} are satisfied and $M>0$ is sufficiently large.
Then the control-to-state mapping $\mathcal{S}$ defined by \eqref{ctsS} is well-defined and bounded on $[0,T]$.
Moreover, the operator $\mathcal{S}$ is Lipschitz continuous from $\wt{\mathcal{U}}$ into the space
\begin{align}
\mathcal{W}:= & \big[C([0, T]; \mathbf{V})\cap L^2(0, T; \mathbf{H}^2(\Omega))\big]
\times  \big[C([0, T]; \mathbf{H}^2(\Omega)) \cap L^2(0, T; \mathbf{H}^3(\Omega))\big].
\label{Wcal}
\end{align}
\ep

\subsection{The existence of an optimal boundary control}

In what follows, we prove the existence of an optimal boundary control for problem (\textbf{CP}):
\bt[Existence of an optimal boundary control]\label{existence}
Let $n=2$. Suppose that the assumptions (A1)--(A2) are satisfied and $(\mathbf{v}_0, \mathbf{d}_0)\in \mathbf{V}\times
 \mathbf{H}^2(\Omega)$ is a given initial datum. Let $M>0$ be sufficiently large. Then the optimal control problem (\textbf{CP})
 admits a solution $((\bfv,\bfd), \mathbf{h})$ such that $\mathbf{h}\in \wt{\mathcal{U}}_{\mathrm{ad}}^M$
 and $(\bfv,\bfd)=\mathcal{S}(\mathbf{h})$ is the unique global strong solution to problem \eqref{1}--\eqref{5}
 with initial condition $(\mathbf{v}_0, \mathbf{d}_0)$.
\et
\begin{proof}
The proof essentially follows from the convexity/coercivity of the nonnegative cost functional $\mathcal{J}((\bfv, \bfd), \mathbf{h})$ (see \eqref{costJ}).
To this end, we introduce the reduced cost functional
$$\wt{\mathcal{J}}(\mathbf{h}):  \wt{\mathcal{U}}\to [0,+\infty)\quad \text{such that}\quad \wt{\mathcal{J}}(\mathbf{h}) :=\mathcal{J}((\bfv, \bfd),\mathbf{h})$$
 for any $\mathbf{h}\in \wt{\mathcal{U}}$, where $(\mathbf{v},\mathbf{d})=\mathcal{S}(\mathbf{h})$ is the unique global strong solution to problem
 \eqref{1}--\eqref{5} with initial datum $(\mathbf{v}_0, \mathbf{d}_0)$.
 It is easy to see that the optimal control problem (\textbf{CP}) is equivalent to the minimization problem (\textbf{CP})':
$$\min_{\mathbf{h}\in \wt{\mathcal{U}}^M_{\mathrm{ad}}} \wt{\mathcal{J}}(\mathbf{h}).$$
It follows from Remark \ref{adrem} that the admissible set $\wt{\mathcal{U}}_{\mathrm{ad}}^M$ is non-empty under our assumption. Then there exists a bounded minimizing sequence
$\{(\bfv^{(i)}, \bfd^{(i)}, \mathbf{h}^{(i)})\}$ $(i=1,2,...)$ of $\mathcal{J}$ such that $\mathbf{h}^{(i)}\in \wt{\mathcal{U}}_{\mathrm{ad}}^M$
and $(\bfv^{(i)},\bfd^{(i)})=\mathcal{S}(\mathbf{h}^{(i)})\in \mathcal{H}$ is the unique strong solution to problem \eqref{1}--\eqref{5}
subject to the initial condition $(\mathbf{v}_0, \mathbf{d}_0)$ and boundary condition $\mathbf{h}^{(i)}$. Moreover,
\begin{equation}
\lim_{i\to+\infty} \wt{\mathcal{J}}(\mathbf{h}^{(i)})=\inf_{\mathbf{h}\in \wt{\mathcal{U}}^M_{\mathrm{ad}}} \wt{\mathcal{J}}(\mathbf{h}).\label{infG}
\end{equation}

Since $\wt{\mathcal{U}}_{\mathrm{ad}}^M$ is a bounded subset of $\mathcal{U}$,
then there exists a weakly convergent subsequence of $\{\mathbf{h}^{(i)}\}$ satisfying the compatibility condition $\mathbf{h}^{(i)}|_{t=0}=\bfd_0|_\Gamma$,
which is still denoted by $\{\mathbf{h}^{(i)}\}$ without loss of generality, such that
\begin{align}
\mathbf{h}^{(i)} \rightharpoonup \mathbf{h}^\sharp\quad\text{in\ } \mathcal{U},\non
\end{align}
for some $\mathbf{h}^\sharp\in \wt{\mathcal{U}}$.
Besides, since $\wt{\mathcal{U}}_{\mathrm{ad}}^M$ is convex and closed, it is also weakly closed and we can infer that
$\mathbf{h}^\sharp\in \wt{\mathcal{U}}_{\mathrm{ad}}^M$.

For the weakly convergent subsequence $\{\mathbf{h}^{(i)}\}$, we consider $(\bfv^{(i)}, \bfd^{(i)})=\mathcal{S}(\mathbf{h}^{(i)})$.
From Proposition \ref{STC} and the definition of $\wt{\mathcal{U}}^M_{\mathrm{ad}}$ it follows that $\mathcal{S}(\mathbf{h}^{(i)})$ are
uniformly bounded (with respect to the index $i$) in the space $\mathcal{H}$. Hence, there exists a weakly convergent subsequence of $\{(\bfv^{(i)}, \bfd^{(i)})\}$,
which is still denoted by $\{(\bfv^{(i)}, \bfd^{(i)})\}$ without loss of generality, such that
\begin{align}
&\bfv^{(i)}\stackrel{*}\rightharpoonup \bfv^\sharp\quad \text{in\ } L^\infty(0, T; \mathbf{V}),\non\\
&\bfv^{(i)}\rightharpoonup \bfv^\sharp\quad \text{in\ } L^2(0, T; \mathbf{H}^2(\Omega))\cap H^1(0,T; \mathbf{H}),\non\\
&\bfd^{(i)}\stackrel{*}\rightharpoonup \bfd^\sharp\quad \text{in\ } L^\infty(0, T; \mathbf{H}^2(\Omega)),\non\\
&\bfd^{(i)}\rightharpoonup \bfd^\sharp\quad \text{in\ } L^2(0, T; \mathbf{H}^3(\Omega))\cap H^1(0,T; \mathbf{H}^1(\Omega)),\non
\end{align}
for some limit functions $(\bfv^\sharp, \bfd^\sharp)\in \mathcal{H}$.

Next, we proceed to verify that $(\bfv^\sharp, \bfd^\sharp)=\mathcal{S}(\mathbf{h}^\sharp)$.
By the Aubin--Lions lemma and the compact embedding theorems for $n=2$ (see, e.g., \cite{SI}),
we have the following strong convergence results for $\{(\bfv^{(i)}, \bfd^{(i)})\}$ (again up to a subsequence)
\begin{align}
 &\bfv^{(i)}\to \bfv^\sharp\quad \text{in\ } C([0, T]; \mathbf{H}^{1-s}(\Omega))\cap L^2(0,T; \mathbf{H}^{2-s}(\Omega)),\non\\
 &\bfd^{(i)}\to \bfd^\sharp\quad \text{in\ } C([0, T]; \mathbf{H}^{2-s}(\Omega))\cap L^2(0, T; \mathbf{H}^{3-s}(\Omega)),\non
\end{align}
for some $s\in (0,1)$. As a consequence, we also have $\bfv^{(i)}\to \bfv^\sharp$, $\bfd^{(i)}\to \bfd^\sharp$ a.e. in $\Omega\times[0,T]$
and the convergence for the initial data
\begin{align}
&\bfv^{(i)}(0)\to \bfv^\sharp(0)=\bfv_0\quad \text{in\ } \mathbf{H}^{1-s}(\Omega),\quad \bfd^{(i)}(0)\to \bfd^\sharp(0)=\bfd_0
\quad \text{in\ } \mathbf{H}^{2-s}(\Omega).\non
\end{align}
Using the above strong convergence results we are able to show the following convergence for nonlinear terms
\begin{align}
&\mathbf{v}^{(i)}\cdot\nabla \mathbf{v}^{(i)} \rightharpoonup \mathbf{v}^{\sharp}\cdot\nabla \mathbf{v}^{\sharp}
\qquad \text{in\ } L^2(0, T; \mathbf{L}^2(\Omega)),\non\\
&\nabla\cdot(\nabla \bfd^{(i)}\odot \nabla \bfd^{(i)})\rightharpoonup \nabla\cdot(\nabla \bfd^{\sharp}\odot \nabla \bfd^{\sharp})
\qquad \text{in\ } L^2(0, T; \mathbf{L}^2(\Omega)),\non\\
&\mathbf{v}^{(i)}\cdot\nabla \mathbf{d}^{(i)} \to \mathbf{v}^{\sharp}\cdot\nabla \mathbf{d}^{\sharp}
\qquad \text{in\ } L^2(0, T; \mathbf{L}^2(\Omega)),\non\\
&\mathbf{f}(\mathbf{d}^{(i)}) \to \mathbf{f}(\mathbf{d}^{\sharp}) \qquad \text{in\ } C(\overline{\Omega}\times [0,T]).\non
\end{align}
The above convergence results easily enable us  to pass to the limit (up to a subsequence) in the weak formulation of problem \eqref{1}--\eqref{5}
for every $(\bfv^{(i)},\bfd^{(i)})$ such that the limit function $(\bfv^\sharp, \bfd^\sharp)$ satisfies the weak formulation of problem \eqref{1}--\eqref{5}
with initial condition $(\mathbf{v}_0, \mathbf{d}_0)$ and boundary condition $\mathbf{h}^{\sharp}$. Since $(\bfv^\sharp, \bfd^\sharp)\in \mathcal{H}$,
we can conclude that $(\bfv^\sharp, \bfd^\sharp)=\mathcal{S}(\mathbf{h}^\sharp)$.

Therefore, the limit $((\bfv^\sharp, \bfd^\sharp), \mathbf{h}^\sharp)$ is admissible for the control problem (\textbf{CP}).
Since the cost functional  $\mathcal{J}$ is weakly lower semi-continuous in $\mathcal{H}\times \mathcal{U}$, by the weak convergent results
(up to a subsequence) obtained before,  it holds
$$  \liminf_{i\to+\infty}\mathcal{J}((\bfv^{(i)},\bfd^{(i)}), \mathbf{h}^{(i)})
\geq \mathcal{J}((\bfv^\sharp, \bfd^\sharp),\mathbf{h}^\sharp).$$
Recalling the definition of $\wt{\mathcal{J}}(\mathbf{h})$ and \eqref{infG}, we conclude
$$\mathcal{J}((\bfv^\sharp, \bfd^\sharp),\mathbf{h}^\sharp)= \min_{\mathbf{h}\in \wt{\mathcal{U}}^M_{\mathrm{ad}}} \mathcal{J}((\bfv, \bfd), \mathbf{h}),$$
which yields that
$((\bfv^\sharp, \bfd^\sharp), \mathbf{h}^\sharp)$ is a solution to the optimal control problem (\textbf{CP}).

The proof is complete.
\end{proof}

\section{Differentiability of the Control-to-State Operator $\mathcal{S}$}
\setcounter{equation}{0}

In this section, we aim to prove the differentiability of the control-to-state operator $\mathcal{S}$ with respect
to the boundary control $\mathbf{h}$ in $\wt{\mathcal{U}}$.

 \subsection{The linearized system}
Let $\mathbf{h}^*\in \wt{\mathcal{U}}$ be an arbitrary but given vector. We denote by $(\bfv^*, \bfd^*)=\mathcal{S}(\mathbf{h}^*)$
the associate unique global strong solution to the state system \eqref{1}--\eqref{5} given by Theorem \ref{exe2d}.

Below we investigate the linearization of the state system \eqref{1}--\eqref{5} around $((\bfv^*, \bfd^*),\mathbf{h}^*)$ for the unknowns denoted by
$(\bom, \bphi)$ with an arbitrary given vector $\bxi\in \wt{\mathcal{U}}-\{\mathbf{h}^*\}$ (i.e., a perturbation with respect to $\mathbf{h}^*$), that is
\begin{align}
 &\p_t\bom-\Delta \bom +\nabla \hat{P}+(\bfv^*\cdot\nabla)\bom + (\bom\cdot\nabla) \mathbf{v}^*\non\\
 &\qquad =-\nabla\cdot(\nabla \bphi\odot\nabla \mathbf{d}^*)-\nabla\cdot(\nabla \mathbf{d}^*\odot\nabla \bphi), & \text{in\ } \Omega\times (0,T),\label{L1}\\
 &\nabla \cdot \bom = 0,&\text{in\ } \Omega\times (0,T),\label{L2}\\
 &\p_t \bphi-\Delta \bphi+(\bfv^*\cdot\nabla)\bphi+(\bom\cdot\nabla) \bfd^*=-\mathbf{f}'(\mathbf{d}^*)\bphi,&\text{in\ } \Omega\times (0,T),\label{L3}
 \end{align}
 subject to the following boundary and initial conditions
 \begin{align}
 & \bom=\mathbf{0},\quad \bphi=\bxi,\quad \text{on\ } \Gamma\times(0,T),\label{L4}\\
 & \bom|_{t=0}=\mathbf{0},\quad \bphi|_{t=0}=\mathbf{0},\quad \text{in\ }\Omega.\label{L5}
 \end{align}
In \eqref{L3} we have denoted $$\mathbf{f}'(\mathbf{d}^*)\bphi:=2(\mathbf{d}^*\cdot\bphi)\mathbf{d}^*+|\mathbf{d}^*|^2\bphi-\bphi.$$
\bp\label{Lin}
Let $n=2$. For any  $\bxi\in \wt{\mathcal{U}}-\{\mathbf{h}^*\}$, problem \eqref{L1}--\eqref{L5} admits a unique weak solution $(\bom, \bphi)$ such that
\begin{align}
&\bom\in C([0, T]; \mathbf{H})\cap L^2(0, T; \mathbf{V})\cap H^1(0,T; \mathbf{V}'),\label{Linreg1}\\
&\bphi\in C([0, T]; \mathbf{H}^1(\Omega)) \cap L^2(0, T; \mathbf{H}^2(\Omega))\cap H^1(0,T; \mathbf{L}^2(\Omega))\label{Linreg2}.
\end{align}
\ep
\begin{proof}
Let $\{\mathbf{u}_i\}_{i=1}^{\infty}$, $\{\bpsi_i\}_{i=1}^\infty$ be the basis of the Hilbert spaces $\mathbf{V}$ and $\mathbf{H}_0^1(\Omega)$,
respectively, which are given by the eigenfunctions of the Stokes problem
$$ (\nabla \mathbf{u}_i, \nabla \mathbf{w})=\lambda_i(\mathbf{u}_i, \mathbf{w}),\quad \forall\, \mathbf{w}\in \mathbf{V},
\ \text{with}\ \|\mathbf{u}_i\|_{\mathbf{L}^2}=1,$$
and the Laplace equation
$$ (\nabla \bpsi_i, \nabla \bm{\zeta})=\mu_i(\bpsi_i, \bm{\zeta}),\quad \forall\, \bm{\zeta}\in \mathbf{H}_0^1(\Omega),
\ \text{with}\ \|\bpsi_i\|_{\mathbf{L}^2}=1.$$
Here $\lambda_i$ is the eigenvalue corresponding to $\mathbf{u}_i$ and $\mu_i$ is the eigenvalue corresponding to $\bpsi_i$.
It is well-known that $0<\lambda_1<\lambda_2<...\nearrow +\infty$
is a monotone increasing sequence, $\{\mathbf{u}_i\}_{i=1}^{\infty}$ forms a complete orthonormal basis in $\mathbf{H}$ and it is
also orthogonal in $\mathbf{V}$ (see e.g., \cite{Te1}). Similarly,  $0<\mu_1<\mu_2<...\nearrow +\infty$ is a monotone increasing sequence,
$\{\bpsi_i\}_{i=1}^{\infty}$ forms a complete orthonormal basis in $\mathbf{L}^2(\Omega)$ and it is also orthogonal in $\mathbf{H}^1_0(\Omega)$.
By the classical elliptic regularity theory, we have $\mathbf{u}_i, \bpsi_i\in C^\infty(\overline{\Omega})$ for all $i\in \mathbb{N}$.
For every $m\in\mathbb{N}$, we denote by $\mathbf{V}_m=\mathrm{span}\{\mathbf{u}_1,...,\mathbf{u}_m\}$ and
$\mathbf{W}_m=\mathrm{span}\{\bpsi_1,...,\bpsi_m\}$ the finite dimensional subspaces of $\mathbf{V}$ and $\mathbf{H}^1_0(\Omega)$ spanned
by their first $m$ basis functions, respectively. Moreover, we use $\Pi_m$ for the orthogonal projection from $\mathbf{H}$ onto $\mathbf{V}_m$
and $\widetilde{\Pi}_m$ for the orthogonal projection from $\mathbf{L}^2(\Omega)$ onto $\mathbf{W}_m$.

Given a vector $\bxi\in \widetilde{\mathcal{U}}-\{\mathbf{h}^*\}$, we denote by $\bphi_E$ the unique solution to the elliptic problem
\be
 \begin{cases}
 -\Delta \bphi_E=\mathbf{0},\qquad\ \ \ \ \text{ in } \Omega\times (0,T),\\
 \bphi_E=\bxi(x,t),\qquad \ \ \text{ on } \Gamma\times (0,T).
 \end{cases}\label{LEphi}
 \ee
 It is easy to see that $\bphi_E|_{t=0}=\mathbf{0}$ according to the definition of $\bxi$.
 Then, for every integer $m\geq 1$, we look for solutions of the form:
\begin{align}
& \bom^m=\sum_{i=1}^m a_i^m(t)\mathbf{u}_i(x), \quad \bphi^m=\bphi_E+\wha{\bphi}^m =\bphi_E+\sum_{i=1}^m b_i^m(t)\bpsi_i(x),\non
\end{align}
solving the following approximate problem of \eqref{L1}--\eqref{L5}, a.e. in $[0,T]$ and for $i=1,2,...,m$:

\begin{align}
 &\langle \partial_t \bom^m, \mathbf{u}_i\rangle_{\mathbf{V}', \mathbf{V}}+ \int_\Omega \nabla \bom^m : \nabla \mathbf{u}_i dx
           + \int_\Omega [(\bfv^*\cdot\nabla)\bom^m + (\bom^m\cdot\nabla) \mathbf{v}^*]\cdot \mathbf{u}_i dx\non\\
 &\qquad =\int_\Omega (\nabla \wha{\bphi}^m \odot\nabla \mathbf{d}^*): \nabla \mathbf{u}_i dx
 + \int_\Omega (\nabla \mathbf{d}^*\odot\nabla \wha{\bphi}^m):\nabla \mathbf{u}_i dx\non \\
  &\qquad \quad + \int_\Omega (\nabla \bphi_E \odot\nabla \mathbf{d}^*): \nabla \mathbf{u}_i dx
  + \int_\Omega (\nabla \mathbf{d}^*\odot\nabla \bphi_E):\nabla \mathbf{u}_i dx, \label{APL1}\\
 &\int_\Omega \p_t \wha{\bphi}^m \cdot \bpsi_i dx+\int_\Omega \nabla \wha{\bphi}^m :\nabla \bpsi_i dx
           +\int_\Omega [(\bfv^*\cdot\nabla)\wha{\bphi}^m+(\bom^m\cdot\nabla) \bfd^*]\cdot \bpsi_i dx\non\\
 &\qquad \qquad +\int_\Omega \mathbf{f}'(\mathbf{d}^*)\wha{\bphi}^m \cdot \bpsi_i dx\non\\
 &\qquad = -\int_\Omega \mathbf{f}'(\mathbf{d}^*)\bphi_E \cdot \bpsi_i dx -\int_\Omega \p_t \bphi_E\cdot \bpsi_i dx
 -\int_\Omega (\bfv^*\cdot\nabla)\bphi_E\cdot \bpsi_i dx,\label{APL3}\\
 & \bom^m|_{t=0}=\mathbf{0}, \quad \wha{\bphi}^m|_{t=0}=\mathbf{0},\quad \text{in}\ \Omega.\label{APLini}
 \end{align}
  From the fact $(\bfv^*, \bfd^*)=\mathcal{S}(\mathbf{h}^*)\in \mathcal{H}$ and applying Lemma \ref{Ap1} to $\bphi_E$, we infer that
  \eqref{APL1}--\eqref{APLini} is indeed a Cauchy problem of a linear ODE system for  $\mathbf{a}^m(t)=(a_1^m (t),...,a_m^m(t))^{\mathrm{tr}}$
  and $\mathbf{b}^m(t)=(b_1^m (t),...,b_m^m(t))^{\mathrm{tr}}$ with all the coefficients belonging to $L^2(0,T)$.
  Then it is standard to conclude that the approximate problem \eqref{APL1}--\eqref{APLini} admits a unique solution
  $(\mathbf{a}^m(t), \mathbf{b}^m(t))\in H^1(0,T;\mathbb{R}^{2m})$.

Next, we derive some \textit{a priori} estimates for the approximate solutions $(\bom^m, \bphi^m)$ that are uniform with respect to the parameter $m$.
For $i=1,...,m$, multiplying \eqref{APL1} by $a_i^m(t)$ and \eqref{APL3} by $\mu_i b_i^m(t)$, respectively, and then adding the resultants together we get
\begin{align}
&\frac12 \frac{d}{dt}\big(\|\bom^m\|_{\mathbf{L}^2}^2+\|\nabla \wha{\bphi}^m\|_{\mathbf{L}^2}^2\big)+\|\nabla \bom^m\|_{\mathbf{L}^2}^2
+\|\Delta \wha{\bphi}^m\|_{\mathbf{L}^2}^2\non\\
&\quad = - \int_\Omega (\bom^m\cdot\nabla) \mathbf{v}^*\cdot \bom^m dx
         + \int_\Omega (\bom^m\cdot\nabla) \bfd^* \cdot \Delta \wha{\bphi}^m dx\non\\
&\quad\ \quad +\int_\Omega (\bfv^*\cdot\nabla)\wha{\bphi}^m\cdot \Delta \wha{\bphi}^m dx
              +\int_\Omega \mathbf{f}'(\mathbf{d}^*)\wha{\bphi}^m \cdot \Delta \wha{\bphi}^m dx\non\\
  &\quad\ \quad +\int_\Omega (\nabla \wha{\bphi}^m \odot\nabla \mathbf{d}^*): \nabla \bom^m  dx
                + \int_\Omega (\nabla \mathbf{d}^*\odot\nabla \wha{\bphi}^m):\nabla \bom^m  dx\non \\
  &\quad\ \quad + \int_\Omega (\nabla \bphi_E \odot\nabla \mathbf{d}^*): \nabla \bom^m  dx
                + \int_\Omega (\nabla \mathbf{d}^*\odot\nabla \bphi_E):\nabla \bom^m  dx\non\\
  &\quad \ \quad +\int_\Omega [\mathbf{f}'(\mathbf{d}^*)\bphi_E +\p_t \bphi_E+ (\bfv^*\cdot\nabla)\bphi_E]\cdot \Delta \wha{\bphi}^m  dx\non\\
  &:=\sum_{k=1}^9 J_k,\label{Lines}
\end{align}
where we have used the conditions $\nabla \cdot \bfv^*=0$ and $\wha{\bphi}^m|_{\Gamma}=\Delta \wha{\bphi}^m|_{\Gamma}=\mathbf{0}$.

For the sake of simplicity, we shall denote by $C$ the constants that may depend on the global strong solution $(\bfv^*, \bfd^*)\in \mathcal{H}$
(cf. \eqref{globstrong1}--\eqref{globstrong3}), $\Omega$ and $T$, but not on $m$.

Now, we estimate the right-hand side of \eqref{Lines} term by term. Using the H\"older inequality, Young's inequality and the estimates
\eqref{globstrong1}--\eqref{globstrong3} for $(\bfv^*, \bfd^*)$, we obtain
\begin{align}
J_1&\leq \|\bom^m\|_{\mathbf{L}^4}\|\nabla \mathbf{v}^*\|_{\mathbf{L}^4}\|\bom^m\|_{\mathbf{L}^2}\non\\
   &\leq C\|\mathbf{v}^*\|_{\mathbf{H}^2}\|\nabla \bom^m\|_{\mathbf{L}^2}\|\bom^m\|_{\mathbf{L}^2}\non\\
   &\leq \varepsilon \|\nabla \bom^m\|_{\mathbf{L}^2}^2+ C\varepsilon^{-1}\|\mathbf{v}^*\|^2_{\mathbf{H}^2}\|\bom^m\|_{\mathbf{L}^2}^2,\non
\end{align}
\begin{align}
J_2&\leq \|\bom^m\|_{\mathbf{L}^4}\|\nabla \bfd^*\|_{\mathbf{L}^4}\|\Delta \wha{\bphi}^m\|_{\mathbf{L}^2}\non\\
   &\leq C\|\bfd^*\|_{\mathbf{H}^2}\|\bom^m\|_{\mathbf{L}^2}^\frac12\|\nabla \bom^m\|_{\mathbf{L}^2}^\frac12\|\Delta \wha{\bphi}^m\|_{\mathbf{L}^2}\non\\
   &\leq \varepsilon \|\nabla \bom^m\|_{\mathbf{L}^2}^2+ \varepsilon \|\Delta \wha{\bphi}^m\|_{\mathbf{L}^2}^2
   +C\varepsilon^{-3} \|\bfd^*\|_{\mathbf{H}^2}^4\|\bom^m\|_{\mathbf{L}^2}^2,\non
\end{align}
\begin{align}
J_3&\leq \|\bfv^*\|_{\mathbf{L}^\infty}\|\nabla\wha{\bphi}^m\|_{\mathbf{L}^2}\|\Delta \wha{\bphi}^m\|_{\mathbf{L}^2}\non\\
   &\leq \varepsilon\|\Delta \wha{\bphi}^m\|_{\mathbf{L}^2}^2+ C \varepsilon^{-1} \|\bfv^*\|_{\mathbf{H}^2}^2\|\nabla\wha{\bphi}^m\|_{\mathbf{L}^2}^2,\non
\end{align}
\begin{align}
J_4&\leq \|\mathbf{f}'(\mathbf{d}^*)\|_{\mathbf{L}^\infty}\|\wha{\bphi}^m\|_{\mathbf{L}^2} \|\Delta \wha{\bphi}^m\|_{\mathbf{L}^2}\non\\
   &\leq \varepsilon \|\Delta \wha{\bphi}^m\|_{\mathbf{L}^2}^2+ C \varepsilon^{-1}(\|\bfd^*\|_{\mathbf{H}^2}^4+1)\|\nabla \wha{\bphi}^m\|_{\mathbf{L}^2}^2,\non
\end{align}
\begin{align}
J_5+J_6&\leq \|\nabla \wha{\bphi}^m\|_{\mathbf{L}^4}\|\nabla \mathbf{d}^*\|_{\mathbf{L}^4}\|\nabla \bom^m\|_{\mathbf{L}^2}\non\\
       &\leq \varepsilon \|\nabla \bom^m\|_{\mathbf{L}^2}^2+ \varepsilon \|\Delta \wha{\bphi}^m\|_{\mathbf{L}^2}^2
             + C\varepsilon^{-2} \|\mathbf{d}^*\|_{\mathbf{H}^2}^2 \|\nabla \wha{\bphi}^m\|_{\mathbf{L}^2}^2,\non
\end{align}
\begin{align}
J_7+J_8&\leq \|\nabla \bphi_E\|_{\mathbf{L}^4} \|\nabla \mathbf{d}^*\|_{\mathbf{L}^4}\| \nabla \bom^m \|_{\mathbf{L}^2}\non\\
       &\leq   \varepsilon \|\nabla \bom^m\|_{\mathbf{L}^2}^2+C\varepsilon^{-1} \|\bphi_E\|_{\mathbf{H}^2}^2 \|\mathbf{d}^*\|_{\mathbf{H}^2}^2\non\\
       &\leq   \varepsilon \|\nabla \bom^m\|_{\mathbf{L}^2}^2+C\varepsilon^{-1} \|\bxi\|_{\mathbf{H}^\frac{3}{2}(\Gamma)}^2 \|\mathbf{d}^*\|_{\mathbf{H}^2}^2,\non
\end{align}
\begin{align}
J_9&\leq \big(\|\mathbf{f}'(\mathbf{d}^*)\|_{\mathbf{L}^\infty}\|\bphi_E\|_{\mathbf{L}^2} + \|\p_t \bphi_E\|_{\mathbf{L}^2}
+ \|\bfv^*\|_{\mathbf{L}^4}\|\nabla\bphi_E\|_{\mathbf{L}^4}\big)
          \|\Delta \wha{\bphi}^m\|_{\mathbf{L}^2}\non\\
   &\leq \varepsilon \|\Delta \wha{\bphi}^m\|_{\mathbf{L}^2}^2+ C\varepsilon^{-1}(\|\mathbf{d}^*\|_{\mathbf{L}^\infty}^2+1)\|\bphi_E\|_{\mathbf{L}^2}^2
   +C\varepsilon^{-1}\|\p_t \bphi_E\|_{\mathbf{L}^2}^2\non\\
   &\qquad  + C\varepsilon^{-1}  \|\bfv^*\|_{\mathbf{H}^1}^2\|\bphi_E\|_{\mathbf{H}^2}^2\non\\
   &\leq \varepsilon \|\Delta \wha{\bphi}^m\|_{\mathbf{L}^2}^2+ C\varepsilon^{-1}(\|\mathbf{d}^*\|_{\mathbf{H}^2}^2+1)\|\bxi\|_{\mathbf{H}^\frac{3}{2}(\Gamma)}^2
            +C\varepsilon^{-1}\|\p_t \bxi\|_{\mathbf{H}^\frac12(\Gamma)}^2\non\\
   &\qquad  + C\varepsilon^{-1}  \|\bfv^*\|_{\mathbf{H}^1}^2\|\bxi\|_{\mathbf{H}^\frac{3}{2}(\Gamma)}^2.\non
\end{align}
Taking $\varepsilon$ sufficiently small, from \eqref{Lines} and the estimates \eqref{globstrong1} for
$\|\bfv^*\|_{L^\infty(0,T;\mathbf{H}^1)}$, $\|\bfd^*\|_{L^\infty(0,T;\mathbf{H}^2)}$  we infer from the above estimates that
\begin{align}
& \frac{d}{dt}\big(\|\bom^m\|_{\mathbf{L}^2}^2+\|\nabla \wha{\bphi}^m\|_{\mathbf{L}^2}^2\big)+\|\nabla \bom^m\|_{\mathbf{L}^2}^2
+\|\Delta \wha{\bphi}^m\|_{\mathbf{L}^2}^2\non\\
&\quad \leq C(1+\|\mathbf{v}^*\|^2_{\mathbf{H}^2})\big(\|\bom^m\|_{\mathbf{L}^2}^2+\|\nabla \wha{\bphi}^m\|_{\mathbf{L}^2}^2\big)
            +C\big( \|\bxi\|_{\mathbf{H}^\frac{3}{2}(\Gamma)}^2+\|\p_t \bxi\|_{\mathbf{H}^\frac12(\Gamma)}^2\big).\label{Lines1}
\end{align}
Besides, it follows from \eqref{globstrong2} that $\|\mathbf{v}^*\|_{L^2(0,T;\mathbf{H}^2)}$ is bounded.
Then, by Gronwall's lemma and the condition $\|\bom^m(0)\|_{\mathbf{L}^2}^2+\|\nabla \wha{\bphi}^m(0)\|_{\mathbf{L}^2}^2=0$, we deduce, for any $m\geq 1$,
\begin{align}
&\|\bom^m(t)\|_{\mathbf{L}^2}^2+\|\nabla \wha{\bphi}^m(t)\|_{\mathbf{L}^2}^2+\int_0^t\big(\|\nabla \bom^m(\tau)\|_{\mathbf{L}^2}^2
+\|\Delta \wha{\bphi}^m(\tau)\|_{\mathbf{L}^2}^2\big)d\tau \non\\
&\quad \leq C\int_0^t\big( \|\bxi(\tau)\|_{\mathbf{H}^\frac{3}{2}(\Gamma)}^2+\|\p_t \bxi(\tau)\|_{\mathbf{H}^\frac12(\Gamma)}^2\big)d\tau\non\\
&\quad \leq C\|\bxi\|_{\mathcal{U}}^2,\quad \forall\, t \in [0,T].\label{Lines2}
\end{align}
Furthermore, by means of the linear equations \eqref{APL1}, \eqref{APL3}, we obtain the following uniform estimate for the time derivatives of
$\bom^m$ and $\wha{\bphi}^m$
\begin{align}
\|\p_t\bom^m\|_{L^2(0,T;\mathbf{V}')}+\|\p_t\wha{\bphi}^m\|_{L^2(0,T;\mathbf{L}^2)}\leq C\|\bxi\|_{\mathcal{U}}.\label{Lines3}
\end{align}

As a consequence, from \eqref{Lines2} and \eqref{Lines3} it follows that there exists a pair $(\bom, \wha{\bphi})$ satisfying
\begin{align}
&\bom\in L^\infty(0, T; \mathbf{H})\cap L^2(0, T; \mathbf{V})\cap H^1(0,T; \mathbf{V}'),\non\\
&\wha{\bphi}\in L^\infty(0, T; \mathbf{H}^1(\Omega)) \cap L^2(0, T; \mathbf{H}^2(\Omega))\cap H^1(0,T; \mathbf{L}^2(\Omega)),\non
\end{align}
which is the weak (or weak star) limit of convergent subsequences $\{\bom^m\}$, $\{\wha{\bphi}^m\}$ (not relabelled for simplicity) in the
corresponding spaces in \eqref{Linreg1}, \eqref{Linreg2}, as $m\to+\infty$.
Then, by a standard compactness argument, we are able to pass to the limit as $m\to+\infty$ (up to a subsequence) in the approximate system
\eqref{APL1}, \eqref{APL3} and to verify that the pair $(\bom, \bphi)$ with $\bphi=\wha{\bphi}+\bphi_E$ is indeed a weak solution to the
linearized problem \eqref{L1}--\eqref{L5} satisfying \eqref{Linreg1}, \eqref{Linreg2}.
Here, we also use the elliptic estimates for $\bphi_E$ (cf. Lemma \ref{Ap1}) and the assumption $\bxi\in \wt{\mathcal{U}}-\{\mathbf{h}^*\}$.
Besides, by the interpolation theorem (cf. \cite{SI}), we conclude that
$\bom\in C([0, T]; \mathbf{H})$, $\bphi\in C([0, T]; \mathbf{H}^1(\Omega))$. Finally, the uniqueness of weak solutions to the linearized problem
\eqref{L1}--\eqref{L5} follows from the standard energy method, which is omitted here.

The proof is complete.
\end{proof}

\noindent Define the lower-order function space (compare with $\mathcal{H}$ given by \eqref{Hcal})
\begin{align}
\mathcal{H}_1:= & \big[C([0, T]; \mathbf{H})\cap L^2(0, T; \mathbf{V})\cap H^1(0,T; \mathbf{V}')\big]\non\\
&\times  \big[C([0, T]; \mathbf{H}^1(\Omega)) \cap L^2(0, T; \mathbf{H}^2(\Omega))\cap H^1(0,T; \mathbf{L}^2(\Omega))\big].
\label{H1cal}
\end{align}
We introduce the linear mapping $\mathcal{L}_{\mathbf{h}^*}$ associated with the given vector $\mathbf{h}^*$ as well as $(\bfv^*, \bfd^*)=\mathcal{S}(\mathbf{h}^*)$,
such that
\begin{align}
\mathcal{L}_{\mathbf{h}^*}:\wt{\mathcal{U}}-\{\mathbf{h}^*\} \to \mathcal{H}_1,
\quad \bxi \in \wt{\mathcal{U}} -\{\mathbf{h}^*\} \mapsto \mathcal{L}_{\mathbf{h}^*}(\bxi):=(\bom,\bphi)\in \mathcal{H}_1.
\label{LinO}
\end{align}
In particular, in the definition \eqref{LinO},  $(\bom,\bphi)$ is the unique global weak solution to the linearized problem \eqref{L1}--\eqref{L5} on $[0,T]$ given
by Proposition \ref{Lin}.

For the linear parabolic problem \eqref{L1}--\eqref{L5}, it follows from the simple energy estimates like \eqref{Lines2}, \eqref{Lines3} that
\bc
The linear mapping $\mathcal{L}_{\mathbf{h}^*}(\bxi)=(\bom,\bphi)$ is a (Lipschitz) continuous mapping from $\wt{\mathcal{U}}-\{\mathbf{h}^*\}$ to $\mathcal{H}_1$.
\ec

\subsection{Differentiability of  $\mathcal{S}$}

First, let us introduce the precise definition of differentiability for the control-to-state operator $\mathcal{S}$:

\bd\label{diff}
Let
\begin{align}
\mathcal{W}_1:= & \big[C([0, T]; \mathbf{H})\cap L^2(0, T; \mathbf{H}^1(\Omega))\big]
\times  \big[C([0, T]; \mathbf{H}^1(\Omega)) \cap L^2(0, T; \mathbf{H}^2(\Omega))\big]
\label{W1cal}
\end{align}
and $\wt{\mathcal{U}}$ and $\mathcal{H}$ be defined as in \eqref{U} and \eqref{Hcal}, respectively.
Consider the control-to-state operator $\mathcal{S}: \wt{\mathcal{U}} \to \mathcal{H}$ associated with a given initial datum
$(\mathbf{v}_0, \mathbf{d}_0)\in \mathbf{V}\times \mathbf{H}^2(\Omega)$ (see Definition \ref{ctsdef}).
In particular, here we view $\mathcal{S}$ as a mapping from $\wt{\mathcal{U}}$ to the weaker space $\mathcal{W}_1$.
We say that $\mathcal{S}: \wt{\mathcal{U}}\to \mathcal{W}_1$ is Fr\'echet differentiable in $\wt{\mathcal{U}}$ if, for any
$\mathbf{h}^*\in \wt{\mathcal{U}}$, there exists a linear operator denoted by $\mathcal{S}'(\mathbf{h}^*): \wt{\mathcal{U}} -\{\mathbf{h}^*\} \to \mathcal{W}_1$
such that
\begin{align}
\lim_{\|\bxi\|_{\mathcal{U}}\to 0}
\frac{\|\mathcal{S}(\mathbf{h}^*+\bxi)-\mathcal{S}(\mathbf{h}^*)-\mathcal{S}'(\mathbf{h}^*)(\bxi)\|_{\mathcal{W}_1}}{\|\bxi\|_{\mathcal{U}}}=0,\label{ddiff}
\end{align}
where $\bxi \in \wt{\mathcal{U}} -\{\mathbf{h}^*\}$ is an arbitrary (small) perturbation of $\mathbf{h}^*$.
\ed
Then we can prove the following result:
\bt\label{SD}
Let $n=2$. Suppose that $(\mathbf{v}_0, \mathbf{d}_0)\in \mathbf{V}\times \mathbf{H}^2(\Omega)$ is a given initial datum.
Then the control-to-state operator $\mathcal{S}$ is Fr\'echet differentiable in the sense of Definition \ref{diff}.
 Moreover, for any $\mathbf{h}^*\in \wt{\mathcal{U}}$, its Fr\'echet derivative $\mathcal{S}'(\mathbf{h}^*)$ is given by
 \begin{align}
 \mathcal{S}'(\mathbf{h}^*)\bxi=(\bom,\bphi),\quad \forall\, \bxi\in \wt{\mathcal{U}} -\{\mathbf{h}^*\},
 \end{align}
where $(\bom,\bphi)$ is the unique global weak solution to the linearized problem \eqref{L1}--\eqref{L5} obtained in Proposition \ref{Lin}.
Namely, we have $\mathcal{S}'(\mathbf{h}^*)=\mathcal{L}_{\mathbf{h}^*}$, where $\mathcal{L}_{\mathbf{h}^*}$ is the linear operator defined in \eqref{LinO}.
\et
\begin{proof}
Let $\mathbf{h}^*\in \wt{\mathcal{U}}$ be an arbitrary but fixed vector. We consider any of its perturbation $\mathbf{h}=\mathbf{h}^*+\bxi$ with
$\bxi\in \wt{\mathcal{U}} -\{\mathbf{h}^*\}$ such that
\begin{align}
\|\bxi\|_{\mathcal{U}}\leq \kappa,\label{kappa}
\end{align}
for some fixed constant $\kappa>0$. For the given initial datum $(\mathbf{v}_0, \mathbf{d}_0)\in \mathbf{V}\times \mathbf{H}^2(\Omega)$, we denote
\begin{align}
&(\bfv, \bfd)=\mathcal{S}(\mathbf{h}),\quad (\bfv^*, \bfd^*)=\mathcal{S}(\mathbf{h}^*),\quad (\bom,\bphi)=\mathcal{L}_{\mathbf{h}^*}(\bxi).\non
\end{align}

It follows from Theorem \ref{exe2d}, Remark \ref{rmexe2d} and Proposition \ref{Lin} that the following uniform estimates hold for
$(\bfv, \bfd)$, $(\bfv^*, \bfd^*)$ and $(\bom,\bphi)$, respectively:
\begin{align}
\|(\bfv, \bfd)\|_{\mathcal{H}}\leq K_1,\quad \|(\bfv^*, \bfd^*)\|_{\mathcal{H}}\leq K_2,
\quad \|(\bom, \bphi)\|_{\mathcal{H}_1}\leq K_3\|\bxi\|_{\mathcal{U}},\label{UES}
\end{align}
where $K_1$, $K_2$, $K_3$ are positive constants that may depend on $\|\mathbf{v}_0\|_{\mathbf{H}^1}$,
$\|\mathbf{d}_0\|_{\mathbf{H}^2}$, $\|\mathbf{h}^*\|_{L^2(0,T;\mathbf{H}^\frac52(\Gamma))}$,
$\|\partial_t\mathbf{h}^*\|_{L^4(0,T; \mathbf{H}^\frac12(\Gamma))}$, $\Omega$ and $T$.
The constant $K_1$ also depends on $\kappa$ but it is independent of $\bxi$.

Besides, from the continuous dependence result for global strong solutions to problem \eqref{1}--\eqref{5} (see Proposition \ref{strconti}), we infer that
    \begin{align}
&\|\bfv-\bfv^*\|_{C([0,T];\mathbf{H}^1)}^2+\|\bfd-\bfd^*\|_{C([0,T];\mathbf{H}^2)}^2
+ \|\bfv-\bfv^*\|_{L^2(0,T;\mathbf{H}^2)}^2+\|\bfd-\bfd^*\|_{L^2(0,T;\mathbf{H}^3)}^2 \non\\
&\quad \leq K_4\|\bxi\|_{\mathcal{U}}^2,\label{contixi}
\end{align}
where $K_4>0$ is a constant depending on $\|\mathbf{v}_0\|_{\mathbf{H}^1}$,
$\|\mathbf{d}_0\|_{\mathbf{H}^2}$,
$\|\mathbf{h}^*\|_{L^2(0,T;\mathbf{H}^\frac52(\Gamma))}$, $\|\partial_t\mathbf{h}^*\|_{L^4(0,T; \mathbf{H}^\frac12(\Gamma))}$, $\kappa$, $\Omega$ and $T$.

Set the difference functions
\begin{align}
&(\mathbf{w},\bfe):=(\bfv-\bfv^*-\bom,\ \bfd-\bfd^*-\bphi).\non
\end{align}
We can easily see that $(\mathbf{w}, \bfe)\in \mathcal{H}_1$, i.e.,
\begin{align}
\mathbf{w} & \in  C([0, T]; \mathbf{H})\cap L^2(0, T; \mathbf{V})\cap H^1(0,T; \mathbf{V}'),\non\\
\bfe &\in  C([0, T]; \mathbf{H}^1(\Omega)) \cap L^2(0, T; \mathbf{H}^2(\Omega))\cap H^1(0,T; \mathbf{L}^2(\Omega)),\non
\end{align}
and $(\mathbf{w}, \bfe)$ turns out to be a weak solution to the following system in $ \Omega\times (0,T)$
\begin{align}
 &\p_t \mathbf{w}-\Delta \mathbf{w} +\nabla \tilde{P}+[(\bfv-\bfv^*)\cdot\nabla](\bfv-\bfv^*)+ (\bfv^*\cdot\nabla)\bfw + (\bfw\cdot\nabla) \mathbf{v}^*\non\\
 &\qquad =-\nabla \cdot[\nabla (\bfd-\bfd^*)\odot\nabla(\bfd-\bfd^*)] -\nabla\cdot(\nabla \mathbf{d}^*\odot\nabla \bfe)
 -\nabla\cdot(\nabla \bfe\odot\nabla \mathbf{d}^*), \label{D1}\\
 &\nabla \cdot \mathbf{w} = 0,\label{D2}\\
 &\p_t\bfe-\Delta \bfe+[(\bfv-\bfv^*)\cdot\nabla](\bfd-\bfd^*)+(\bfv^*\cdot\nabla)\bfe+(\bfw\cdot\nabla) \bfd^*\non\\
 &\qquad =-\mathbf{f}(\bfd)+\mathbf{f}(\bfd^*)+\mathbf{f}'(\mathbf{d}^*)\bphi,\label{D3}
 \end{align}
 subject to the homogeneous boundary and initial conditions
 \begin{align}
 & \bfw=\mathbf{0},\quad \bfe=\mathbf{0},\quad \text{on\ } \Gamma\times(0,T),\label{D4}\\
 & \bfw|_{t=0}=\mathbf{0},\quad \bfe|_{t=0}=\mathbf{0},\quad \text{in\ }\Omega.\label{D5}
 \end{align}
In \eqref{D1}, the pressure is given by $\tilde{P}=P-P^*-\hat{P}$, with $P$, $P^*$ and $\hat{P}$ being associated with $(\bfv, \bfd)$, $(\bfv^*, \bfd^*)$
and $(\bom, \bphi)$, respectively.
Besides, the right-hand side of \eqref{D3} can be re-written into the following explicit form:
\begin{align}
&-\mathbf{f}(\bfd)+\mathbf{f}(\bfd^*)+\mathbf{f}'(\mathbf{d}^*)\bphi\non\\
&\quad = -|\bfd-\bfd^*|^2(\bfd-\bfd^*)-|\bfd-\bfd^*|^2\bfd^*-2[(\bfd-\bfd^*)\cdot\bfd^*](\bfd-\bfd^*)\non\\
&\qquad \ -2(\bfd^*\cdot \bfe) \bfd^*-|\bfd^*|^2\bfe+\bfe.\label{dfaa}
\end{align}

Testing \eqref{D1} by $\bfw$, integrating over $\Omega$, and using the incompressibility condition for $\bfv^*$ as well as $\bfw$,
after integration by parts, we deduce
\begin{align}
&\frac12 \frac{d}{dt} \|\bfw\|_{\mathbf{L}^2}^2+\|\nabla \bfw \|_{\mathbf{L}^2}^2\non\\
&\quad = -\int_\Omega [(\bfv-\bfv^*)\cdot\nabla](\bfv-\bfv^*) \cdot \bfw dx
   -\int_\Omega  (\bfw\cdot\nabla) \mathbf{v}^*\cdot \bfw dx\non\\
&\qquad \ +\int_\Omega [\nabla (\bfd-\bfd^*)\odot\nabla(\bfd-\bfd^*)]: \nabla \bfw dx
          + \int_\Omega (\nabla \mathbf{d}^*\odot\nabla \bfe+\nabla \bfe\odot\nabla \mathbf{d}^*):\nabla\bfw dx\non\\
&\quad := \sum_{j=1}^4 r_j. \label{Dres1}
\end{align}
In a similar manner, testing \eqref{D3} by $-\Delta \bfe$, we obtain
\begin{align}
&\frac12\frac{d}{dt}\|\nabla \bfe\|_{\mathbf{L}^2}^2+\|\Delta  \bfe\|_{\mathbf{L}^2}^2\non\\
&\quad = \int_\Omega [(\bfv-\bfv^*)\cdot\nabla](\bfd-\bfd^*) \cdot \Delta  \bfe dx
          + \int_\Omega [(\bfv^*\cdot\nabla)\bfe+(\bfw\cdot\nabla) \bfd^*]\cdot \Delta  \bfe dx \non\\
 &\qquad\ +\int_\Omega [|\bfd-\bfd^*|^2(\bfd-\bfd^*)+|\bfd-\bfd^*|^2\bfd^*+2[(\bfd-\bfd^*)\cdot\bfd^*](\bfd-\bfd^*)]\cdot \Delta  \bfe dx\non\\
 &\qquad\ + \int_\Omega [2(\bfd^*\cdot \bfe) \bfd^*+|\bfd^*|^2\bfe-\bfe]\cdot  \Delta  \bfe dx\non\\
 &\quad := \sum_{j=5}^8 r_j. \label{Dres3}
\end{align}

In what follows, we estimate the remainder terms $r_j$ ($j=1,...,8$) in \eqref{Dres1} and \eqref{Dres3}, by means of the uniform estimates
\eqref{UES} and the stability estimate \eqref{contixi}. More precisely, we can deduce that
\begin{align}
r_1&\leq \|\bfv-\bfv^*\|_{\mathbf{L}^4}\|\nabla(\bfv-\bfv^*)\|_{\mathbf{L}^2}\|\bfw\|_{\mathbf{L}^4}\non\\
   &\leq C\|\bfv-\bfv^*\|_{\mathbf{H}^1}^2\|\bfw\|_{\mathbf{L}^2}^\frac12\|\nabla \bfw\|_{\mathbf{L}^2}^\frac12\non\\
   &\leq \frac{1}{10} \|\nabla \bfw\|_{\mathbf{L}^2}^2  + \|\bfw\|_{\mathbf{L}^2}^2+ C\|\bfv-\bfv^*\|_{\mathbf{H}^1}^4\non\\
   &\leq \frac{1}{10} \|\nabla \bfw\|_{\mathbf{L}^2}^2+ \|\bfw\|_{\mathbf{L}^2}^2 +C\|\bxi\|_{\mathcal{U}}^4,\non
\end{align}
\begin{align}
 r_2&\leq \|\bfw\|_{\mathbf{L}^4}^2\|\nabla \mathbf{v}^*\|_{\mathbf{L}^2}
 \leq C\|\bfw\|_{\mathbf{L}^2}\|\nabla \bfw\|_{\mathbf{L}^2} \|\mathbf{v}^*\|_{\mathbf{H}^1}\non\\
    &\leq \frac{1}{10} \|\nabla \bfw\|_{\mathbf{L}^2}^2+ C \|\bfw\|_{\mathbf{L}^2}^2,\non
\end{align}
\begin{align}
r_3&\leq \|\nabla (\bfd-\bfd^*)\|_{\mathbf{L}^4}^2\| \nabla \bfw\|_{\mathbf{L}^2}
         \leq C\|\bfd-\bfd^*\|_{\mathbf{H}^2}^2\| \nabla \bfw\|_{\mathbf{L}^2}\non\\
   &\leq \frac{1}{10} \|\nabla \bfw\|_{\mathbf{L}^2}^2+ C\|\bxi\|_{\mathcal{U}}^4,\non
\end{align}
\begin{align}
r_4&\leq 2\|\nabla \mathbf{d}^*\|_{\mathbf{L}^4}\|\nabla \bfe\|_{\mathbf{L}^4} \| \nabla \bfw\|_{\mathbf{L}^2}\non\\
   &\leq C\|\mathbf{d}^*\|_{\mathbf{H}^2}\|\nabla \bfe\|_{\mathbf{L}^2}^\frac12\|\Delta \bfe\|_{\mathbf{L}^2}^\frac12\| \nabla \bfw\|_{\mathbf{L}^2}\non\\
   &\leq \frac{1}{10} \|\nabla \bfw\|_{\mathbf{L}^2}^2 + \frac{1}{10} \|\Delta \bfe\|_{\mathbf{L}^2}^2+ C\|\nabla \bfe\|_{\mathbf{L}^2}^2,\non
\end{align}
\begin{align}
r_5&\leq \|\bfv-\bfv^*\|_{\mathbf{L}^4}\|\nabla(\bfd -\bfd^*)\|_{\mathbf{L}^4}\|\Delta  \bfe\|_{\mathbf{L}^2}\non\\
   &\leq C\|\bfv-\bfv^*\|_{\mathbf{H}^1}\|\bfd-\bfd^*\|_{\mathbf{H}^2}\|\Delta  \bfe\|_{\mathbf{L}^2}\non\\
   &\leq \frac{1}{10} \|\Delta \bfe\|_{\mathbf{L}^2}^2 +C\|\bxi\|_{\mathcal{U}}^4,\non
\end{align}
\begin{align}
r_6&\leq \|\bfv^*\|_{\mathbf{L}^4}\|\nabla\bfe\|_{\mathbf{L}^4}\|\Delta  \bfe\|_{\mathbf{L}^2}
          +\|\bfw\|_{\mathbf{L}^4}\|\nabla \bfd^*\|_{\mathbf{L}^4}\|\Delta  \bfe\|_{\mathbf{L}^2}\non\\
   &\leq C\|\bfv^*\|_{\mathbf{H}^1}\|\nabla\bfe\|_{\mathbf{L}^2}^\frac12\|\Delta  \bfe\|_{\mathbf{L}^2}^\frac32
         + C\|\bfd^*\|_{\mathbf{H}^2} \|\bfw\|_{\mathbf{L}^2}^\frac12\|\nabla \bfw\|_{\mathbf{L}^2}^\frac12\|\Delta  \bfe\|_{\mathbf{L}^2}\non\\
   &\leq \frac{1}{10} \|\Delta \bfe\|_{\mathbf{L}^2}^2+\frac{1}{10} \|\nabla \bfw\|_{\mathbf{L}^2}^2
   + C\|\nabla\bfe\|_{\mathbf{L}^2}^2+C\|\bfw\|_{\mathbf{L}^2}^2,\non
\end{align}
\begin{align}
r_7&\leq (\|\bfd-\bfd^*\|_{\mathbf{L}^6}^3+3\|\bfd-\bfd^*\|_{\mathbf{L}^4}^2\|\bfd^*\|_{\mathbf{L}^\infty})\|\Delta  \bfe \|_{\mathbf{L}^2}\non\\
   &\leq \frac{1}{10} \|\Delta \bfe\|_{\mathbf{L}^2}^2+C(\|\bfd-\bfd^*\|_{\mathbf{H}^1}^6+\|\bfd-\bfd^*\|_{\mathbf{H}^1}^4)\non\\
   &\leq \frac{1}{10} \|\Delta \bfe\|_{\mathbf{L}^2}^2+C(\|\bxi\|_{\mathcal{U}}^2+1)\|\bxi\|_{\mathcal{U}}^4,\non
\end{align}
\begin{align}
r_8&\leq C(1+\|\bfd^*\|_{\mathbf{L}^\infty}^2)\|\bfe\|_{\mathbf{L}^2}\|\|\Delta  \bfe\|_{\mathbf{L}^2}\non\\
   &\leq C(1+\|\bfd^*\|_{\mathbf{H}^2}^2)\|\nabla \bfe\|_{\mathbf{L}^2}\|\|\Delta  \bfe\|_{\mathbf{L}^2}\non\\
   &\leq \frac{1}{10} \|\Delta \bfe\|_{\mathbf{L}^2}^2+C\|\nabla \bfe\|_{\mathbf{L}^2}^2.\non
\end{align}
Collecting the above estimates,  from the assumption \eqref{kappa} and the inequalities \eqref{Dres1}, \eqref{Dres3} we infer that
\begin{align}
&\frac{d}{dt}\big(\|\bfw\|_{\mathbf{L}^2}^2+\|\nabla \bfe\|_{\mathbf{L}^2}^2\big)+\|\nabla \bfw\|_{\mathbf{L}^2}^2+\|\Delta \bfe\|_{\mathbf{L}^2}^2\non\\
&\quad \leq C\big(\|\bfw\|_{\mathbf{L}^2}^2+\|\nabla \bfe\|_{\mathbf{L}^2}^2\big)+C(\kappa^2+1)\|\bxi\|_{\mathcal{U}}^4,\label{wees}
\end{align}
where the positive constant $C$ depends on $K_1,...,K_4$ and $\Omega$.

Taking \eqref{D5} into account, from \eqref{wees} and the Gronwall's lemma we deduce
\begin{align}
\|\bfw\|_{C([0,T];\mathbf{L}^2)}^2+\|\bfe\|_{C([0,T]; \mathbf{H}^1)}^2+\| \bfw\|_{L^2(0,T;\mathbf{H}^1)}^2
+\|\bfe\|_{L^2(0,T;\mathbf{H}^2)}^2\leq C_T\|\bxi\|_{\mathcal{U}}^4,\non
\end{align}
where the constant $C_T$ depends on $K_1,...,K_4$, $\kappa$, $\Omega$ and $T$. As a consequence, it follows that
\begin{align}
\frac{\|\mathcal{S}(\mathbf{h}^*+\bxi)-\mathcal{S}(\mathbf{h}^*)-(\bom,\bphi)\|_{\mathcal{W}_1}}{\|\bxi\|_{\mathcal{U}}}
& = \frac{\|(\bfw, \bfe)\|_{\mathcal{W}_1}}{\|\bxi\|_{\mathcal{U}}}\non\\
& \leq C_T\|\bxi\|_{\mathcal{U}} \to 0,\quad \text{as}\ \|\bxi\|_{\mathcal{U}}\to 0.\non
\end{align}
Recalling the fact $(\bom, \bphi)=\mathcal{L}_{\mathbf{h}^*}(\bxi)$, where $\mathcal{L}_{\mathbf{h}^*}$ is the linear operator defined in \eqref{LinO},
we can conclude that the control-to-state operator $\mathcal{S}$ is Fr\'echet differentiable at $\mathbf{h}^*$ in the sense of Definition \ref{diff}.
Moreover, its derivative at $\mathbf{h}^*$ is given by $\mathcal{S}'(\mathbf{h}^*)=\mathcal{L}_{\mathbf{h}^*}$. Since $\mathbf{h}^*$ is an arbitrary
vector in $\wt{\mathcal{U}}$, we arrive at the conclusion of Theorem \ref{SD}.

The proof is complete.
\end{proof}

\section{The First-order Necessary Optimality Condition} \label{necessary}
\setcounter{equation}{0}

Based on the differentiability of the control-to-state operator $\mathcal{S}$ (see Theorem \ref{SD}),
it is straightforward to derive the first-order necessary optimality condition for the optimal control problem (\textbf{CP}).
By a direct calculation, we have

\bt\label{nece1}
Let $n=2$. Suppose that the assumptions (A1)--(A2) are satisfied and $(\mathbf{v}_0, \mathbf{d}_0)\in \mathbf{V}\times
 \mathbf{H}^2(\Omega)$ is a given initial datum. Assume that $\mathbf{h}^\sharp$ is an optimal control for problem (\textbf{CP}) in the admissible set
 $\wt{\mathcal{U}}_{\mathrm{ad}}^M$ with the associate state $(\bfv^\sharp,\bfd^\sharp)=\mathcal{S}(\mathbf{h}^\sharp)$.
 Then, for any $\mathbf{h}\in \wt{\mathcal{U}}^M_{\mathrm{ad}}$, denoting by $(\bom, \bphi)=\mathcal{L}_{\mathbf{h}^\sharp}(\mathbf{h}-\mathbf{h}^\sharp)$
 the unique global weak solution to the linearized problem \eqref{L1}--\eqref{L5}  corresponding to the perturbation $\mathbf{h}-\mathbf{h}^\sharp$, we have the following variational inequality
\begin{align}
&\beta_1\int_{Q} (\bfv^\sharp -\bfv_Q)\cdot \bom dxdt
  + \beta_2\int_{Q} (\bfd^\sharp-\bfd_Q)\cdot \bphi dxdt
  + \beta_3\int_\Omega (\bfv^\sharp(T)-\bfv_\Omega)\cdot \bom(T)dx\non\\
&\quad  +\beta_4\int_\Omega (\bfd^\sharp(T)-\bfd_\Omega)\cdot \bphi(T) dx
        +\gamma\int_{\Sigma} \mathbf{h}^\sharp\cdot(\mathbf{h}-\mathbf{h}^\sharp)dSdt  \geq 0,
\quad \forall\,\mathbf{h}\in \wt{\mathcal{U}}^M_{\mathrm{ad}}.\label{FNOC1}
\end{align}
\et
\begin{proof}
Recall the reduced cost functional $\wt{\mathcal{J}}(\mathbf{h})$ introduced in the proof of Theorem \ref{existence} such that
$\wt{\mathcal{J}}(\mathbf{h}) :=\mathcal{J}((\bfv, \bfd),\mathbf{h})$ for any $\mathbf{h}\in \wt{\mathcal{U}}$,
where $(\mathbf{v},\mathbf{d})=\mathcal{S}(\mathbf{h})$.
For the cost functional $\mathcal{J}((\bfv, \bfd), \mathbf{h})$, here we regard it as $\mathcal{J}:\mathcal{W}_1\times \wt{\mathcal{U}}\to[0,+\infty)$
with spaces $\mathcal{W}_1$ and $\wt{\mathcal{U}}$ being defined in \eqref{W1cal} and \eqref{U}, respectively.

From the convexity of the admissible set $\wt{\mathcal{U}}_{\mathrm{ad}}^M$ and the well-known result \cite[Lemma 2.21]{TR10} on the necessary condition
for optimal control problems, it follows that
\begin{equation}
\wt{\mathcal{J}}'(\mathbf{h}^\sharp)(\mathbf{h}-\mathbf{h}^\sharp)\geq 0, \quad \forall\, \mathbf{h}\in \wt{\mathcal{U}}^M_{\mathrm{ad}}.\label{Jnece}
\end{equation}
To determine the operator $\wt{\mathcal{J}}'(\mathbf{h}^\sharp)$, we infer from the chain rule that
\begin{align}
\wt{\mathcal{J}}'(\mathbf{h}^\sharp)=
\mathcal{J}'_\mathbf{\mathcal{S}(\mathbf{h}^\sharp)}(\mathcal{S}(\mathbf{h}^\sharp), \mathbf{h}^\sharp)\circ \mathcal{S}'(\mathbf{h}^\sharp)
+\mathcal{J}'_{\mathbf{h}^\sharp}(\mathcal{S}(\mathbf{h}^\sharp), \mathbf{h}^\sharp).\label{FN1}
\end{align}
For every fixed $\mathbf{g}\in \wt{\mathcal{U}}$, $\mathcal{J}'_\mathbf{z}(\mathbf{z}, \mathbf{g})$ stands for the Fr\'echet derivative of
$\mathcal{J}(\mathbf{z}, \mathbf{g})$ with respect to $\mathbf{z}=(\mathbf{z}_1, \mathbf{z}_2)$, at $\mathbf{z}\in \mathcal{W}_1$.
By the definition of $\mathcal{J}$ (see \eqref{costJ}), we deduce that, for any $\mathbf{y}=(\mathbf{y}_1, \mathbf{y}_2)\in \mathcal{W}_1$, it holds
 \begin{align}
 \mathcal{J}'_\mathbf{z}(\mathbf{z}, \mathbf{g})(\mathbf{y})
& = \beta_1\int_{Q} (\mathbf{z}_1 -\bfv_Q)\cdot \mathbf{y}_1 dxdt
          + \beta_2\int_{Q}(\mathbf{z}_2-\bfv_Q)\cdot\mathbf{y}_2 dxdt\non\\
&\qquad   + \beta_3\int_\Omega (\mathbf{z}_1(T)-\bfv_\Omega)\cdot \mathbf{y}_1(T)dx
               +\beta_4\int_\Omega (\mathbf{z}_2(T)-\bfd_\Omega)\cdot \mathbf{y}_2(T) dx.\label{FN2}
\end{align}
 In a similar manner, for every fixed $\mathbf{z}\in \mathcal{W}_1$, $\mathcal{J}'_\mathbf{g}(\mathbf{z}, \mathbf{g})$ stands for the Fr\'echet derivative
 of $\mathcal{J}(\mathbf{z}, \mathbf{g})$ with respect to $\mathbf{g}$, at $\mathbf{g}\in \wt{\mathcal{U}}\subset\mathcal{U}$, namely,
 \begin{align}
 \mathcal{J}'_\mathbf{h}(\mathbf{z}, \mathbf{g})\bm{\zeta}
 =\gamma\int_{\Sigma}\mathbf{g}\cdot \bm{\zeta} dSdt,\quad \forall\, \bm{\zeta}\in \wt{\mathcal{U}}-\{\mathbf{g}\}.\label{FN3}
 \end{align}

Returning to \eqref{FN1}, here we take
 $$\mathbf{g}=\mathbf{h}^\sharp,\quad \mathbf{z}=\mathcal{S}(\mathbf{h}^\sharp)=(\bfv^\sharp,\bfd^\sharp),\quad  \bm{\zeta}=\mathbf{h}-\mathbf{h}^\sharp,$$
and by Theorem \ref{SD} we further choose
$$\mathbf{y}=\mathcal{S}'(\mathbf{h}^\sharp)(\mathbf{h}-\mathbf{h}^\sharp)=(\bom,\bphi)\in \mathcal{W}_1,$$
where $(\bom, \bphi)$ is the unique global weak solution to the linearized problem \eqref{L1}--\eqref{L5} corresponding to  $\mathbf{h}-\mathbf{h}^\sharp$.
As a consequence, it follows from Theorem \ref{SD} and the abstract chain rule \cite[Theorem 2.20]{TR10} that the operator $\wt{\mathcal{J}}'(\mathbf{h}^\sharp)$ in \eqref{FN1}
is well-defined and then we can deduce from \eqref{Jnece}--\eqref{FN3} that the variational inequality \eqref{FNOC1} holds.

The proof is complete.
\end{proof}

\section{First-order necessary optimality condition via adjoint states}
\setcounter{equation}{0}
In this section we aim to eliminate the pair $(\bom, \bphi)$ from the variational inequality \eqref{FNOC1}
and derive a first-order necessary optimality condition in terms of the optimal solution together with its adjoint states.
For this purpose, we shall first derive the corresponding adjoint system for the control problem (\textbf{CP}) (also referred to as optimality system)
that serves as the basis for computing numerical approximations of optimal solutions.
On the other hand, it turns out that the adjoint states serve as Lagrange multipliers associated with the state problem \eqref{1}--\eqref{5}.

Generally speaking, optimality systems and necessary optimality conditions can be directly deduced from the well-known abstract
Lagrange multiplier principle, i.e., the Karush--Kuhn--Tucker (KKT) theory for optimization problems in Banach spaces (see \cite{ATF87,TR10, HPUU}).
However, as it has been pointed out in \cite[Section 2.10]{TR10}, direct application of the KKT theory is rather
difficult in many cases, requiring extensive experience in matching the operators, functionals,
and involved function spaces. More precisely, the given control-to-state operators have to be differentiable in
suitably chosen Banach spaces, the adjoint operators have to be determined, and the
Lagrange multipliers must exist in the right spaces.
For instance, we can refer to \cite{FGH98,FGH05} for its nontrivial applications to optimal Dirichlet boundary control problems for the evolutionary
Navier--Stokes equations.

Since our state problem \eqref{1}--\eqref{5} is a highly nonlinear PDE system with strongly coupling structures,
it is rather difficult to apply the above mentioned abstract Lagrange multiplier principle (see e.g., \cite[Section 3.2]{ATF87})
to derive the optimality system and the existence of Lagrange multipliers.
Instead, we shall extend the formal Lagrange method described in \cite{TR10}, which turns out to be an effective tool for finding the
form of the partial differential equations from which the adjoint operator can be determined. This method  is successfully illustrated for various
elliptic and parabolic control problems (see \cite[Chapters 2-5]{TR10} for details). Roughly speaking, our strategy is as follows:
first, we derive the adjoint system by using the formal Lagrange method,
then we prove in a rigorous way the existence of adjoint states in suitable function spaces and we derive the necessary optimality condition.

\subsection{Formal derivation of the adjoint system}

For the control problem (\textbf{CP}), we introduce the Lagrange functional $\mathcal{G}$ (in a formal way)
\begin{align}
& \mathcal{G}((\mathbf{v}, \mathbf{d}), \mathbf{h}, (\tilde{\mathbf{p}}, \tilde{\mathbf{p}}_1, \tilde{P}, \tilde{\mathbf{q}}, \tilde{\mathbf{q}}_1))\non\\
&\quad :=  \mathcal{J}((\mathbf{v}, \mathbf{d}), \mathbf{h})
          - \int_Q [\p_t\mathbf{v}+ (\mathbf{v}\cdot\nabla) \mathbf{v} -\Delta \bfv
          +\nabla P+ \nabla\cdot(\nabla \mathbf{d}\odot\nabla \mathbf{d})] \cdot \tilde{\mathbf{p}} dxdt\non\\
&\qquad\ -\int_Q (\nabla\cdot \mathbf{v})\tilde{P} dxdt
         -\int_Q [\p_t\mathbf{d} +(\mathbf{v}\cdot\nabla) \mathbf{d}-\Delta \mathbf{d}+\mathbf{f}(\mathbf{d})]\cdot \tilde{\mathbf{q}} dxdt\non\\
&\qquad\ -\int_\Sigma \bfv\cdot \tilde{\mathbf{p}}_1 dSdt -\int_\Sigma (\bfd-\mathbf{h})\cdot \tilde{\mathbf{q}}_1 dSdt,
 \label{costJL}
\end{align}
for any $\mathbf{h}\in \wt{\mathcal{U}}^M_{\mathrm{ad}}$ and $(\bfv,\bfd)\in \mathcal{H}$.
\begin{remark}
Here, in the expression of $\mathcal{G}$ we only try to eliminate the five constraints due to the state problem \eqref{1}--\eqref{5} with the corresponding
Lagrange multipliers being denoted by $\tilde{\mathbf{p}}$, $\tilde{\mathbf{p}}_1$, $\tilde{P}$, $\tilde{\mathbf{q}}$, $\tilde{\mathbf{q}}_1$, respectively.
The constraint on the boundary control \eqref{Uad} is easy to handle by using variational inequalities.
\end{remark}

Let $((\bfv^\sharp,\bfd^\sharp), \mathbf{h}^\sharp)$ be a solution to the optimal control problem (\textbf{CP}) such that
$\mathbf{h}^\sharp\in \wt{\mathcal{U}}^M_{\mathrm{ad}}$, $(\bfv^\sharp,\bfd^\sharp)=\mathcal{S}(\mathbf{h}^\sharp)\in \mathcal{H}$.
Then, by the Lagrange principle, we expect that $(\mathbf{v}^\sharp, \mathbf{d}^\sharp)$ and $\mathbf{h}^\sharp$ together with
the corresponding Lagrange multipliers $\tilde{\mathbf{p}}$, $\tilde{\mathbf{p}}_1$, $\tilde{P}$, $\tilde{\mathbf{q}}$, $\tilde{\mathbf{q}}_1$
satisfy the optimality conditions associated with the minimization problem for the Lagrange functional $\mathcal{G}$ denoted by (\textbf{MG}):
\begin{align}
(\textbf{MG})\quad \quad
\min \mathcal{G}((\mathbf{v}, \mathbf{d}), \mathbf{h}, (\tilde{\mathbf{p}}, \tilde{\mathbf{p}}_1, \tilde{P}, \tilde{\mathbf{q}}, \tilde{\mathbf{q}}_1)),
\quad \text{with}\  (\mathbf{v}, \mathbf{d})\ \textit{unconstrained}\ \text{and}\ \mathbf{h} \in \wt{\mathcal{U}}^M_\mathrm{ad}.\non
\end{align}
In particular, since the pair $(\mathbf{v}, \mathbf{d})$ is now formally unconstrained, then it follows that
\begin{align}
\mathcal{G}'_{(\mathbf{v}, \mathbf{d})}((\mathbf{v}^\sharp, \mathbf{d}^\sharp), \mathbf{h}^\sharp,
(\tilde{\mathbf{p}}, \tilde{\mathbf{p}}_1, \tilde{P}, \tilde{\mathbf{q}}, \tilde{\mathbf{q}}_1))(\bom, \bphi)=0,
\label{ddd}
\end{align}
for all smooth functions $(\bom, \bphi)$ satisfying
\begin{align}
& \bom|_{t=0}=\mathbf{0},\quad \bphi|_{t=0}=\mathbf{0}, \quad \text{in}\ \Omega.\label{conss}
\end{align}
Here, the zero initial conditions for $(\bom, \bphi)$ in \eqref{conss} follow from the fact that in the control problem (\textbf{CP})
the initial datum $(\mathbf{v}_0,\mathbf{d}_0)$ for the state problem \eqref{1}--\eqref{5} is given and fixed.

A direct (formal) calculation yields that the condition \eqref{ddd} is equivalent to
\begin{align}
& \beta_1\int_Q (\bfv^\sharp -\bfv_Q)\cdot \bom dxdt
   + \beta_2\int_Q (\bfd^\sharp-\bfd_Q)\cdot \bphi dxdt\non\\
&\quad +\beta_3\int_\Omega (\bfv^\sharp(T)-\bfv_\Omega)\cdot \bom(T)dx
  +\beta_4\int_\Omega (\bfd^\sharp(T)-\bfd_\Omega)\cdot \bphi(T) dx\non\\
&\quad -\int_Q [\p_t\bom+(\bfv^\sharp\cdot\nabla)\bom -\Delta \bom+\nabla \hat{P}+ (\bom\cdot\nabla) \mathbf{v}^\sharp]\cdot \tilde{\mathbf{p}}dxdt \non\\
&\quad  -\int_Q [\nabla\cdot(\nabla \bphi\odot\nabla \mathbf{d}^\sharp)
+\nabla\cdot(\nabla \mathbf{d}^\sharp\odot\nabla \bphi)]\cdot \tilde{\mathbf{p}} dxdt
       -\int_Q (\nabla\cdot \bom)\tilde{P} dxdt \non\\
&\quad  -\int_Q [\p_t \bphi-\Delta \bphi+(\bfv^\sharp\cdot\nabla)\bphi
+(\bom\cdot\nabla) \bfd^\sharp+\mathbf{f}'(\mathbf{d}^\sharp)\bphi]\cdot \tilde{\mathbf{q}}dxdt\non\\
&\quad  -\int_\Sigma \bom\cdot \tilde{\mathbf{p}}_1 dSdt -\int_\Sigma \bphi\cdot \tilde{\mathbf{q}}_1 dSdt\non\\
&=0.\label{ddvd}
\end{align}
After integration by parts and keeping the constraint \eqref{conss} in mind, we get
\begin{align}
& \beta_1\int_Q (\bfv^\sharp -\bfv_Q)\cdot \bom dxdt
   + \beta_2\int_Q (\bfd^\sharp-\bfd_Q)\cdot \bphi dxdt\non\\
&\quad +\beta_3\int_\Omega (\bfv^\sharp(T)-\bfv_\Omega)\cdot \bom(T)dx
  +\beta_4\int_\Omega (\bfd^\sharp(T)-\bfd_\Omega)\cdot \bphi(T) dx\non\\
&\quad  + \int_Q \partial_t \tilde{\mathbf{p}} \cdot \bom dxdt
        -\int_\Omega \tilde{\mathbf{p}}(T)\cdot \bom(T) dx
        +\int_Q (\Delta \tilde{\mathbf{p}}) \cdot \bom dxdt\non\\
&\quad  +\int_\Sigma \frac{\partial \bom}{\partial \mathbf{n}}\cdot \tilde{\mathbf{p}} dSdt
        -\int_\Sigma \frac{\partial \tilde{\mathbf{p}}}{\partial \mathbf{n}}\cdot \bom dSdt
        +\int_Q \hat{P}(\nabla \cdot \tilde{\mathbf{p}}) dxdt \non\\
&\quad  - \int_\Sigma \hat{P}(\tilde{\mathbf{p}}\cdot\mathbf{n}) dSdt
        +\int_Q\nabla \tilde{P}\cdot \bom dxdt
        -\int_\Sigma (\tilde{P}\mathbf{n})\cdot \bom dSdt\non\\
&\quad  + \int_Q(\bfv^\sharp\cdot \nabla)\tilde{\mathbf{p}}\cdot \bom dxdt
        -\int_Q[(\nabla \bfv^\sharp) \tilde{\mathbf{p}}]\cdot \bom dxdt\non\\
&\quad  -\int_Q  \nabla_i(\nabla_j d^\sharp_k\nabla_j \tilde{p}_i) \phi_k dxdt
        -\int_Q  \nabla_j(\nabla_i d^\sharp_k\nabla_j \tilde{p}_i) \phi_k dxdt\non\\
&\quad  -\int_\Sigma [(\nabla \bphi \odot \nabla \bfd^\sharp)\mathbf{n}]\cdot \tilde{\mathbf{p}} dSdt
        -\int_\Sigma [(\nabla \bfd^\sharp \odot \nabla \bphi)\mathbf{n}]\cdot \tilde{\mathbf{p}} dSdt\non\\
&\quad +\int_\Sigma (\nabla_j d^\sharp_k\nabla_j \tilde{p}_i)n_i\phi_k dSdt
       +\int_\Sigma (\nabla_i d^\sharp_k\nabla_j \tilde{p}_i)n_j\phi_k dSdt \non\\
&\quad +\int_Q \partial_t \tilde{\mathbf{q}} \cdot \bphi dxdt
       -\int_\Omega \tilde{\mathbf{q}}(T) \cdot \bphi(T)dx
       + \int_Q (\Delta \tilde{\mathbf{q}})\cdot \bphi dxdt
       +\int_\Sigma \frac{\partial \bphi}{\partial \mathbf{n}}\cdot \tilde{\mathbf{q}} dsdt \non\\
&\quad -\int_\Sigma \frac{\partial \tilde{\mathbf{q}}}{\partial \mathbf{n}}\cdot \bphi dSdt
       +\int_Q (\bfv^\sharp\cdot\nabla) \tilde{\mathbf{q}} \cdot \bphi dxdt
       -\int_Q [(\nabla \bfd^\sharp)\tilde{\mathbf{q}}]\cdot \bom dxdt\non\\
&\quad -\int_Q \mathbf{f}'(\bfd^\sharp)\tilde{\mathbf{q}}\cdot \bphi dxdt
       -\int_\Sigma \tilde{\mathbf{p}}_1 \cdot \bom dSdt
       -\int_\Sigma \tilde{\mathbf{q}}_1 \cdot \bphi dSdt\non\\
&=0, \label{Iccc}
\end{align}
where $\mathbf{n}=(n_1,n_2)^\mathrm{tr}$ denotes the unit outer normal vector on $\Gamma$.

We hope to determine the Lagrange multiplies, i..e., the \emph{adjoint states}
$(\tilde{\mathbf{p}}, \tilde{\mathbf{p}}_1, \tilde{P}, \tilde{\mathbf{q}}, \tilde{\mathbf{q}}_1)$ from the equality \eqref{Iccc}.
To this end, grouping the terms in \eqref{Iccc} with respect to the test functions $(\bom, \hat{P}, \bphi)$ that can be chosen in an arbitrary way,
after integration by parts, we (at least formally) arrive at the \emph{adjoint system} in terms of partial differential equations.

More precisely, the triple $(\tilde{\mathbf{p}}, \tilde{P}, \tilde{\mathbf{q}})$ satisfies the following linear PDE system in $Q$
\begin{align}
&\p_t \tilde{\mathbf{p}}+\Delta \tilde{\mathbf{p}}
      +\nabla \tilde{P}+(\bfv^\sharp\cdot \nabla)\tilde{\mathbf{p}}
      -(\nabla \bfv^\sharp)\tilde{\mathbf{p}}
      -(\nabla\bfd^\sharp)\tilde{\mathbf{q}}
        +\beta_1 (\bfv^\sharp -\bfv_Q)
      =\mathbf{0},\label{ADJ1}\\
&\nabla \cdot \tilde{\mathbf{p}}=0,\label{ADJ2}\\
&\partial_t \tilde{\mathbf{q}}+\Delta \tilde{\mathbf{q}}
   +(\bfv^\sharp\cdot\nabla) \tilde{\mathbf{q}}-\mathbf{f}'(\bfd^\sharp)\tilde{\mathbf{q}}
   -\tilde{\mathbf{r}}(\bfd^\sharp,\tilde{\mathbf{p}})
   +\beta_2(\bfd^\sharp-\bfd_Q)=\mathbf{0},\label{ADJ3}
\end{align}
where the vector $ \tilde{\mathbf{r}}(\bfd^\sharp,\tilde{\mathbf{p}})$ is given by
 \begin{align}
 \tilde{\mathbf{r}}(\bfd^\sharp,\tilde{\mathbf{p}})
 =\nabla \cdot[\nabla^{\mathrm{tr}} \bfd^\sharp\odot(\nabla \tilde{\mathbf{p}}+\nabla^{\mathrm{tr}}\tilde{\mathbf{p}})],\label{wadj3}
 \end{align}
while on  $\Gamma$ and at time $t=T$, we have the following boundary and \emph{final} conditions
\begin{align}
&\tilde{\mathbf{p}}=\mathbf{0}, \quad \tilde{\mathbf{q}}=\mathbf{0},\quad \text{on}\ \Sigma,\label{ADJ4}\\
&\tilde{\mathbf{p}}|_{t=T}=\beta_3(\bfv^\sharp(T)-\bfv_\Omega),\quad \tilde{\mathbf{q}}|_{t=T}=\beta_4(\bfd^\sharp(T)-\bfd_\Omega),
\quad \text{in}\ \Omega.\label{ADJ5}
\end{align}
Moreover, the other two Lagrange multipliers $(\tilde{\mathbf{p}}_1,\tilde{\mathbf{q}}_1)$ corresponding to the boundary constraints, can be
uniquely determined by $(\tilde{\mathbf{p}}, \tilde{P}, \tilde{\mathbf{q}})$ such that
\begin{align}
&\tilde{\mathbf{p}}_1+\frac{\partial\tilde{\mathbf{p}}}{\partial\mathbf{n}}+\tilde{P}\mathbf{n}=\mathbf{0},\quad \text{on}\ \Sigma,\label{ADJ6a}\\
&\tilde{\mathbf{q}}_1+\frac{\partial \tilde{\mathbf{q}}}{\partial \mathbf{n}}
 =[\nabla^{\mathrm{tr}} \bfd^\sharp\odot(\nabla \tilde{\mathbf{p}}
 +\nabla^{\mathrm{tr}}\tilde{\mathbf{p}})]\mathbf{n},\quad \text{on}\ \Sigma.\label{ADJ6}
\end{align}

\subsection{Unique solvability of the adjoint system}

We proceed to prove the existence of the adjoint states $(\tilde{\mathbf{p}}, \tilde{P}, \tilde{\mathbf{q}})$ that satisfy a proper variational
formulation of problem \eqref{ADJ1}--\eqref{ADJ5}.

For this purpose, making the change of variable $t\to T-t$
to the \emph{adjoint system} \eqref{ADJ1}--\eqref{ADJ5} and denoting the new variables
\begin{align}
\mathbf{p}(t)=\tilde{\mathbf{p}}(T-t),\quad \mathbf{q}(t)=\tilde{\mathbf{q}}(T-t),\quad P(t)=\tilde{P}(T-t),\label{change}
\end{align}
we derive the following linear parabolic system for the new variables $(\mathbf{p}, P, \mathbf{q})$  in $Q$:
\begin{align}
&\p_t \mathbf{p}-\Delta \mathbf{p}-\nabla P-(\bfv^\sharp\cdot \nabla)\mathbf{p}
+(\nabla\bfv^\sharp)\mathbf{p}+(\nabla\bfd^\sharp)\mathbf{q} =\beta_1 (\bfv^\sharp -\bfv_Q),\label{AD1}\\
&\nabla \cdot \mathbf{p}=0,\label{AD2}\\
&\partial_t \mathbf{q}-\Delta \mathbf{q}-(\bfv^\sharp\cdot\nabla) \mathbf{q}+\mathbf{f}'(\bfd^\sharp)\mathbf{q} +\mathbf{r}(\bfd^\sharp,\mathbf{p})
=\beta_2(\bfd^\sharp-\bfd_Q),\label{AD3}
\end{align}
 subject to the boundary and \emph{initial} conditions
\begin{align}
&\mathbf{p}|_\Gamma=\mathbf{0}, \quad \mathbf{q}|_\Gamma=\mathbf{0},\label{AD4}\\
&\mathbf{p}|_{t=0}=\beta_3(\bfv^\sharp(T)-\bfv_\Omega),\quad \mathbf{q}|_{t=0}=\beta_4(\bfd^\sharp(T)-\bfd_\Omega),\label{AD5}
\end{align}
where in \eqref{AD3} we have $\mathbf{r}(\bfd^\sharp,\mathbf{p})=\nabla \cdot[\nabla^{\mathrm{tr}} \bfd^\sharp\odot(\nabla \mathbf{p}+\nabla^{\mathrm{tr}}\mathbf{p})]$.

Then we can prove
\bp\label{LinAD}
Let $n=2$. Suppose that $(\bfv^\sharp,\bfd^\sharp)\in \mathcal{H}$ and the assumptions (A1)--(A2) are satisfied.
Besides, we assume
\begin{align}
& \bfv_\Omega \in \mathbf{V},\quad \mathrm{if}\  \beta_3>0,\non\\
& \bfd_\Omega \in \mathbf{H}^1(\Omega)\quad \text{with}\ (\bfd^\sharp(T)-\bfd_\Omega)|_{\Gamma}=\mathbf{0},\quad \mathrm{if}\ \beta_4>0.\non
 \end{align}
 Then the linear parabolic problem \eqref{AD1}--\eqref{AD5} admits a unique weak solution $(\mathbf{p}, P, \mathbf{q})$ such that
\begin{align}
&\mathbf{p}\in C([0, T]; \mathbf{V})\cap L^2(0, T; \mathbf{H}^2(\Omega)\cap\mathbf{V})\cap H^1(0,T; \mathbf{H}),\label{ADLinreg1}\\
&P\in L^2(0,T; H^1(\Omega)) \quad \text{with} \ \  \int_\Omega Pdx=0,\label{ADLinreg1a}\\
&\mathbf{q}\in
C([0, T]; \mathbf{L}^2(\Omega)) \cap L^2(0, T; \mathbf{H}^{1}_0(\Omega))\cap H^1(0,T; \mathbf{H}^{-1}(\Omega)).\label{ADLinreg2}
\end{align}
In particular, $(\mathbf{p}, \mathbf{q})$ satisfies the following weak formulation, for a.e. $t\in (0,t)$,
\begin{align}
&\int_\Omega \p_t \mathbf{p}\cdot \mathbf{u} dx+\int_\Omega \nabla  \mathbf{p}: \nabla \mathbf{u} dx
       -\int_\Omega (\bfv^\sharp\cdot \nabla)\mathbf{p} \cdot \mathbf{u} dx\non\\
&\qquad \quad +\int_\Omega [(\nabla\bfv^\sharp)\mathbf{p}]\cdot \mathbf{u} dx
              +\int_\Omega [(\nabla\bfd^\sharp)\mathbf{q}]\cdot \mathbf{u} dx\non\\
&\quad = \beta_1 \int_\Omega (\bfv^\sharp -\bfv_Q)\cdot \mathbf{u}dx,\quad \forall\,\mathbf{u}\in \mathbf{V},\label{AD1w}\\
&\langle \partial_t \mathbf{q},\bpsi\rangle_{\mathbf{H}^{-1},\mathbf{H}^1_0}+\int_\Omega \nabla \mathbf{q}:\nabla \bpsi dx
        -\int_\Omega (\bfv^\sharp\cdot\nabla) \mathbf{q} \cdot \bpsi dx\non\\
&\qquad \quad +\int_\Omega \mathbf{f}'(\bfd^\sharp)\mathbf{q} \cdot \bpsi dx
              -\int_\Omega [\nabla^{\mathrm{tr}} \bfd^\sharp\odot(\nabla \mathbf{p}+\nabla^{\mathrm{tr}}\mathbf{p})]:\nabla \bpsi dx\non\\
&\quad  =\beta_2\int_\Omega (\bfd^\sharp-\bfd_Q)\cdot \bpsi dx,\quad \forall\, \bpsi\in \mathbf{H}^1_0(\Omega),\label{AD3w}
\end{align}
and the initial conditions \eqref{AD5}.
Furthermore, the following additional regularity for $\mathbf{q}$ holds, for some $s\in(\frac12, 1)$,
\begin{align}
\mathbf{q}\in C([0, T]; \mathbf{H}^s(\Omega)) \cap L^2(0, T; \mathbf{H}^{1+s}(\Omega))\cap H^1(0,T; \mathbf{H}^{s-1}(\Omega)).\label{ADLinreg3}
\end{align}
\ep
\begin{proof}
The proof follows from a similar argument used in Proposition \ref{Lin}, by means of the Faedo--Galerkin method.
Therefore, we simply omit the approximation scheme and just perform the necessary \emph{a priori} estimates. Here below $C_T$ stands for a positive constant
that may depend on $\|(\bfv^\sharp,\bfd^\sharp)\|_\mathcal{H}$, $\beta_1\|\bfv^\sharp-\bfv_Q\|_{L^2(0,T;\mathbf{L}^2)}$,
$\beta_2\|\bfd^\sharp- \bfd_Q\|_{L^2(0,T;\mathbf{L}^2)}$, $\beta_3\|\bfv_\Omega\|_{\mathbf{H}^1}$, $\beta_4\|\bfd_\Omega\|_{\mathbf{H}^1}$,
$\Omega$ and $T$, while $C$ denotes a constant depending on $\Omega$.

Testing \eqref{AD1w} by $S \mathbf{p}$ ($S$ being the Stokes operator) over $\Omega$ and then integrating by parts, from Young's inequality,
H\"older's inequality, the Gagliardo--Nirenberg inequality and Lemma~\ref{S}, we infer that
\begin{align}
&\frac{1}{2}\frac{d}{dt}\|\nabla \mathbf{p}\|_{\mathbf{L}^2}^2
          +\|S \mathbf{p}\|_{\mathbf{L}^2}^2\non\\
&\quad =\int_\Omega (\bfv^\sharp\cdot \nabla)\mathbf{p} \cdot  S \mathbf{p} dx
           -\int_\Omega [(\nabla\bfv^\sharp)\mathbf{p}]\cdot S \mathbf{p} dx
           -\int_\Omega [(\nabla\bfd^\sharp)\mathbf{q}] \cdot  S \mathbf{p} dx \non\\
&\qquad\quad +\int_\Omega \beta_1 (\bfv^\sharp -\bfv_Q)\cdot  S \mathbf{p} dx\non\\
&\quad \leq C \|\bfv^\sharp\|_{\mathbf{L}^4}\|\nabla \mathbf{p}\|_{\mathbf{L}^4}\|S \mathbf{p}\|_{\mathbf{L}^2}
              +C \|\mathbf{p}\|_{\mathbf{L}^4} \|\nabla \bfv^\sharp\|_{\mathbf{L}^4}\|S \mathbf{p}\|_{\mathbf{L}^2}
       \non\\
&\qquad \quad +C \|\mathbf{q}\|_{\mathbf{L}^4} \|\nabla \bfd^\sharp\|_{\mathbf{L}^4}\|S \mathbf{p}\|_{\mathbf{L}^2}
              +C\beta_1 \|\bfv^\sharp-\bfv_Q\|_{\mathbf{L}^2}  \| S \mathbf{p}\|_{\mathbf{L}^2}\non\\
&\quad \leq C \|\bfv^\sharp\|_{\mathbf{H}^1}\|\nabla \mathbf{p}\|_{\mathbf{L}^2}^\frac12\|S \mathbf{p}\|_{\mathbf{L}^2}^\frac32
+C \|\mathbf{p}\|_{\mathbf{L}^2}^\frac12 \|\nabla \mathbf{p}\|_{\mathbf{L}^2}^\frac12 \|\bfv^\sharp\|_{\mathbf{H}^2} \|S \mathbf{p}\|_{\mathbf{L}^2}\non\\
&\qquad \quad +C \|\mathbf{q}\|_{\mathbf{L}^2}^\frac12\|\nabla \mathbf{q}\|_{\mathbf{L}^2}^\frac12 \| \bfd^\sharp\|_{\mathbf{H}^2}\|S \mathbf{p}\|_{\mathbf{L}^2}
+C\beta_1 \|\bfv^\sharp-\bfv_Q\|_{\mathbf{L}^2}  \| S \mathbf{p}\|_{\mathbf{L}^2}\non\\
&\quad \leq \frac{1}{2}\| S \mathbf{p}\|_{\mathbf{L}^2}^2
             +\frac14\|\nabla \mathbf{q}\|_{\mathbf{L}^2}^2
             +C(\|\bfv^\sharp\|_{\mathbf{H}^2}^2+ \|\bfv^\sharp\|_{\mathbf{H}^1}^4+1)\|\nabla \mathbf{p}\|_{\mathbf{L}^2}^2 \non\\
&\qquad \quad +C\|\bfd^\sharp\|_{\mathbf{H}^2}^4\|\mathbf{q}\|_{\mathbf{L}^2}^2 +C \beta_1^2 \|\bfv^\sharp-\bfv_Q\|_{\mathbf{L}^2}^2.\label{ADes1}
\end{align}
Next, testing \eqref{AD3w} by $\mathbf{q}$ and  integrating over $\Omega$, and using \eqref{ADJ4}, we get
\begin{align}
&\frac12\frac{d}{dt} \|\mathbf{q}\|_{\mathbf{L}^2}^2+\|\nabla \mathbf{q}\|^2_{\mathbf{L}^2}\non\\
 &\quad = -\int_\Omega \mathbf{f}'(\bfd^\sharp)\mathbf{q}\cdot \mathbf{q} dx
 -\int_\Omega \mathbf{r}(\bfd^\sharp,\mathbf{p})\cdot \mathbf{q} dx+\int_\Omega \beta_2(\bfd^\sharp-\bfd_Q)\cdot \mathbf{q}dx\non\\
&\quad \leq C(\|\mathbf{d}^\sharp\|_{\mathbf{L}^\infty}^2+1)\|\mathbf{q}\|_{\mathbf{L}^2}^2
            +C \|\nabla \bfd^\sharp\|_{\mathbf{L}^\infty}\|\nabla \mathbf{p}\|_{\mathbf{L}^2}\|\nabla \mathbf{q}\|_{\mathbf{L}^2}\non\\
&\qquad \quad +C\beta_2 \|\mathbf{q}\|_{\mathbf{L}^2}\|\bfd^\sharp-\mathbf{d}_Q\|_{\mathbf{L}^2}\non\\
&\quad \leq \frac14\|\nabla \mathbf{q}\|_{\mathbf{L}^2}^2 +C \| \bfd^\sharp\|_{\mathbf{H}^3}^2\|\nabla \mathbf{p}\|_{\mathbf{L}^2}^2
            +C (\|\mathbf{d}^\sharp\|_{\mathbf{H}^2}^2+1)\|\mathbf{q}\|_{\mathbf{L}^2}^2 \non\\
&\qquad \quad  +C\beta_2^2 \|\bfd^\sharp-\mathbf{d}_Q\|_{\mathbf{L}^2}^2.\label{ADes2}
\end{align}
Then adding \eqref{ADes1} and \eqref{ADes2} together, since $(\bfv^\sharp,\bfd^\sharp)\in \mathcal{H}$ (recall the global estimates
\eqref{globstrong1}, \eqref{globstrong2} for $(\bfv^\sharp,\bfd^\sharp)$),  from the Gronwall lemma and Lemma~\ref{S}, it easily follows that
\begin{align}
&\| \mathbf{p}(t)\|_{\mathbf{H}^1}^2+\|\mathbf{q}(t)\|_{\mathbf{L}^2}^2+\int_0^t(\| \mathbf{p}(\tau)\|_{\mathbf{H}^2}^2
+\|\mathbf{q}(\tau)\|_{\mathbf{H}^1}^2)d\tau \leq C_T,\quad \forall\, t\in[0,T],\label{espq}
\end{align}
and $P\in L^2(0,T; H^1(\Omega))$.
Estimate \eqref{espq} combined with equations \eqref{AD1w} and \eqref{AD3w} also yields uniform estimates for time derivatives of $(\mathbf{p}, \mathbf{q})$
such that $\p_t\mathbf{p}\in L^2(0,T;\mathbf{H})$, $\p_t \mathbf{q}\in L^2(0,T; \mathbf{H}^{-1}(\Omega))$.

Therefore, by the standard compactness argument we can conclude the existence of a weak solution $(\mathbf{p}, P, \mathbf{q})$ to problem \eqref{AD1}--\eqref{AD5}.
Besides, the uniqueness is a straightforward consequence of the energy method for the linear parabolic system.

Next, we improve the spatial regularity of $\mathbf{q}$. For $s\in (\frac12,1)$, by the Sobolev embedding theorem ($n=2$), we see that
$$\mathbf{H}^1(\Omega)\hookrightarrow {\mathbf L}^\frac{2}{1-s}(\Omega),
\quad \mathbf{H}^{1-s}(\Omega)\hookrightarrow \mathbf{L}^{\sigma'}(\Omega),\quad \text{with}\ \ \frac{1}{\sigma'}=\frac12-\frac{1-s}{2}$$
and thus
$$\mathbf{L}^{\sigma}(\Omega)\hookrightarrow \mathbf{H}^{s-1}(\Omega),
\quad\text{with}\ \ \frac{1}{\sigma}=\frac12+\frac{1-s}{2},\ \ \text{i.e.},\ \ \sigma=\frac{2}{2-s}.$$
As a consequence, from the H\"older inequality it follows
\begin{align}
\|\mathbf{r}(\bfd^\sharp,\mathbf{p})\|_{\mathbf{H}^{s-1}}
&\leq\ C\|\mathbf{r}(\bfd^\sharp,\mathbf{p})\|_{\mathbf{L}^\frac{2}{2-s}}\non\\
&\leq\ C\|\nabla \bfd^\sharp\|_{\mathbf{L}^\frac{2}{1-s}}\|\mathbf{p}\|_{\mathbf{H}^2}
       +C\|\Delta\bfd^\sharp\|_{\mathbf{L}^\frac{2}{1-s}} \|\nabla \mathbf{p}\|_{\mathbf{L}^2}\non\\
&\leq \ C\|\bfd^\sharp\|_{\mathbf{H}^2}\|\mathbf{p}\|_{\mathbf{H}^2}
       +C\|\bfd^\sharp\|_{\mathbf{H}^3} \|\nabla \mathbf{p}\|_{\mathbf{L}^2},\label{rr}
\end{align}
and
\begin{align}
\|(\bfv^\sharp\cdot\nabla) \mathbf{q}\|_{\mathbf{H}^{s-1}}
&\leq\ C\|(\bfv^\sharp\cdot\nabla) \mathbf{q}\|_{\mathbf{L}^\frac{2}{2-s}}\non\\
&\leq\ C\|\bfv^\sharp\|_{\mathbf{L}^\frac{2}{1-s}}\|\nabla \mathbf{q}\|_{\mathbf{L}^2}\non\\
&\leq\ C\|\bfv^\sharp\|_{\mathbf{H}^1}\|\nabla \mathbf{q}\|_{\mathbf{L}^2}.\label{rr1}
\end{align}
Then from \eqref{espq}, \eqref{rr}, \eqref{rr1} and the fact $(\bfv^\sharp,\bfd^\sharp)\in \mathcal{H}$  we infer
\begin{align}
& \int_0^T\|\mathbf{r}(\bfd^\sharp,\mathbf{p})\|_{\mathbf{H}^{s-1}}^2 dt\non\\
&\quad  \leq\   C\sup_{t\in[0,T]}\|\bfd^\sharp(t)\|^2_{ \mathbf{H}^2} \int_0^T\|\mathbf{p}(t)\|_{\mathbf{H}^2}^2 dt
         +C\sup_{t\in[0,T]}\|\nabla \mathbf{p}(t)\|_{\mathbf{L}^2}^2 \int_0^T \|\bfd^\sharp(t)\|_{\mathbf{H}^3}^2 dt\non\\
&\quad \leq \ C_T,\label{esr}
\end{align}
and
\begin{align}
\int_0^T \|(\bfv^\sharp\cdot\nabla) \mathbf{q}\|_{\mathbf{H}^{s-1}}^2 dt
&\leq\ C\sup_{t\in[0,T]}\|\bfv^\sharp(t)\|_{\mathbf{H}^1}^2\int_0^T \|\nabla \mathbf{q}(t)\|_{\mathbf{L}^2}^2dt\non\\
&\leq\ C_T.\label{esr1}
\end{align}
Hence, testing \eqref{AD3w} by $(-\Delta)^s\mathbf{q}$ in $\mathbf{H}^{s}(\Omega)$, from the H\"older inequality and Young's inequality we obtain
\begin{align}
&\frac12\frac{d}{dt}\|\mathbf{q}\|_{\mathbf{H}^s}^2+ \|\mathbf{q}\|_{\mathbf{H}^{1+s}}^2\non\\
&\quad = \int_\Omega [(\bfv^\sharp\cdot\nabla) \mathbf{q}] \cdot (-\Delta)^s\mathbf{q} dx
         -\int_\Omega \mathbf{f}'(\bfd^\sharp)\mathbf{q}\cdot (-\Delta)^s\mathbf{q} dx
         -\int_\Omega \mathbf{r}(\bfd^\sharp,\mathbf{p})\cdot (-\Delta)^s\mathbf{q} dx\non\\
&\qquad \quad +\int_\Omega \beta_2(\bfd^\sharp-\bfd_Q)\cdot (-\Delta)^s\mathbf{q}dx\non\\
&\quad \leq  C\|(\bfv^\sharp\cdot\nabla) \mathbf{q}\|_{\mathbf{H}^{s-1}}\|\mathbf{q}\|_{\mathbf{H}^{1+s}}
            + C\|\mathbf{r}(\bfd^\sharp,\mathbf{p})\|_{\mathbf{H}^{s-1}}\|\mathbf{q}\|_{\mathbf{H}^{1+s}}
            + C\|\mathbf{f}'(\bfd^\sharp)\mathbf{q}\|_{\mathbf{L}^2}\|\mathbf{q}\|_{\mathbf{H}^{2s}}\non\\
&\qquad \quad
            + C\beta_2\|\mathbf{q}\|_{\mathbf{H}^{2s}}\|\bfd^\sharp-\mathbf{d}_Q\|_{\mathbf{L}^2}
            \non\\
&\quad \leq \frac12\|\mathbf{q}\|_{\mathbf{H}^{1+s}}^2
             + C\|(\bfv^\sharp\cdot\nabla) \mathbf{q}\|_{\mathbf{H}^{s-1}}^2
             + C\|\mathbf{r}(\bfd^\sharp,\mathbf{p})\|_{\mathbf{H}^{s-1}}^2\non\\
&\qquad \quad
             + C(\|\mathbf{d}^\sharp\|_{\mathbf{L}^\infty}^2+1)\|\mathbf{q}\|_{\mathbf{L}^2}^2
             + C\beta_2^2 \|\bfd^\sharp-\mathbf{d}_Q\|_{\mathbf{L}^2}^2.\label{Hsn}
\end{align}
From \eqref{esr}--\eqref{Hsn} and the Gronwall lemma, it follows
\begin{align}
\|\mathbf{q}(t)\|_{\mathbf{H}^s}^2+\int_0^t \|\mathbf{q}(\tau)\|_{\mathbf{H}^{1+s}}^2d\tau\leq C_T,\quad t\in[0,T].\label{Hsn1}
\end{align}
Then by a comparison argument in \eqref{AD3} we also have $\p_t \mathbf{q}\in L^2(0,T; \mathbf{H}^{s-1}(\Omega))$, which combined with
the interpolation theory further implies $\mathbf{q}\in C([0,T]; \mathbf{H}^s(\Omega))$.

The proof is complete.
\end{proof}
Recalling the change of variables \eqref{change} and the relations \eqref{ADJ6a}--\eqref{ADJ6}, we have
\bc\label{adstate}
Under the same assumptions of Proposition~\ref{LinAD}, the adjoint system \eqref{ADJ1}--\eqref{ADJ5} admits a unique weak solution
$(\tilde{\mathbf{p}}, \tilde{P}, \tilde{\mathbf{q}})$, satisfying the same properties as for the weak solution $(\mathbf{p}, P, \mathbf{q})$
to problem \eqref{AD1}--\eqref{AD5} stated in Proposition~\ref{LinAD}. Moreover, the Lagrange multipliers $(\mathbf{p}_1, \mathbf{q}_1)$ are
uniquely determined by \eqref{ADJ6a}, \eqref{ADJ6} such that
\begin{align}
\mathbf{p}_1\in L^2(0,T;\mathbf{H}^\frac12(\Gamma)),\quad \mathbf{q}_1\in L^\frac43(0,T;\mathbf{H}^\frac12(\Gamma)).\label{regpq1}
\end{align}
\ec
\begin{proof}
We only need to prove \eqref{regpq1}. It easily follows from \eqref{ADJ6a}, \eqref{change}, \eqref{ADLinreg1}, \eqref{ADLinreg1a} and the trace
theorem that $\mathbf{p}_1\in L^2(0,T;\mathbf{H}^\frac12(\Gamma))$. Next, by \eqref{change}, \eqref{ADLinreg2} and the trace theorem we also have
$\frac{\partial \tilde{\mathbf{q}}}{\partial \mathbf{n}}\in L^2(0,T;\mathbf{H}^\frac12(\Gamma))$. On the other hand,  using Sobolev embedding theorem
and Agmon's inequality ($n=2$), we deduce
\begin{align}
&\|\nabla^{\mathrm{tr}} \bfd^\sharp\odot(\nabla \tilde{\mathbf{p}}
 +\nabla^{\mathrm{tr}}\tilde{\mathbf{p}})\|_{\mathbf{H}^1}\non\\
&\quad \leq \|\nabla^{\mathrm{tr}} \bfd^\sharp\odot\nabla \tilde{\mathbf{p}} \|_{\mathbf{H}^1}
+\|\nabla^{\mathrm{tr}} \bfd^\sharp\odot \nabla^{\mathrm{tr}}\tilde{\mathbf{p}}\|_{\mathbf{H}^1}\non\\
&\quad\leq  C\|\bfd^\sharp\|_{\mathbf{W}^{2,4}}\|\tilde{\mathbf{p}}\|_{\mathbf{W}^{1,4}}
      +C\|\bfd^\sharp\|_{\mathbf{W}^{1,\infty}}\|\tilde{\mathbf{p}}\|_{\mathbf{H}^2}
      +C\|\bfd^\sharp\|_{\mathbf{W}^{1,4}}\|\tilde{\mathbf{p}}\|_{\mathbf{W}^{1,4}}\non\\
&\quad \leq C\|\bfd^\sharp\|_{\mathbf{H}^{3}}^\frac12\|\tilde{\mathbf{p}}\|_{\mathbf{H}^{2}}^\frac12\|\bfd^\sharp\|_{\mathbf{H}^{2}}^\frac12
\|\tilde{\mathbf{p}}\|_{\mathbf{H}^{1}}^\frac12
       +C\|\bfd^\sharp\|_{\mathbf{H}^{3}}^\frac12\|\bfd^\sharp\|_{\mathbf{H}^{1}}^\frac12 \|\tilde{\mathbf{p}}\|_{\mathbf{H}^2}
       +C\|\bfd^\sharp\|_{\mathbf{H}^{2}}\|\tilde{\mathbf{p}}\|_{\mathbf{H}^{2}}.\non
\end{align}
As a consequence, we infer from the fact $(\bfv^\sharp, \bfd^\sharp)\in \mathcal{H}$ and \eqref{ADLinreg1} that
$$\nabla^{\mathrm{tr}} \bfd^\sharp\odot(\nabla \tilde{\mathbf{p}}
 +\nabla^{\mathrm{tr}}\tilde{\mathbf{p}})\in L^\frac43(0,T;\mathbf{H}^1(\Omega)),$$
which implies
$$[\nabla^{\mathrm{tr}} \bfd^\sharp\odot(\nabla \tilde{\mathbf{p}}
 +\nabla^{\mathrm{tr}}\tilde{\mathbf{p}})]\mathbf{n}\in L^\frac43(0,T;\mathbf{H}^\frac12(\Gamma)).$$
Therefore, it follows from \eqref{ADJ6} that  $\mathbf{q}_1\in L^\frac43(0,T;\mathbf{H}^\frac12(\Gamma))$.

The proof is complete.
\end{proof}

\begin{remark}
We can easily see from the regularity of the boundary control $\mathbf{h}$, the state constraints $(\bfv, \bfd)$ and the
Lagrange multipliers  $(\tilde{\mathbf{p}}, \tilde{\mathbf{p}}_1, \tilde{P}, \tilde{\mathbf{q}}, \tilde{\mathbf{q}}_1)$
that the Lagrange functional $\mathcal{G}$ in \eqref{costJL} is well-defined.
\end{remark}

\subsection{The first-order necessary optimality condition via adjoint states}

Now we are able to eliminate the pair $(\bom, \bphi)$ from the variational inequality \eqref{FNOC1} and, alternatively, form the first-order
necessary optimality condition by the state system \eqref{1}--\eqref{5} for $(\bfv^\sharp, \bfd^\sharp)$ together with the adjoint system
\eqref{ADJ1}--\eqref{ADJ6}:

\bt\label{nece2}
Let $n=2$. Assume that (A1)--(A2) are satisfied, $(\mathbf{v}_0, \mathbf{d}_0)\in \mathbf{V}\times \mathbf{H}^2(\Omega)$  and
\begin{align}
& \bfv_\Omega \in \mathbf{V},\quad \mathrm{if}\  \beta_3>0,\non\\
& \bfd_\Omega \in \mathbf{H}^1(\Omega)\quad \text{with}\ (\bfd^\sharp(T)-\bfd_\Omega)|_{\Gamma}=\mathbf{0},\quad \mathrm{if}\ \beta_4>0.\non
 \end{align}
Besides, suppose that $\mathbf{h}^\sharp$ is an optimal boundary control for the control problem (\textbf{CP}) in the admissible set $\wt{\mathcal{U}}_{\mathrm{ad}}^M$
with the associate state
$(\bfv^\sharp,\bfd^\sharp)=\mathcal{S}(\mathbf{h}^\sharp)$ as well as the adjoint state $(\tilde{\mathbf{p}}, \tilde{\mathbf{q}})$.
Then, for any $\mathbf{h}\in \wt{\mathcal{U}}^M_{\mathrm{ad}}$, the following variational inequality holds
\begin{align}
&\gamma\int_\Sigma\mathbf{h}^\sharp\cdot(\mathbf{h}-\mathbf{h}^\sharp)dSdt
-\int_\Sigma [(\nabla \tilde{\mathbf{q}})  \mathbf{n})]\cdot (\mathbf{h}-\mathbf{h}^\sharp) dSdt\non\\
&\qquad +\int_\Sigma ([\nabla^{\mathrm{tr}} \bfd^\sharp\odot(\nabla \tilde{\mathbf{p}}
+\nabla^{\mathrm{tr}}\tilde{\mathbf{p}})]\mathbf{n})\cdot (\mathbf{h}-\mathbf{h}^\sharp) dSdt\non\\
&\quad \geq 0,
\quad \forall\,\mathbf{h}\in \wt{\mathcal{U}}^M_{\mathrm{ad}}.\label{FNOC2}
\end{align}
\et
\begin{proof} Concerning the minimization problem (\textbf{MG}) for the Lagrange functional $\mathcal{G}$, since the control function $\mathbf{h}$ is
constrained such that $\mathbf{h} \in \wt{\mathcal{U}}^M_\mathrm{ad}$, then, by the Lagrange principle, we have
\begin{align}
\mathcal{G}'_{\mathbf{h}}((\mathbf{v}^\sharp, \mathbf{d}^\sharp), \mathbf{h}^\sharp, (\tilde{\mathbf{p}}, \tilde{\mathbf{p}}_1, \tilde{P}, \tilde{\mathbf{q}},
\tilde{\mathbf{q}}_1))(\mathbf{h}-\mathbf{h}^\sharp)\geq 0, \quad \forall\,\mathbf{h}\in \wt{\mathcal{U}}^M_{\mathrm{ad}}.
\label{hhh}
\end{align}
A direction calculation yields that
\begin{align}
\gamma \int_\Sigma \mathbf{h}^\sharp\cdot(\mathbf{h}-\mathbf{h}^\sharp)dSdt
       +\int_\Sigma \tilde{\mathbf{q}}_1 \cdot(\mathbf{h}-\mathbf{h}^\sharp)dSdt \geq 0,
\quad \forall\,\mathbf{h}\in \wt{\mathcal{U}}^M_{\mathrm{ad}}.\label{FNOC3}
\end{align}
On the other hand, we recall that $\tilde{\mathbf{q}}_1$ can be uniquely determined by \eqref{ADJ6}, i.e.,
\begin{align}
\tilde{\mathbf{q}}_1=-(\nabla \tilde{\mathbf{q}}) \mathbf{n}+ [\nabla^{\mathrm{tr}} \bfd^\sharp\odot(\nabla \tilde{\mathbf{p}}
 +\nabla^{\mathrm{tr}}\tilde{\mathbf{p}})]\mathbf{n},\label{FNOC3a}
 \end{align}
where $(\tilde{\mathbf{p}}, \tilde{\mathbf{q}})$ are the adjoint states associated with $\mathbf{h}^\sharp$ that can be obtained in Corollary \ref{adstate}.
Hence, from the variational inequality \eqref{FNOC3} and \eqref{FNOC3a} we arrive at our conclusion \eqref{FNOC2}.

The proof is complete.
\end{proof}
\begin{remark}
Based on the Lagrange principle, an alternative proof of \eqref{FNOC2} can be given as follows.
For any $\mathbf{h}\in \wt{\mathcal{U}}^M_{\mathrm{ad}}$, let $(\bom, \bphi)=\mathcal{L}_{\mathbf{h}^\sharp}(\mathbf{h}-\mathbf{h}^\sharp)$ be
the unique global weak solution to the linearized problem \eqref{L1}--\eqref{L5} corresponding to  $\bxi=\mathbf{h}-\mathbf{h}^\sharp$.
Then from the optimality condition \eqref{ddvd} and the facts $\bom|_\Gamma =\mathbf{0}$, $\bphi|_\Gamma=\mathbf{h}-\mathbf{h}^\sharp$, it follows that
\begin{align}
& \beta_1\int_Q (\bfv^\sharp -\bfv_Q)\cdot \bom dxdt
   + \beta_2\int_Q (\bfd^\sharp-\bfd_Q)\cdot \bphi dxdt\non\\
&\quad +\beta_3\int_\Omega (\bfv^\sharp(T)-\bfv_\Omega)\cdot \bom(T)dx
  +\beta_4\int_\Omega (\bfd^\sharp(T)-\bfd_\Omega)\cdot \bphi(T) dx\non\\
&= \int_\Sigma \tilde{\mathbf{q}}_1 \cdot(\mathbf{h}-\mathbf{h}^\sharp)dSdt, \quad \forall\,\mathbf{h}\in \wt{\mathcal{U}}^M_{\mathrm{ad}},\label{eqdd}
\end{align}
which together with the variational inequality \eqref{FNOC1} immediately yields our conclusion \eqref{FNOC2}.
\end{remark}

Finally, we remark that \eqref{FNOC2} allows for the interpretation of the optimal boundary control via the following projection formula:

\bc
Suppose that $\gamma>0$ and $\mathbf{h}^\sharp\in \wt{\mathcal{U}}^M_{\mathrm{ad}}$ is an optimal boundary control to problem (\textbf{CP}).
Then $\mathbf{h}^\sharp$ together with its associated adjoint states $(\tilde{\mathbf{p}}, \tilde{\mathbf{q}})$ satisfies the projection formula
$$
\mathbf{h}^\sharp= \frac{1}{\gamma}\mathbb{P}_{\wt{\mathcal{U}}^M_{\mathrm{ad}}}\left\{\big[\nabla \tilde{\mathbf{q}}
-\nabla^{\mathrm{tr}} \bfd^\sharp\odot(\nabla \tilde{\mathbf{p}}+\nabla^{\mathrm{tr}}\tilde{\mathbf{p}})\big] \mathbf{n}\right\},
$$
where $\mathbb{P}_{\wt{\mathcal{U}}^M_{\mathrm{ad}}}$ is the orthogonal projector in $\mathbf{L}^2(\Sigma)$ onto the convex set $\wt{\mathcal{U}}^M_{\mathrm{ad}}$.
\ec

\section{Appendix}
\label{App}
 \setcounter{equation}{0}

\subsection{Lifting functions and variable transformation}

The natural energy functional of the problem \eqref{1}--\eqref{5} consists of the kinetic and elastic potential energies, which is given by
 \be
  \mathcal{E}(t)=\frac12\|\mathbf{v}(t)\|_{\mathbf{L}^2}^2+\frac12
 \|\nabla\mathbf{d}(t)\|_{\mathbf{L}^2}^2+\int_\Omega F(\mathbf{d}(t)) dx, \quad t\in[0,T].\label{OE}
 \ee
However, due to the time-dependent boundary condition $\mathbf{h}(t,x)$ for the molecule director $\mathbf{d}$,
one cannot expect the total energy $\mathcal{E}(t)$ to be decreasing in time as in the autonomous case \cite{LL95}.

In order to overcome this difficulty and obtain proper energy estimates for global weak solutions,
suitable lifting functions were introduced in the literature (see, e.g., \cite{B,C06,C09, GW13}).
The first lifting function
$\mathbf{d}_E=\mathbf{d}_E(x,t)$ is of \emph{elliptic} type (see, e.g., \cite{B,C06,C09})
 \be
 \begin{cases}
 -\Delta \mathbf{d}_E=\mathbf{0},\qquad\ \ \ \ \text{ in } \Omega\times (0,T),\\
 \mathbf{d}_E=\mathbf{h}(x,t),\qquad \ \ \text{ on } \Gamma\times (0,T).
 \end{cases}\label{LE}
 \ee
The second lifting function $\mathbf{d}_P=\mathbf{d}_P(x,t)$ is of \emph{parabolic} type (see, e.g., \cite{GW13})
 \be
 \begin{cases}
 \p_t\mathbf{d}_P-\Delta \mathbf{d}_P=\mathbf{0},\qquad \text{ in } \Omega\times (0,T),\\
 \mathbf{d}_P=\mathbf{h}(x,t),\qquad \qquad  \text{ on } \Gamma\times (0,T),\\
 \mathbf{d}_P|_{t=0}= \mathbf{d}_{E0}, \qquad \quad\ \  \text{ in } \Omega,
 \end{cases}\label{LP}
 \ee
 where the initial data $\mathbf{d}_{E0}$ can be viewed as a lifting function for
the original initial datum $\bfd_0$:
 \be
 \begin{cases}
 -\Delta \mathbf{d}_{E0}=\mathbf{0},\qquad \text{ in } \Omega,\\
 \mathbf{d}_{E0}=\mathbf{d}_0,\qquad \quad \ \text{ on } \Gamma.
 \end{cases}\label{iE0}
 \ee
\begin{remark}
(i) Since we assume $\mathbf{d}_0|_\Gamma=\mathbf{h}|_{t=0}$, then the elliptic lifting functions \eqref{LE} and \eqref{iE0} are compatible
at $t=0$, that is $\bfd_E|_{t=0}=\bfd_{E0}$.

(ii) The elliptic lifting function $\mathbf{d}_E$ will be helpful to get uniform lower-order energy estimates for the solution
$(\bfv(t),\bfd(t))$ in $\mathbf{L}^2(\Omega)\times\mathbf{H}^1(\Omega)$ (see Lemma \ref{BEL}), which is crucial to prove the existence of global
weak solutions to problem \eqref{1}--\eqref{5}.

(iii) The parabolic lifting function $\mathbf{d}_P$ implies the property $(\Delta(\mathbf{d}-\mathbf{d}_P)-\mathbf{f}(\mathbf{d}))|_\Gamma=\mathbf{0}$,
which enables us to perform integration by parts and derive some higher-order differential inequalities for problem \eqref{1}--\eqref{5} (see Lemma \ref{H2D}).
We note that the lifting function $\mathbf{d}_P$ is slightly different from the one introduced in \cite{B,C09} since they have different initial values.
Both choices work for the existence of strong solutions \cite{B,C09,GW13}, while the current definition of $\mathbf{d}_P$ also turns out to be convenient
for the study of long-time behavior \cite{GW13}.
\end{remark}

Below we report some properties of the lifting functions $\mathbf{d}_E$ and $\mathbf{d}_P$ that have been
used in the previous sections. Here, we denote by $c$ a generic
positive constant which depends on the spatial dimension $n$ and $\Omega$ at most.

From the classical elliptic regularity theory \cite{LM,Tay} it follows  (see \cite[Lemmas A.8, A.9]{B})
 \bl\label{Ap1}
Suppose that $\mathbf{h}$ satisfies \eqref{hyp2} and \eqref{hyp3}, then the lifting problem \eqref{LE} admits a unique strong solution
$$\bfd_E\in L^2(0,T; \mathbf{H}^3(\Omega))\cap L^\infty(0,T;\mathbf{H}^2(\Omega)) \cap W^{1,4}(0,T; \mathbf{H}^1(\Omega))$$
such that, for $t\in[0,T]$ and $k=0,1$,
\begin{align}
 &\int_0^t \|\mathbf{d}_E(\tau)\|_{\mathbf{H}^{k+2}}^2 d\tau
 \leq c\int_0^t\|\mathbf{h}(\tau)\|_{\mathbf{H}^{\frac32+k}(\Gamma)}^2d\tau,\nonumber\\
 &\int_0^t \|\partial_t\mathbf{d}_E(\tau)\|_{\mathbf{H}^{k}}^4 dt
 \leq c\int_0^t\|\partial_t\mathbf{h}(\tau)\|_{\mathbf{H}^{k-\frac12}(\Gamma)}^4d\tau.\nonumber
  \end{align}
 \el
 Concerning the parabolic lifting problem \eqref{LP}, we have the following result (see, e.g., \cite[Lemmas A.2-A.4]{B} and also \cite[Lemma 6.2]{GW13})
 \bl \label{Ap2}
 (1) Suppose that $\mathbf{d}_0\in \mathbf{H}^1(\Omega)$ 
 and the assumptions
 \eqref{hyp2}--\eqref{hyp4} are satisfied.
 Then the lifting problem \eqref{LP} admits a unique weak solution
 $$\bfd_P\in L^\infty(0,T;\mathbf{H}^1(\Omega))\cap L^2(0,T;\mathbf{H}^2(\Omega))\cap H^1(0,T;\mathbf{L}^2(\Omega))$$
 such that, for $t\in[0,T]$,
 \begin{align}
  &\|\mathbf{d}_P(t)\|^2_{\mathbf{H}^1}
     +\int_0^t\Big(\|\mathbf{d}_P(\tau)\|^2_{\mathbf{H}^2}+ \|\p_t\mathbf{d}_P(\tau)\|_{\mathbf{L}^2}^2\Big)d\tau\non\\
  &\quad    \leq \|\bfd_{E0}\|_{\mathbf{H}^1}^2
       + c\int_0^t\Big(\|\mathbf{h}(\tau)\|_{\mathbf{H}^{\frac32}(\Gamma)}^2
       +\|\partial_t\mathbf{h}(\tau)\|_{\mathbf{H}^{-\frac12}(\Gamma)}^2\Big)d\tau.\label{A3h1}
 \end{align}

(2) If $\mathbf{d}_0\in \mathbf{H}^2(\Omega)$ and the assumptions
 \eqref{hyp2s}, \eqref{hyp3s} and \eqref{hyp4} are satisfied, then the lifting problem \eqref{LP} admits a unique strong solution
 $$\bfd_P\in L^\infty(0,T;\mathbf{H}^2(\Omega))\cap L^2(0,T;\mathbf{H}^3(\Omega))\cap W^{1,\infty}(0,T;\mathbf{L}^2(\Omega))\cap H^1(0,T;\mathbf{H}^1(\Omega))$$
  such that, for $t\in[0,T]$,
 \begin{align}
  &\|\mathbf{d}_P(t)\|^2_{\mathbf{H}^2}+\int_0^t \Big(\|\mathbf{d}_P(\tau)\|^2_{\mathbf{H}^3} + \|\p_t\mathbf{d}_P(\tau)\|_{\mathbf{H}^1}^2\Big) d\tau\nonumber\\
  &\quad\ \ \leq \|\bfd_{E0}\|_{\mathbf{H}^2}^2
       + c\int_0^t\Big(\|\mathbf{h}(\tau)\|_{\mathbf{H}^{\frac52}(\Gamma)}^2
       +\|\partial_t\mathbf{h}(\tau)\|_{\mathbf{H}^{\frac12}(\Gamma)}^2\Big)d\tau,\label{A4h2}\\
 &\|\p_t\mathbf{d}_P(t)\|_{\mathbf{L}^2}^2 \leq c\int_0^t\|\partial_t\mathbf{h}(\tau)\|^2_{\mathbf{H}^{\frac12}(\Gamma)}d\tau.\label{A5h2}
 \end{align}
 Besides, if \eqref{hyp5} is fulfilled, then the solution
  $\bfd_P$ satisfies the weak maximum principle such that $$|\bfd_P|_{\mathbb{R}^n}\leq 1,\quad \text{a.e. in}\ \Omega\times[0,T].$$

\el
The following lemma states the relation between the two lifting functions $\bfd_E$ and $\bfd_P$  (see \cite[Lemma 6.2]{GW13}):
\bl\label{Ap3}
Under the assumptions of Lemmas \ref{Ap1}, \ref{Ap2}, we have
\be
  \|\mathbf{d}_P(t)-\mathbf{d}_E(t)\|_{\mathbf{H}^k}^2
 \leq c\int_0^t\|\partial_t\mathbf{h}(\tau)\|^2_{\mathbf{H}^{k-\frac32}(\Gamma)}d\tau,\quad k=1,2,\quad \forall\, t \in [0,T].\label{A5}
 \ee
\el

Using the lifting functions introduced above, problem \eqref{1}--\eqref{5} can be written into different equivalent forms,
which will be convenient to obtain suitable \emph{a priori} estimates for the solution via energy method.

 Let $\bfd_E$ be defined as in \eqref{LE}. Set
 \be
 \wha{\mathbf{d}}=\mathbf{d}-\mathbf{d}_E.\label{dhat}
 \ee
 Then the system \eqref{1}--\eqref{3} can be rewritten into the following form for $(\bfv, \wha{\bfd})$:
  \bea
 \mathbf{v}_t+\mathbf{v}\cdot\nabla \mathbf{v}- \Delta \mathbf{v}
 +\nabla P&=&- (\nabla \mathbf{d})^{\mathrm{tr}} \Delta \wha{\mathbf{d}},\label{1E}\\
 \nabla \cdot \mathbf{v} &=& 0,\label{2E}\\
 \wha{\mathbf{d}}_t+\mathbf{v}\cdot\nabla \mathbf{d}
 &=&\Delta \wha{\mathbf{d}}-\mathbf{f}(\mathbf{d})-\p_t\mathbf{d}_E(t),\label{3E}
 \eea
 subject to homogeneous Dirichlet boundary conditions and initial conditions
 \bea
 && \mathbf{v} =\mathbf{0},\quad \wha{\mathbf{d}} = \mathbf{0},\qquad \text{ on } \Gamma\times
 (0,T),
 \label{4E}\\
 &&
 \mathbf{v}|_{t=0}=\mathbf{v}_0 ,\quad
 \wha{\mathbf{d}}|_{t=0}=\mathbf{d}_0-\mathbf{d}_{E0},\qquad \text{ in }\Omega.\label{5E}
 \eea
 In \eqref{1E} we have used the identity
 $\nabla \cdot(\nabla \mathbf{d}\odot\nabla \mathbf{d})=\frac12\nabla \left(|\nabla \mathbf{d}|^2\right)
 + (\nabla \mathbf{d})^{\mathrm{tr}}\Delta \mathbf{d}$ and absorbed the gradient term into the pressure (cf. \cite{LL95}).

 Next, let $\bfd_P$ be defined as in \eqref{LP}. Set
 \begin{equation}
 \wt{\mathbf{d}}=\mathbf{d}-\mathbf{d}_P. \label{dtil}
 \end{equation}
  Then system \eqref{1}--\eqref{3} can be rewritten into the following form for the pair $(\bfv, \wt{\bfd})$:
  \bea
 \mathbf{v}_t+\mathbf{v}\cdot\nabla \mathbf{v}-\Delta \mathbf{v}+\nabla P&=&-(\nabla \mathbf{d})^{\mathrm{tr}} \Delta \mathbf{d},\label{1P}\\
 \nabla \cdot \mathbf{v} &=& 0,\label{2P}\\
 \wt{\mathbf{d}}_t+\mathbf{v}\cdot\nabla \mathbf{d}&=&\Delta \wt{\mathbf{d}}-\mathbf{f}(\mathbf{d}),\label{3P}
 \eea
 subject to homogeneous Dirichlet boundary conditions and initial conditions
 \bea
 &&\mathbf{v} =\mathbf{0},\quad \wt{\mathbf{d}} = \mathbf{0},\qquad \text{ on } \Gamma\times
 (0,T),
 \label{4P}
 \\
 &&
 \mathbf{v}|_{t=0}=\mathbf{v}_0,\quad
 \wt{\mathbf{d}}|_{t=0}=\mathbf{d}_0-\mathbf{d}_{E0},\qquad \text{ in } \Omega.\label{5P}
 \eea

 The following lemma on equivalent norms will be useful in our proofs (cf. e.g., \cite{C09})
 \bl \label{dpe}
 Let $\wha{\bfd}$, $\wt{\bfd}$ be the functions defined in \eqref{dhat} and \eqref{dtil}, respectively.
 The following equivalences between norms hold (provided that they are smooth enough)
 \bea
 &&  \|\wt{\mathbf{d}}\|_{\mathbf{H}^1}
 \approx \|\nabla \wt{\mathbf{d}}\|_{\mathbf{L}^2},
 \quad \|\wha{\mathbf{d}}\|_{\mathbf{H}^1}
 \approx \|\nabla \wha{\mathbf{d}}\|_{\mathbf{L}^2} ,\quad\ \mathrm{in}\ \mathbf{H}_0^1(\Omega),\non\\
 &&  \|\wt{\mathbf{d}}\|_{\mathbf{H}^2}
 \approx \|\Delta \wt{\mathbf{d}}\|_{\mathbf{L}^2},\quad \|\wha{\mathbf{d}}\|_{\mathbf{H}^2}
 \approx \|\Delta \wha{\mathbf{d}}\|_{\mathbf{L}^2},\quad\ \mathrm{in}\ \mathbf{H}_0^1(\Omega)
 \cap \mathbf{H}^2(\Omega),\non\\
 && \|\wt{\mathbf{d}}\|_{\mathbf{H}^3}\approx
 \|\nabla(\Delta\wt{\mathbf{d}})\|_{\mathbf{L}^2}+\|\Delta \wt{\mathbf{d}}\|_{\mathbf{L}^2},\quad
 \|\wha{\mathbf{d}}\|_{\mathbf{H}^3}\approx
 \|\nabla(\Delta\wha{\mathbf{d}})\|_{\mathbf{L}^2}+\|\Delta \wha{\mathbf{d}}\|_{\mathbf{L}^2}
 \quad  \mathrm{in}\ \mathbf{H}_0^1(\Omega)\cap \mathbf{H}^3(\Omega).\non
 \eea
  Besides, we have
 \begin{align}
 &\|\Delta \mathbf{d}\|_{\mathbf{L}^2}\leq \|\Delta \mathbf{d}_P\|_{\mathbf{L}^2}+\|\Delta \wt{\mathbf{d}}-\mathbf{f}(\mathbf{d})\|_{\mathbf{L}^2}
     +\|\mathbf{f}(\mathbf{d})\|_{\mathbf{L}^2},\non\\
 &\|\Delta \mathbf{d}\|_{\mathbf{L}^2}\leq \|\Delta \mathbf{d}_E\|_{\mathbf{L}^2}+\|\Delta \wha{\mathbf{d}}-\mathbf{f}(\mathbf{d})\|_{\mathbf{L}^2}
     +\|\mathbf{f}(\mathbf{d})\|_{\mathbf{L}^2}.\non
 \end{align}
 \el

\subsection{Proof of Theorem \ref{we}: existence of global weak solutions}

Introduce the lifted energy functional (cf. \eqref{OE})
 \be
  \wha{\mathcal{E}}(t)=\frac12\|\mathbf{v}(t)\|_{\mathbf{L}^2}^2+\frac12
 \|\nabla\wha{\mathbf{d}}(t)\|_{\mathbf{L}^2}^2+\int_\Omega F(\mathbf{d}(t)) dx, \quad t\in [0,T],\label{E}
 \ee
 where $\wha{\bfd}$ is given by \eqref{dtil}.  Then we can derive the following \emph{basic energy inequality} for problem  \eqref{1}--\eqref{5}:
  \bl \label{BEL} Let $(\bfv, \bfd)$ be a smooth solution to problem \eqref{1}--\eqref{5} on $[0,T]$.
  The following differential inequality holds
  \begin{align}
 &\frac{d}{dt}\wha{\mathcal{E}}(t)
 + \|\nabla \mathbf{v}(t)\|_{\mathbf{L}^2}^2+ \frac12 \|\Delta \wha{\mathbf{d}}(t)-\mathbf{f}(\mathbf{d}(t))\|_{\mathbf{L}^2}^2\non\\
 &\quad \leq \frac32\|\p_t\mathbf{d}_E(t)\|_{\mathbf{L}^2}^2+\frac14\|\p_t\bfd_E(t)\|_{\mathbf{L}^4}^4
 +\frac{13}{3}\int_\Omega F(\bfd(t))dx+\frac14|\Omega|,\quad \forall\,t\in [0,T]. \label{EN1}
 \end{align}
 \el
\begin{proof}
  Multiplying \eqref{1E} and \eqref{3E} by $\mathbf{v}$ and
  $-\Delta \wha{\mathbf{d}}+\mathbf{f}(\mathbf{d})$, respectively,
  integrating over $\Omega$ and adding the results together, we get
 \begin{align}
 &\frac{d}{dt}\left(\frac12\|\mathbf{v}\|_{\mathbf{L}^2}^2
 +\frac12 \|\nabla\wha{\mathbf{d}}\|_{\mathbf{L}^2}^2+\int_\Omega F(\mathbf{d}) dx\right)
 + \|\nabla \mathbf{v}\|_{\mathbf{L}^2}^2+\|\Delta \wha{\mathbf{d}}-\mathbf{f}(\mathbf{d})\|_{\mathbf{L}^2}^2\non\\
 &\quad = \int_\Omega \p_t\mathbf{d}_E\cdot \Delta \wha{\mathbf{d}} dx. \label{EN2}
 \end{align}
Here, we have used the facts $\int_\Omega (\mathbf{v}\cdot \nabla) \mathbf{v}\cdot \mathbf{v}dx =0, \, \int_\Omega\nabla P \cdot
\mathbf{v} dx= 0,\, \int_\Omega (\mathbf{v}\cdot\nabla) \mathbf{d}\cdot \mathbf{f}(\mathbf{d})dx =0$, due to the incompressibility condition
$\nabla \cdot\mathbf{v}=0$.
 Besides, recalling that $F(\bfd)=\frac14(|\bfd|^2-1)^2$ (as mentioned before, we simply set $\epsilon=1$), and using Young's inequality, we have
 \begin{align}
 &|\bfd|^2\leq \frac23 F(\bfd)+\frac12\quad\text{and}\quad  |\bfd|^4\leq \frac{16}{3}F(\bfd).\non
 \end{align}
 As a consequence, the right-hand side of \eqref{EN2} can be estimated as follows
 \begin{align}
 &\left|\int_\Omega \p_t\mathbf{d}_E\cdot \Delta \wha{\mathbf{d}} dx\right|\non\\
 &\quad \leq  \left|\int_\Omega \p_t\mathbf{d}_E\cdot(\Delta \wha{\mathbf{d}}-\mathbf{f}(\mathbf{d})) dx\right|
 + \left|\int_\Omega \p_t\mathbf{d}_E \cdot \mathbf{f}(\mathbf{d})dx\right|\non\\
 &\quad \leq \|\Delta \wha{\mathbf{d}}-\mathbf{f}(\mathbf{d})\|_{\mathbf{L}^2}
 \|\p_t\mathbf{d}_E\|_{\mathbf{L}^2}+ \int_\Omega\big(|\bfd|^3+|\bfd|\big)|\p_t\mathbf{d}_E| dx \non\\
 &\quad \leq \frac14\|\Delta \wha{\mathbf{d}}-\mathbf{f}(\mathbf{d})\|_{\mathbf{L}^2}^2+\frac32\|\p_t\mathbf{d}_E\|_{\mathbf{L}^2}^2
 +\frac14\|\p_t\bfd_E\|_{\mathbf{L}^4}^4+\frac12\int_\Omega |\bfd|^2dx+\frac34\int_\Omega |\bfd|^4dx\non\\
 &\quad \leq \frac14\|\Delta \wha{\mathbf{d}}-\mathbf{f}(\mathbf{d})\|_{\mathbf{L}^2}^2+\frac32\|\p_t\mathbf{d}_E\|_{\mathbf{L}^2}^2
 +\frac14\|\p_t\bfd_E\|_{\mathbf{L}^4}^4 +\frac{13}{3}\int_\Omega F(\bfd)dx+\frac14|\Omega|,\label{ggg1}
 \end{align}
which easily yields our conclusion \eqref{EN1}. The proof is complete.
\end{proof}

\textbf{Proof of Theorem \ref{we}}. We sketch the proof for the existence of global weak solutions to problem \eqref{1}--\eqref{5}.
 It follows from the arguments in \cite{LL95, B, C09} with some necessary modifications. In particular, here we do not assume the condition \eqref{hyp5} that yields
 the weak maximum principle for $\bfd$.

Let the family $\{\mathbf{u}_i\}_{i=1}^{\infty}$ be a basis of the Hilbert space $\mathbf{V}$ as introduced in the proof of Proposition \ref{Lin}.
We denote by $\mathbf{V}_m=\mathrm{span}\{\mathbf{u}_1,...,\mathbf{u}_m\}$ the finite dimensional subspaces of $\mathbf{V}$  spanned by the first
$m$ basis functions and by $\Pi_m$  the orthogonal projection from $\mathbf{H}$ onto $\mathbf{V}_m$.

\textit{Step 1. Semi-Galerkin approximation}. For every $m\in\mathbb{N}$ and $T\in(0,+\infty)$, we consider the following approximate problem (AP):
\vskip0.2truecm
{\it Determine the vectorial functions $\bfv^m(t,x)=\sum_{i=1}^m g^m_i(t)\mathbf{v}_i(x)$ and $\bfd^m(t,x)$ such that}
\be
\mathrm{(AP)}
 \begin{cases}
  \langle \bfv^m_t, \mathbf{w}\rangle_{\mathbf{V}', \mathbf{V}}
 + \int_\Omega (\bfv^m\cdot\nabla)\bfv^m\cdot \mathbf{w} dx
 + \int_\Omega \nabla \bfv^m : \nabla \mathbf{w} dx \\
  \qquad = \int_\Omega (\nabla\bfd^m\odot\nabla\bfd^m) : \nabla \mathbf{w} dx,\quad \forall\,
 \mathbf{w}\in \mathbf{V}_m,\\
  \bfd^m_t+\bfv^m\cdot \nabla\bfd^m=\Delta \bfd^m-\mathbf{f}(\bfd^m),\quad \text{a.e.\ in}\   \Omega\times (0, T), \\
  \bfd^m=\mathbf{h}(x,t),\quad \text{a.e. on} \  \Gamma\times (0,T),\\
  \bfv^m|_{t=0}=\bfv_{0}^m:=\Pi_m \bfv_0, \quad \bfd^m|_{t=0}=\bfd_0,\quad \text{in}\ \Omega.
  \end{cases}
  \non
 \ee
\bp\label{ppn3}
Let $n=2, 3$. Under the assumption of Theorem \ref{we}, for every $m\in\mathbb{N}$, there is a time $T_m\in(0,T]$ depending on
$\bfv_0$, $\bfd_0$, $m$ and $\Omega$ such that the approximate problem (AP) admits a unique solution $(\bfv^m,\bfd^m)$ on $[0,T_m]$ satisfying
$$\bfv^m\in H^1(0,T_m; \mathbf{V}_m), \quad \bfd^m\in L^\infty(0, T_m; \mathbf{H}^1(\Omega)) \cap L^2(0,T_m; \mathbf{H}^2(\Omega)).$$
\ep
\begin{proof}
Let us fix an arbitrary vector
$$\tilde{\mathbf{v}}^m(t,x)=\sum_{i=1}^m \tilde{g}^m_i(t)\mathbf{u}_i\in C([0,T]; \mathbf{V}_m),$$
where $$\tilde{g}_i^m \in C([0,T]), \quad \tilde{g}_i^m(0)=\int_\Omega \bfv_0\cdot \mathbf{u}_idx,\quad \sup_{t\in[0,T]}\sum_{i=1}^m |\tilde{g}^m_i(t)|^2\leq L,$$
with $L=2+2\|\bfv_0\|_{\mathbf{L}^2}^2$. It is obvious that
\be
 \sup_{t\in[0,T]}\|\tilde{\mathbf{v}}^m(t,x)\|_{\mathbf{L}^2}^2\leq L,\quad
 \sup_{t\in[0,T]}\|\tilde{\mathbf{v}}^m(t,x)\|^2_{\mathbf{L}^\infty}\leq L \max_{1\leq i\leq m}\|\mathbf{u}_i\|_{\mathbf{L}^\infty}^2\leq LC_m .\label{aesw}
\ee

(1) We first consider the following semilinear parabolic problem with convection term for $\bfd^m$, under the given fluid velocity $\tilde{\mathbf{v}}^m$:
\begin{equation}
\begin{cases}
\bfd^m_t+\tilde{\mathbf{v}}^m\cdot \nabla\bfd^m=\Delta \bfd^m-\mathbf{f}(\bfd^m),\quad \text{a.e.\ in}\   \Omega\times (0,T),\\
\bfd^m=\mathbf{h},\qquad  \text{a.e. on} \ \Gamma\times(0,T),\\
\bfd^m|_{t=0}=\bfd_0,\quad \text{for}\ x\in \Omega.
\end{cases}\label{apphi1}
\end{equation}
By a standard fixed point argument \cite{LL95,B,C09} one can prove that  problem \eqref{apphi1} is well-posed on $[0,T]$ and admits a unique weak
solution $\bfd^m\in L^\infty(0, T; \mathbf{H}^1(\Omega)) \cap L^2(0,T; \mathbf{H}^2(\Omega))$.
As in \cite{C09}, it suffices to show that under the assumption \eqref{aesw} there exists a constant $K=K(T,m,L)$ depending on
$\|\bfd_0\|_{\mathbf{H}^1}$, $\Vert \mathbf{h}\Vert_{L^2(0,T;\mathbf{H}^{\frac32}(\Gamma))}$,
  $\|\partial_t\mathbf{h}\|_{ L^4(0,T; \mathbf{H}^{\frac12}(\Gamma))}$, $\Omega$, $T$, $m$ and $L$ such that
\begin{align}
\|\bfd^m\|^2_{L^\infty(0,T; \mathbf{H}^1(\Omega))}+\|\bfd^m\|^2_{L^2(0,T; \mathbf{H}^2(\Omega))}\leq K(T,m,L).\label{KTmL}
\end{align}
 For this purpose,  we consider the lifted system  for $\wha{\bfd}^m=\bfd^m-\bfd_E$, where $\bfd_E$ is the lifting function in \eqref{LE} such that
 \begin{equation}
\begin{cases}
\wha{\bfd}^m_t+\tilde{\mathbf{v}}^m\cdot \nabla\bfd^m=\Delta \wha{\bfd}^m-\mathbf{f}(\bfd^m)-\p_t\bfd_E,\quad \text{a.e.\ in}\ \Omega\times (0,T),\\
\wha{\bfd}^m=\mathbf{0},\qquad  \text{a.e. on} \ \Gamma\times(0,T),\\
\wha{\bfd}^m|_{t=0}=\bfd_0-\bfd_{E0},\quad \text{for}\ x\in \Omega.
\end{cases}
\label{apphi1lift}
\end{equation}
Multiplying the equation for $\wha{\bfd}^m$ by $-\Delta \wha{\bfd}^m+\mathbf{f}(\bfd^m)$ and integrating over $\Omega$, using similar estimate as in
\eqref{ggg1}, we obtain
 \begin{align}
 & \frac{d}{dt}\left(\frac12 \|\nabla\wha{\mathbf{d}}^m\|^2_{\mathbf{L}^2}
 +\int_\Omega F(\mathbf{d}^m) dx\right)
 +\|\Delta \wha{\mathbf{d}}^m-\mathbf{f}(\mathbf{d}^m)\|_{\mathbf{L}^2}^2\non\\
 &\quad =-\int_\Omega (\tilde{\mathbf{v}}^m\cdot \nabla\bfd^m)\cdot (-\Delta \wha{\bfd}^m+\mathbf{f}(\bfd^m)) dx
         + \int_\Omega \p_t\mathbf{d}_E\cdot \Delta \wha{\mathbf{d}}^m dx\non\\
 &\quad \leq \frac12\|\Delta \wha{\mathbf{d}}^m-\mathbf{f}(\mathbf{d}^m)\|_{\mathbf{L}^2}^2
             +2C_m L(\|\nabla \wha{\bfd}^m\|_{\mathbf{L}^2}^2+\|\nabla \bfd_E\|_{\mathbf{L}^2}^2)
             +\frac{13}{3}\int_\Omega F(\bfd)dx\non\\
 &\qquad +\frac32\|\p_t\mathbf{d}_E\|_{\mathbf{L}^2}^2
        +\frac14\|\p_t\bfd_E\|_{\mathbf{L}^4}^4 +\frac14|\Omega|. \label{EN2lift}
 \end{align}
From Gronwall's Lemma and Lemma \ref{Ap1} it follows that
\begin{align}
& \frac12 \|\nabla\wha{\mathbf{d}}^m(t)\|^2_{\mathbf{L}^2}+\int_\Omega F(\mathbf{d}^m(t)) dx\non\\
&\quad \leq \Big(\frac12 \|\nabla(\mathbf{d}_0-\bfd_{E0})\|^2_{\mathbf{L}^2}+\int_\Omega F(\bfd_0) dx\Big)e^{(4C_m L+5)t}\non\\
&\qquad \quad +2(C_m L+1)\int_0^t\Big(\|\nabla \bfd_E\|_{\mathbf{L}^2}^2+\|\p_t\mathbf{d}_E\|_{\mathbf{L}^2}^2+\|\p_t\bfd_E\|_{\mathbf{L}^4}^4
 +|\Omega|\Big) d\tau\non\\
 &\quad \leq C(\|\bfd_0\|_{\mathbf{H}^1})e^{(4C_m L+5)t} +C\int_0^t\Big(\|\mathbf{h}\|_{\mathbf{H}^\frac12(\Gamma)}^2
 +\|\p_t\mathbf{h}\|_{\mathbf{H}^\frac12(\Gamma)}^4\Big) d\tau+Ct\non\\
&\quad \leq K_1(T,m,L), \quad \forall\, t\in [0,T].\non
\end{align}
Then, integrating \eqref{EN2lift} with respective to time, we also have
\begin{align}
\int_0^t \|\Delta \wha{\mathbf{d}}^m-\mathbf{f}(\mathbf{d}^m)\|_{\mathbf{L}^2}^2\leq K_2(T,m,L)\quad \forall\, t\in [0,T].\non
\end{align}
The above estimates together with Lemmas \ref{Ap1}, \ref{dpe} and the Sobolev embedding theorem ($n=2,3$) easily yield the estimate \eqref{KTmL}.

Besides, for the semilinear parabolic equation \eqref{apphi1}, it is easy to prove the continuous dependence on the initial data as well as the
given velocity field $\tilde{\mathbf{v}}^m$. Therefore, the solution operator defined by problem
\eqref{apphi1} $\Phi^m: C([0,T]; \mathbf{V}_m)\to  L^\infty(0, T; \mathbf{H}^1(\Omega)) \cap L^2(0, T; \mathbf{H}^2(\Omega))$
such that $\bfd^m=\Phi^m(\tilde{\mathbf{v}}^m)$ is continuous.

(2) Once the solution $\bfd^m$ to problem \eqref{apphi1} is determined, we turn to look for functions
$\mathbf{v}^m(t,x)=\sum_{i=1}^m g^m_i(t)\mathbf{u}_i$
that satisfy the following system, for $i=1,...,m$,
\be
\begin{cases}
  \langle \bfv^m_t, \mathbf{u}_i\rangle_{\mathbf{V}', \mathbf{V}}
 + \int_\Omega (\bfv^m\cdot\nabla)\bfv^m\cdot \mathbf{u}_i dx
 + \int_\Omega \nabla \bfv^m : \nabla \mathbf{u}_i dx \\
  \qquad = \int_\Omega (\nabla\bfd^m\odot \nabla\bfd^m) : \nabla \mathbf{u}_i dx,\\
  \bfv^m|_{t=0}=\Pi_m \bfv_0,\quad \text{in}\ \Omega,
  \end{cases}\label{ODE}
 \ee
which is equivalent to a nonlinear ordinary differential system for the coefficients $\{g^m_i(t)\}_{i=1}^m$.
It is standard to prove that problem \eqref{ODE} admits a unique local solution on $[0,T_*]$ such that
$\mathbf{v}^m(t,x)=\sum_{i=1}^m g^m_i(t)\mathbf{u}_i\in H^1(0,T_*; \mathbf{V}_m)$, where $T_*\in (0,T]$ may depend on $L$, $\bfd^m$ and $m$.
 Similarly to \cite[Lemma A.7]{B}, we have
 \be
    \sup_{t\in [0,T_*]}\|\bfv^m(t)\|_{\mathbf{L}^2}^2+\int_0^{T_*}\|\nabla \bfv^m\|_{\mathbf{L}^2}^2 dt \leq \|\bfv_0\|_{\mathbf{L}^2}^2 + C_m T_*
    \sup_{t\in[0,T_*]}\|\nabla \bfd^m\|_{\mathbf{L}^2}^4.\label{aesum}
 \ee
 Moreover, it is easy to prove the continuous dependence result for the ODE system \eqref{ODE} on its initial data and the given function $\bfd^m$.
 As a consequence, the solution operator defined by problem
 \eqref{ODE} $\Psi^m: L^\infty(0,T_*; \mathbf{H}^1(\Omega))\cap L^2(0,T_*;\mathbf{H}^2(\Omega)) \to  H^1(0,T_*; \mathbf{V}_m)$ such that
 $\bfv^m=\Psi^m(\bfd^m)$ is continuous.

(3) We see that the mapping
 $$\Psi^m\circ\Phi^m:C([0,T_*]; \mathbf{V}_m)\to H^1(0,T_*; \mathbf{V}_m),\quad \Psi^m\circ\Phi^m(\tilde{\mathbf{v}}^m)=\bfv^m$$
 is continuous, where $\bfv^m$ is the solution to problem \eqref{ODE}.
 The compactness of $H^1(0,T_*; \mathbf{V}_m)$ into $C([0,T_*]; \mathbf{V}_m)$ (because $\mathbf{V}_m$ are actually finite dimensional spaces) implies
 that $\Psi^m\circ\Phi^m$ is a compact operator from $C([0,T_*]; \mathbf{V}_m)$ into itself.
 Due to our choice of $L$ and the estimates \eqref{KTmL}, \eqref{aesum}, it holds
 \be
 \sup_{t\in[0,T_*]}\|\bfv^m(t)\|_{\mathbf{L}^2}^2\leq \frac{L}{2}+T_*C_m L K^2.
 \ee
 Hence, we can take $T_m\in(0,T_*)$ to be sufficiently small such that $\|\bfv^m(t)\|_{\mathbf{L}^2}^2\leq L$ for all $t\in [0,T_m]$.
  Then applying the Schauder's fixed point theorem, we can conclude that there exists at least one fixed point $\bfv^m$
 in the bounded closed convex set
 \bea
 &&\Big\{\bfv^m \in C([0,T_m]; \mathbf{V}_m)\ \mid\ \sup_{t\in[0,T_m]} \|\bfv^m(t)\|_{\mathbf{L}^2}^2\leq L \quad \text{with} \ \bfv^m(0)=\Pi_m \bfv_0. \Big\}\non
 \eea
 such that $\bfv^m\in H^1(0,T_m; \mathbf{V}_m)$ and $\bfd^m\in L^\infty(0, T_m; \mathbf{H}^1(\Omega)) \cap L^2(0,
T_m; \mathbf{H}^2(\Omega))$. Finally, uniqueness of the approximate solution $(\bfv^m,\bfd^m)$ is an easy consequence of the energy method.
 The proof of Proposition \ref{ppn3} is complete.
\end{proof}

\textit{Step 2. Uniform estimates and passage to the limit}.
One can easily verify that the approximate solutions $(\bfv^m, \bfd^m)$ ($m\in \mathbb{N}$) obtained in Proposition \ref{ppn3} by the semi-Galerkin scheme
 satisfy an energy inequality like in Lemma \ref{BEL},
which yields that $\|\bfv^m\|_{L^\infty(0,T_m; \mathbf{L}^2)}$, $\|\bfv^m\|_{L^2(0,T_m; \mathbf{H}^1)}$,
$\|\bfd^m\|_{L^\infty(0,T_m;\mathbf{H}^1)}$, $\|\bfd^m\|_{L^2(0,T_m;\mathbf{H}^2)}$,
 $\|\p_t\bfv^m\|_{L^p(0,T_m; \mathbf{V}')}$, $\|\p_t \bfd^m\|_{L^p(0,T_m; \mathbf{L}^2)}$ ($p=2$ if $n=2$, $p=\frac43$ if $n=3$)
 are uniformly bounded with respect to the parameter $m$ and these bounds only depend on $T$, not on $T_m$.
 Therefore,  all the local approximate solutions $(\bfv^m, \bfd^m)$ on $[0,T_m]$ can be extended up to $[0,T]$.
 Moreover, using the same argument as in \cite{B,LL95}, we can pass to the limit as $m\to+\infty$ and prove the existence of global weak solutions
 to problem \eqref{1}--\eqref{5}. The details are omitted here.

Thus, Theorem \ref{we} is proved.

\subsection{Higher-order estimates for global strong solutions in $2D$}

The following higher-order energy inequality will be helpful to derive \textit{a priori} estimates for global strong solutions to problem \eqref{1}--\eqref{5}:
  \bl \label{H2D}
 Let $n=2$. Assume that the assumptions of Theorem \ref{exe2d} are satisfied.
 If $(\mathbf{v}, \mathbf{d})$ is a smooth solution to problem \eqref{1}--\eqref{5}, then it satisfies the following higher-order differential inequality
  \be
  \frac{d}{dt}\mathcal{A}_P(t)+\mathcal{B}_P(t)\leq C_T\big(\mathcal{A}_P^2(t) + R(t)\big), \label{2dI}
  \ee
  with
  \begin{align}
  \mathcal{A}_P(t)&=\|\nabla \mathbf{v}(t)\|_{\mathbf{L}^2}^2 + \|\Delta \wt{\mathbf{d}}(t)-\mathbf{f}(\mathbf{d}(t))\|_{\mathbf{L}^2}^2,\nonumber\\
  \mathcal{B}_P(t)&=\|S \mathbf{v}(t)\|_{\mathbf{L}^2}^2+ \|\nabla (\Delta  \wt{\mathbf{d}}(t)-\mathbf{f}(\mathbf{d}(t)))\|_{\mathbf{L}^2}^2,\non\\
  R(t)&=\| \mathbf{d}_P(t)\|_{\mathbf{H}^2}^4+\|\mathbf{d}_P(t)\|_{\mathbf{H}^3}^2+1.\label{2dr}
  \end{align}
  Here, $\wt{\mathbf{d}}=\bfd-\bfd_P$ with $\bfd_P$ being the parabolic lifting function (see \eqref{dtil}), and $C_T$ is a positive constant
  depending on $\|\mathbf{v}_0\|_{\mathbf{L}^2}$,
  $\|\mathbf{d}_0\|_{\mathbf{H}^1}$, $\|\mathbf{h}\|_{L^2(0,T;\mathbf{H}^\frac52(\Gamma))}$,  $\|\partial_t\mathbf{h}\|_{ L^4(0,T;\mathbf{H}^{\frac12}(\Gamma))}$,
  $\Omega$ and T.
 \el

 \begin{proof}
 The proof mainly follows from \cite[Lemma 2.6]{GW13}. However, since we do not assume \eqref{hyp5}, the weak maximum principle \eqref{max} no longer holds.
 Therefore, some modifications are necessary in the proof, which will be sketched below.

Taking the time derivative of $\mathcal{A}_P(t)$, using the facts
 $\Delta\wt{\mathbf{d}}-\mathbf{f}(\mathbf{d})|_\Gamma=\mathbf{0}$ and
 $\int_\Omega S \mathbf{v}\cdot \mathbf{v}_t dx=\int_\Omega (-\Delta \mathbf{v})\cdot \mathbf{v}_t dx$, then by a
direct calculation we obtain
 \begin{align}
 & \frac12\frac{d}{dt}\mathcal{A}_P(t)+ (\|S \mathbf{v}\|_{{\mathbf L}^2}^2+ \|\nabla(\Delta \wt{\mathbf{d}}-\mathbf{f}(\mathbf{d})\|_{{\mathbf L}^2}^2)\non\\
 &\quad = -\int_\Omega S \mathbf{v}\cdot (\mathbf{v}\cdot\nabla \mathbf{v}) dx
          -\int_\Omega S\mathbf{v}\cdot [(\nabla \mathbf{d})^{\mathrm{tr}}\Delta \mathbf{d}] dx\non\\
 &\qquad -\int_\Omega [\nabla (\mathbf{v}\cdot\nabla) \mathbf{d} ]:  \nabla (\Delta \wt{\mathbf{d}}-\mathbf{f}(\mathbf{d}))dx
         -\int_\Omega \p_t\mathbf{f}(\mathbf{d})\cdot  (\Delta \wt{\mathbf{d}}-\mathbf{f}(\mathbf{d})) dx \non\\
 &\quad :=\sum_{j=1}^{4}I_j.  \label{2dI1}
 \end{align}
Here we also used the fact $\p_t\mathbf{f}(\mathbf{d})=2(\bfd\cdot\p_t\bfd)\bfd+|\bfd|^2\p_t\bfd-\p_t\bfd:=\mathbf{f}'(\bfd)\p_t\bfd$.

By means of the Sobolev embedding theorem ($n=2$), it is easy to see that (see \cite{LL95,GW13})
 \be
 |I_1|\leq  \varepsilon \|S \mathbf{v}\|^2+C\|\nabla \mathbf{v}\|^2.\non
 \ee
Next, on account of the lower-order estimates for $(\bfv, \bfd)$ in Proposition \ref{lowe} and the facts
$\wt{\bfd}|_\Gamma=(\Delta \wt{\bfd}-\mathbf{f}(\bfd))|_{\Gamma}=0$,  from the Sobolev embedding theorem and Agmon's inequality ($n=2$) we infer
  \bea
   \|\bfd\|_{\mathbf{H}^2}
    &\leq& \|\wt{\bfd}\|_{\mathbf{H}^2}+\|\bfd_P\|_{\mathbf{H}^2}\non\\
    &\leq& C\|\Delta \wt{\bfd}-\mathbf{f}(\bfd)\|_{\mathbf{L}^2}+\|\mathbf{f}(\bfd)\|_{\mathbf{L}^2}+\|\bfd_P\|_{\mathbf{H}^2}\non\\
    &\leq& C\|\Delta \wt{\bfd}-\mathbf{f}(\bfd)\|_{\mathbf{L}^2}+\|\bfd_P\|_{\mathbf{H}^2}+C,\non
  \eea
 \bea
 \|\bfd\|_{\mathbf{H}^3}
    &\leq& \|\wt{\bfd}\|_{\mathbf{H}^3}+\|\bfd_P\|_{\mathbf{H}^3}\non\\
    &\leq& C\|\Delta \wt{\bfd}\|_{\mathbf{H}^1}+\|\bfd_P\|_{\mathbf{H}^3}\non\\
    &\leq& C\|\Delta \wt{\bfd}-\mathbf{f}(\bfd)\|_{\mathbf{H}^1}+C\|\mathbf{f}(\bfd)\|_{\mathbf{H}^1}+\|\bfd_P\|_{\mathbf{H}^3}\non\\
    &\leq& C\|\nabla (\Delta \wt{\bfd}-\mathbf{f}(\bfd))\|_{\mathbf{L}^2}+C(\|\bfd\|_{\mathbf{\infty}}^2+1)\|\bfd\|_{\mathbf{H}^1}+\|\bfd_P\|_{\mathbf{H}^3}\non\\
    &\leq& C\|\nabla (\Delta \wt{\bfd}-\mathbf{f}(\bfd))\|_{\mathbf{L}^2}+C\|\bfd\|_{\mathbf{H}^2}+\|\bfd_P\|_{\mathbf{H}^3}\non\\
    &\leq& C\|\nabla (\Delta \wt{\bfd}-\mathbf{f}(\bfd))\|_{\mathbf{L}^2}+C\|\bfd_P\|_{\mathbf{H}^3}+C.\non
 \eea
As a consequence, we have
  \bea
  \|\nabla \mathbf{d}\|^2_{\mathbf{L}^4}
  &\leq& C\|\mathbf{d}\|_{\mathbf{H}^2}\| \mathbf{d}\|_{\mathbf{H}^1}\non\\
  &\leq& C\|\Delta \wt{\mathbf{d}}-\mathbf{f}(\mathbf{d})\|_{\mathbf{L}^2}+C\|\mathbf{d}_P\|_{\mathbf{H}^2}+C,\non
  \eea
  \bea
  \|\nabla \mathbf{d}\|^2_{\mathbf{L}^\infty}
  &\leq& C\|\nabla \mathbf{d}\|_{\mathbf{H}^2}\|\nabla \mathbf{d}\|_{\mathbf{L}^2}\non\\
  &\leq& C\|\nabla (\Delta \wt{\bfd}-\mathbf{f}(\bfd))\|_{\mathbf{L}^2}+C\|\bfd_P\|_{\mathbf{H}^3}+C,\non
  \eea
  \bea
  \|\Delta
  \wt{\mathbf{d}}-\mathbf{f}(\mathbf{d})\|^2_{\mathbf{L}^4}
  &\leq&
  C\|\nabla (\Delta \wt{\mathbf{d}}-\mathbf{f}(\mathbf{d}))\|_{\mathbf{L}^2}
  \|\Delta \wt{\mathbf{d}}-\mathbf{f}(\mathbf{d})\|_{\mathbf{L}^2}.\non
 \eea
 Using the above inequalities, we obtain the estimates for $I_2$ and $I_3$ such that
 \begin{align}
  |I_2|&\leq \left |\int_\Omega S\mathbf{v}\cdot [(\nabla \mathbf{d})^{\mathrm{tr}} (\Delta \wt{\mathbf{d}}-\mathbf{f}(\mathbf{d}))] dx\right|
            + \left |\int_\Omega S\mathbf{v}\cdot [(\nabla \mathbf{d})^{\mathrm{tr}} \mathbf{f}(\mathbf{d})]dx\right|\non\\
 &\quad  + \left|\int_\Omega S\mathbf{v}\cdot [(\nabla \mathbf{d})^{\mathrm{tr}}\Delta\mathbf{d}_P] dx\right|\non\\
 &\leq \|S\mathbf{v}\|_{\mathbf{L}^2}\|\nabla \mathbf{d}\|_{\mathbf{L}^4} \|\Delta \wt{\mathbf{d}}-\mathbf{f}(\mathbf{d})\|_{\mathbf{L}^4}
       +\|S\mathbf{v}\|_{\mathbf{L}^2}\|\nabla \mathbf{d}\|_{\mathbf{L}^\infty} \|\mathbf{d}_P\|_{\mathbf{H}^2}\non\\
 &\leq \varepsilon \|S\mathbf{v}\|_{\mathbf{L}^2}^2+ C\|\nabla \mathbf{d}\|^2_{\mathbf{L}^4} \|\Delta \wt{\mathbf{d}}-\mathbf{f}(\mathbf{d})\|^2_{\mathbf{L}^4}
       + C\|\nabla \mathbf{d}\|^2_{\mathbf{L}^\infty}\|\mathbf{d}_P\|_{\mathbf{H}^2}^2 \non\\
 &\leq \varepsilon \|S\mathbf{v}\|_{\mathbf{L}^2}^2+
 C\|\nabla (\Delta \wt{\mathbf{d}}-\mathbf{f}(\mathbf{d}))\|_{\mathbf{L}^2} \|\Delta
 \wt{\mathbf{d}}-\mathbf{f}(\mathbf{d})\|_{\mathbf{L}^2}(\|\Delta \wt{\mathbf{d}}-\mathbf{f}(\mathbf{d})\|_{\mathbf{L}^2}+\| \mathbf{d}_P\|_{\mathbf{H}^2}+1)\non\\
 &\quad  + C\|\mathbf{d}_P\|_{\mathbf{H}^2}^2(\|\nabla (\Delta \wt{\mathbf{d}}-\mathbf{f}(\mathbf{d}))\|_{\mathbf{L}^2}+\|\mathbf{d}_P\|_{\mathbf{H}^3}+1) \non\\
 &\leq \varepsilon \|S\mathbf{v}\|_{\mathbf{L}^2}^2+ \varepsilon \|\nabla (\Delta \wt{\mathbf{d}}-\mathbf{f}(\mathbf{d}))\|_{\mathbf{L}^2}^2
 +C\|\Delta \wt{\mathbf{d}}-\mathbf{f}(\mathbf{d})\|_{\mathbf{L}^2}^4\non\\
 &\quad + C\| \mathbf{d}_P\|_{\mathbf{H}^2}^4+C\|\mathbf{d}_P\|_{\mathbf{H}^3}^2+C,\non
 \end{align}
 \begin{align}
 |I_3|&\leq  \|\nabla (\Delta \wt{\mathbf{d}}-\mathbf{f}(\mathbf{d}))\|_{\mathbf{L}^2}(\|\nabla \mathbf{v}\|_{\mathbf{L}^4}
 \|\nabla \mathbf{d}\|_{\mathbf{L}^4}+\|\mathbf{v}\|_{\mathbf{L}^\infty}\|\mathbf{d}\|_{\mathbf{H}^2})\non\\
 &\leq \varepsilon \|\nabla (\Delta \wt{\mathbf{d}}-\mathbf{f}(\mathbf{d}))\|_{\mathbf{L}^2}^2+ C\|\nabla \mathbf{v}\|_{\mathbf{L}^4}^2
 \|\nabla \mathbf{d}\|_{\mathbf{L}^4}^2+ C\|\mathbf{v}\|^2_{\mathbf{L}^\infty}\|\mathbf{d}\|_{\mathbf{H}^2}^2\non\\
 &\leq \varepsilon \|\nabla (\Delta \wt{\mathbf{d}}-\mathbf{f}(\mathbf{d}))\|_{\mathbf{L}^2}^2
 +C\|\Delta \mathbf{v}\|_{\mathbf{L}^2}\|\nabla \mathbf{v}\|_{\mathbf{L}^2}(\|\Delta \wt{\mathbf{d}}-\mathbf{f}(\mathbf{d})\|_{\mathbf{L}^2}
 +\| \mathbf{d}_P\|_{\mathbf{H}^2}+1)\non\\
 &\quad +C\|\Delta \mathbf{v}\|_{\mathbf{L}^2}\|\mathbf{v}\|_{\mathbf{L}^2}(\|\Delta \wt{\mathbf{d}}-\mathbf{f}(\mathbf{d})\|_{\mathbf{L}^2}^2
 +\| \mathbf{d}_P\|_{\mathbf{H}^2}^2+1)\non\\
 &\leq \varepsilon \|S\mathbf{v}\|_{\mathbf{L}^2}^2+\varepsilon \|\nabla (\Delta \wt{\mathbf{d}}-\mathbf{f}(\mathbf{d}))\|_{\mathbf{L}^2}^2
 + C\|\nabla \mathbf{v}\|_{\mathbf{L}^2}^4+ C\|\Delta \wt{\mathbf{d}}-\mathbf{f}(\mathbf{d})\|^4 \non\\
 &\quad +C\| \mathbf{d}_P\|_{\mathbf{H}^2}^4+C.\non
 \end{align}
For the term $I_5$,
  \begin{align}
 I_5&=-\int_\Omega \mathbf{f}'(\mathbf{d})\wt{\mathbf{d}}_t\cdot  (\Delta\wt{\mathbf{d}}-\mathbf{f}(\mathbf{d}))dx
 -\int_\Omega (\mathbf{f}'(\mathbf{d})\p_t\mathbf{d}_P)\cdot (\Delta\wt{\mathbf{d}}-\mathbf{f}(\mathbf{d}))dx\non\\
 &:=I_{5a}+I_{5b},
 \label{split}
 \end{align}
then recalling the equation for $\wt{\bfd}$ (see \eqref{3P}), we have
 \begin{align}
 |I_{5a}| &\leq \|\mathbf{f}'(\mathbf{d})\|_{\mathbf{L}^2}\|\mathbf{v}\|_{\mathbf{L}^4}
 \|\nabla \mathbf{d}\|_{\mathbf{L}^4}\|\Delta \wt{\mathbf{d}}-\mathbf{f}(\mathbf{d})\|_{\mathbf{L}^4}
 +\|\mathbf{f}'(\mathbf{d})\|_{\mathbf{L}^2}\|\Delta \wt{\mathbf{d}}-\mathbf{f}(\mathbf{d})\|_{\mathbf{L}^4}^2\non\\
 &\leq C\|\nabla \mathbf{v}\|_{\mathbf{L}^2}\|\mathbf{v}\|_{\mathbf{L}^2}\|\nabla \mathbf{d}\|_{\mathbf{L}^4}^2
 +C\|\nabla(\Delta \wt{\mathbf{d}}-\mathbf{f}(\mathbf{d}))\|_{\mathbf{L}^2}\|\Delta \wt{\mathbf{d}}-\mathbf{f}(\mathbf{d})\|_{\mathbf{L}^2}\non\\
 &\leq \varepsilon \|\nabla (\Delta \wt{\mathbf{d}}-\mathbf{f}(\mathbf{d}))\|_{\mathbf{L}^2}^2+C\|\nabla \mathbf{v}\|^2
 +C\|\Delta \wt{\mathbf{d}}-\mathbf{f}(\mathbf{d})\|^2+C\|\mathbf{d}_P\|_{\mathbf{H}^2}^2,\non
 \end{align}
 \begin{align}
 |I_{5b}|&\leq \|\mathbf{f}'(\mathbf{d})\|_{\mathbf{L}^2}\|\p_t \mathbf{d}_P\|_{\mathbf{L}^4}
 \|\Delta \wt{\mathbf{d}}-\mathbf{f}(\mathbf{d})\|_{\mathbf{L}^4}\non\\
 &\leq \varepsilon \|\nabla (\Delta \wt{\mathbf{d}}-\mathbf{f}(\mathbf{d}))\|_{\mathbf{L}^2}^2+C\|\Delta \wt{\mathbf{d}}-\mathbf{f}(\mathbf{d})\|^2
 +C\|\p_t \mathbf{d}_P\|_{\mathbf{H}^1}^2.\label{2dI5b}
 \end{align}
Collecting the above estimates and taking $\varepsilon$ to be sufficiently small, we deduce
 \be
 \frac{d}{dt}\mathcal{A}_P(t)+ \|S \mathbf{v}\|^2+ \|\nabla (\Delta
 \wt{\mathbf{d}}-\mathbf{f}(\mathbf{d}))\|^2
  \leq  C\| \mathbf{d}_P\|_{\mathbf{H}^2}^4+C\|\mathbf{d}_P\|_{\mathbf{H}^3}^2+C\|\p_t \mathbf{d}_P\|_{\mathbf{H}^1}^2+C,\non
  \ee
  which together with the cfact $\p_t\bfd_P=\Delta \bfd_P$ easily yields the inequality \eqref{2dI}.
\end{proof}
 Next, recall the following analysis lemma (see \cite[Lemma 6.2.1]{Z04})
 \bl \label{SZ}
 Let $T$ be given with $0<T\leq +\infty$. Suppose that $y(t)$ and $h(t)$ are nonnegative continuous functions defined on
 $[0,T]$ and satisfy the following conditions: $$
 \frac{dy}{dt}\leq c_1 y^2+ c_2 +h,
 $$
 with
 $
 \int_0^T y(t) dt\leq c_3$, $ \int_0^T h(t)dt\leq c_4,
 $
 where $c_i\ (i=1,2,3,4)$ are some nonnegative constants. Then for any $\delta\in (0,T)$, the following estimate holds:
 $$
 y(t+\delta)\leq \left(c_3\delta^{-1}+c_2\delta+c_4\right)e^{c_1c_3},$$
 for all $ t\in[0,T-\delta]$.
 Furthermore, if $T=+\infty$, then
 $\displaystyle\lim_{t\to +\infty} y(t)=0$.
 \el
 \noindent We can prove the following
 \bp\label{esAP}
 Under the assumptions of Theorem \ref{exe2d}, it holds
  \begin{align}
  \mathcal{A}_P(t)+\int_0^t\mathcal{B}_P(\tau)d\tau \leq C_T, \quad \forall\, t\in[0,T],\label{esAPBP}
  \end{align}
  where $C_T>0$ is a constant depending on  $\|\mathbf{v}_0\|_{\mathbf{H}^1}$, $\|\mathbf{d}_0\|_{\mathbf{H}^2}$,
  $\Vert \mathbf{h}\Vert_{L^2(0,T;\mathbf{H}^{\frac32}(\Gamma))}$,
  $\|\partial_t\mathbf{h}\|_{ L^4(0,T; \mathbf{H}^{\frac12}(\Gamma))}$, $\Omega$ and $T$.
 \ep
 \begin{proof}
From Lemmas \ref{Ap1}, \ref{Ap2} and the basic energy inequality derived in Lemma  \ref{BEL}  it follows that
 \begin{align}
  &\int_0^T \mathcal{A}_P(t) dt\leq C_T,\qquad   \int_0^T R(t) dt\leq C_T',\label{inter}
  \end{align}
 where $C_T>0$ is a constant depending on  $\|\mathbf{v}_0\|_{\mathbf{L}^2}$, $\|\mathbf{d}_0\|_{\mathbf{H}^1}$,
 $\Vert \mathbf{h}\Vert_{L^2(0,T;\mathbf{H}^{\frac32}(\Gamma))}$, $\|\partial_t\mathbf{h}\|_{ L^4(0,T; \mathbf{H}^{\frac12}(\Gamma))}$, $\Omega$ and $T$,
 while $C_T'>0$ is a constant depending on $\|\mathbf{d}_0\|_{\mathbf{H}^2}$, $\|\mathbf{h}\|_{L^2(0,T;\mathbf{H}^\frac52(\Gamma))}$,
 $\|\partial_t\mathbf{h}\|_{ L^4(0,T;\mathbf{H}^{\frac12}(\Gamma))}$. Taking $\delta>0$ to be sufficiently small, it is easy to see from the differential inequality \eqref{2dI} that $\mathcal{A}_P(t)$ is bounded on $[0,\delta]$.
 On the other hand, using \eqref{2dI}, \eqref{inter} and applying Lemma \ref{SZ}, we deduce that $\mathcal{A}_P(t)$ is also bounded on $[\delta, T]$.
 Finally, integrating \eqref{2dI} with respect to time, we can conclude that \eqref{esAPBP} holds. The proof is complete.
 \end{proof}
\begin{remark}
From the definitions of $\mathcal{A}_P$, $\mathcal{B}_P$, Lemma \ref{Ap2} and the lower-order estimates obtained in Proposition \ref{lowe}, we can easily prove the conclusion of Theorem \ref{exe2d}.
\end{remark}


\section*{Acknowledgment}
C. Cavaterra and E. Rocca are supported by the FP7-IDEAS-ERC-StG \#256872 (EntroPhase)
and by GNAMPA (Gruppo Nazionale per l'Analisi Matematica, la Probabilit\`a
e le loro Applicazioni) of INdAM (Istituto Nazionale di Alta
Matematica). E. Rocca is supported also by IMATI -- C.N.R. Pavia.
H. Wu is partially supported by NNSFC under the grant No. 11371098, 11631011 and Shanghai Center for Mathematical
Sciences of Fudan University.


\end{document}